\documentclass[12pt]{article}
\usepackage{amsmath,enumerate,amsbsy,amsfonts,amssymb,mathabx,amscd,graphicx,multirow,xcolor,subcaption}
\usepackage{caption}
\usepackage[round,authoryear]{natbib}
\usepackage{authblk}
\usepackage{mathrsfs}
\usepackage[flushleft]{threeparttable}
\RequirePackage[colorlinks,citecolor=blue,urlcolor=blue]{hyperref}
\usepackage{epsfig}

\newtheorem{theorem}{Theorem}

\newtheorem{lemma}{Lemma}

\newtheorem{condition}{Condition}

\usepackage{multirow}
\newcommand{\blind}{1}

\newcommand{\bvarepsilon}{\boldsymbol \varepsilon}

\newcommand{\bOmega}{\boldsymbol \Omega}
\newcommand{\bomega}{\boldsymbol \omega}
\newcommand{\bSigma}{\boldsymbol \Sigma}

\newcommand{\wbSigma}{\widehat {\boldsymbol \Sigma}}

\newcommand{\bphi}{\boldsymbol \phi}

\newcommand{\btheta}{\boldsymbol \theta}

\newcommand{\bbf}{{\mathbf f}}

\newcommand{\bs}{{\mathbf s}}

\newcommand{\bD}{{\bf D}}
\newcommand{\bA}{{\bf A}}
\newcommand{\bB}{{\bf B}}

\newcommand{\bI}{{\bf I}}

\newcommand{\bX}{{\bf X}}

\newcommand{\bZ}{{\bf Z}}

\newcommand{\cD}{{\cal D}}

\newcommand{\cN}{{\cal N}}
\newcommand{\cU}{{\cal U}}
\newcommand{\cC}{{\cal C}}

\newcommand{\cS}{{\cal S}}

\newcommand{\eR}{\mathbb{R}}

\newcommand{\tA}{\text{A}}

\newcommand{\tF}{\text{F}}

\newcommand{\cov}{\text{cov}}
\newcommand{\var}{\text{var}}

\newcommand{\err}{\text{err}}
\newcommand{\argmin}{\ensuremath{\operatornamewithlimits{argmin}}}

\newcommand{\A}{{\scriptscriptstyle \textup{A}}}
\newcommand{\tS}{{\scriptscriptstyle \textup{S}}}

\newcommand{\SC}{{\scriptscriptstyle \textup{SC}}}
\newcommand{\AL}{{\scriptscriptstyle \textup{AL}}}
\newcommand{\U}{{\scriptscriptstyle\textup{U}}}
\newcommand{\ADHD}{{\scriptscriptstyle\textup{ADHD}}}
\newcommand{\TDC}{{\scriptscriptstyle\textup{TDC}}}
\newcommand{\low}{{\scriptstyle\textup{low}}}
\newcommand{\high}{{\scriptstyle\textup{high}}}

\newcommand\ceil[1]{\lceil#1\rceil}

\def\T{{ \mathrm{\scriptscriptstyle T} }}

\begin{document}

\def\spacingset#1{\renewcommand{\baselinestretch}%
{#1}\small\normalsize} \spacingset{1}


	\if1\blind
	{
	\spacingset{1.25}
\title{\bf \Large Adaptive Functional Thresholding for Sparse Covariance Function Estimation in High Dimensions}
\author[1]{Qin Fang}
\author[2]{Shaojun Guo}
\author[1]{Xinghao Qiao}
\affil[1]{Department of Statistics, London School of Economics and Political Science, U.K.}
\affil[2]{Institute of Statistics and Big Data, Renmin University of China, P.R. China}
		\setcounter{Maxaffil}{0}
		
		\renewcommand\Affilfont{\itshape\small}
		\date{\vspace{-5ex}}
		\maketitle
	} \fi
	
	\if0\blind
	{\spacingset{2}
		\bigskip
		\bigskip
		\bigskip
		\begin{center}
			{\Large\bf Adaptive Functional Thresholding for Sparse Covariance Function Estimation in High Dimensions}
		\end{center}
		\medskip
	} \fi

\bigskip
\spacingset{1.5}
\begin{abstract}
Covariance function estimation is a fundamental task in multivariate functional data analysis and arises in many applications. In this paper, we consider estimating sparse covariance functions for high-dimensional functional data, where the number of random functions $p$ is comparable to, or even larger than the sample size $n$.  Aided by the Hilbert--Schmidt norm of functions, we introduce a new class of functional thresholding operators that combine functional versions of thresholding and shrinkage, and propose the adaptive functional thresholding estimator by incorporating the variance effects of individual entries of the sample covariance function into functional thresholding. To handle the practical scenario where curves are partially observed with errors, we also develop a nonparametric smoothing approach to obtain the smoothed adaptive functional thresholding estimator and its binned implementation to accelerate the computation. We investigate the theoretical properties of our proposals when $p$ grows exponentially with $n$ under both fully and partially observed functional scenarios. Finally, we demonstrate that the proposed adaptive functional thresholding estimators significantly outperform the competitors through extensive simulations and the functional connectivity analysis of two neuroimaging datasets.
\end{abstract}

\noindent%
{\it Keywords:} Binning; High-dimensional functional data; 
Functional connectivity;  Functional sparsity; Local linear smoothing; Partially observed functional data.
\vfill

\spacingset{1.7} 
\section{Introduction}
\label{sec.intro}
The covariance function estimation plays an important role in functional data analysis, while existing methods are restricted to data with a single or small number of random functions. Recent advances in technology have made multivariate or even high-dimensional functional datasets increasingly common in various applications: e.g., time-course gene expression data in genomics \cite[]{storey2005}, air pollution data in environmental studies \cite[]{kong2016} and different types of brain imaging data in neuroscience \cite[]{li2018,qiao2019a}.
Under such scenarios, suppose we observe $n$ independent samples $\bX_i(\cdot)=\{X_{i1}(\cdot), \dots, X_{ip}(\cdot)\}^{\T}$ $(i=1, \dots,n)$ defined on a compact interval $\cU$ with covariance function $\bSigma(u,v)=\{\Sigma_{jk}(u,v)\}_{p \times p}=\cov\{\bX_i(u), \bX_i(v)\}$ for $u,v \in \cU,$ which can also be seen as a matrix of marginal- and cross-covariance functions. 
Besides being of interest in itself, an estimator of $\bSigma(\cdot,\cdot)$ is useful for many applications including, e.g.,
dimension reduction via multivariate functional principal components analysis (FPCA) \cite[]{happ2018} or functional factor model \cite[]{guo2022}
or functional independent component analysis, and functional classification \cite[]{park2021}

Our paper focuses on estimating $\bSigma$ under high-dimensional scaling, where $p$ can be comparable to, or even larger than $n.$ 
In this setting, the sample covariance function 
$$
\widehat \bSigma(u,v) =\{\widehat \Sigma_{jk}(u,v)\}_{p \times p} =\frac{1}{n-1} \sum_{i=1}^n \{\bX_i(u)-\widebar \bX(u)\}\{\bX_i(v) - \widebar \bX(v)\}^{\T},~~u,v \in \cU,
$$
where $\widebar \bX(\cdot) = n^{-1} \sum_{i = 1}^n \bX_i(\cdot),$ performs poorly, and some lower-dimensional structural assumptions need to be imposed to estimate $\bSigma(u,v)$ consistently. In contrast to extensive work on estimating high-dimensional sparse covariance matrices \cite[]{bickel2008,rothman2009,cai2011,chen2016,avella2018,wang2020}, research on sparse covariance function estimation in high dimensions remains largely unaddressed in the literature.

In this paper, we consider estimating sparse covariance functions via adaptive functional thresholding. 
To achieve this, we introduce a new class of functional thresholding operators that combine functional versions of thresholding and shrinkage based on the Hilbert-Schmidt norm of functions, and develop an adaptive functional thresholding procedure on $\widehat\bSigma(\cdot,\cdot)$ using entry-dependent functional thresholds that automatically adapt to the variability of $\widehat\Sigma_{jk}(\cdot,\cdot)$'s. To provide theoretical guarantees of our method under high-dimensional scaling, it is essential to develop standardized concentration results taking into account the variability adjustment. Compared with adaptive thresholding for non-functional data \cite[]{cai2011}, the intrinsic infinite-dimensionality of each $X_{ij}(\cdot)$ leads to a substantial rise in the complexity of sparsity modeling and theoretical analysis, as one needs to rely on some functional norm of standardized $\widehat\Sigma_{jk}$'s, e.g., the Hilbert--Schmidt norm, to enforce the functional sparsity in $\widehat\bSigma$ and tackle more technical challenges for standardized processes within an abstract Hilbert space. To handle the practical scenario where functions are partially observed with errors, it is desirable to apply nonparametric smoothers in conjunction with adaptive functional thresholding. This poses a computationally intensive task especially when $p$ is large, thus calling for the development of fast implementation strategy.


There are many applications of the proposed sparse covariance function estimation method in
neuroimaging analysis, where brain signals are measured over time at a large number of regions of interest (ROIs) for individuals. Examples include the brain-computer interface classification \cite[]{lotte2018} and the brain functional connectivity identification \cite[]{rogers2007}.
Traditional neuroimaging analysis models brain signals for each subject as multivariate random variables, where each ROI is represented by a random variable, and hence the covariance/correlation matrices of interest are estimated by treating the time-course data of each ROI as repeated observations.
However, due to the non-stationary and dynamic features of signals \cite[]{chang2010}, the strategy of averaging over time fails to characterize the time-varying structure leading to the loss of information in the original space. To overcome these drawbacks, we follow recent proposals to model signals directly as multivariate random functions with each ROI represented by a random function \cite[]{li2018,qiao2019a,zapata2019,lee2021}. The identified functional sparsity pattern in our estimate of $\bSigma$ can be used to recover the functional connectivity network among different ROIs, which is illustrated using examples of functional magnetic resonance imaging (fMRI) datasets in Section~\ref{sec.real}.

Our paper makes useful contributions at multiple fronts. On the method side, it generalizes the thresholding/sparsity concept in multivariate statistics to the functional setting and offers a novel adaptive functional thresholding proposal to handle the heteroscedastic problem of the sparse covariance function estimation 
motivated from neuroimaging analysis and many statistical applications, e.g., those in Section~\ref{subsec.app}.
It also provides an alternative way of identifying correlation-based functional connectivity with no need to
specify the correlation function, the estimation of which poses challenges as the inverses of $\Sigma_{jj}(u,v)$'s are unbounded. In practice when functions are observed with errors at either a dense grid of points or a small subset of points, we also develop a unified local linear smoothing approach to obtain the smoothed adaptive functional thresholding estimator and its fast implementation via binning \cite[]{fan1994} to speed up the computation without sacrificing the estimation accuracy.
On the theory side, we show that the proposed estimators enjoy the convergence and support recovery properties under both fully and partially observed functional scenarios when $p$ grows exponentially fast relative to $n$. The proof relies on tools from empirical process theory due to the infinite-dimensional nature of functional data and some novel standardized concentration bounds in the Hilbert--Schmidt norm to deal with issues of high-dimensionality and variance adjustment. Our theoretical results and adopted techniques are general, and can be applied to other settings in high-dimensional functional data analysis.

The remainder of this paper is organized as follows. Section~\ref{sec.method} introduces a class of functional thresholding operators, based on which we propose the adaptive functional thresholding of the sample covariance function. We then discuss a couple of applications of the sparse covariance function estimation. Section~\ref{sec.theory} presents convergence and support recovery analysis of our proposed estimator. 
In Section~\ref{sec.partial}, we develop a nonparametric smoothing approach and its binned implementation to deal with partially observed functional data, and then investigate its theoretical properties.
In Sections~\ref{sec.sim} and \ref{sec.real}, we demonstrate the uniform superiority of the adaptive functional thresholding estimators over the universal counterparts through an extensive set of simulation studies and the  functional connectivity analysis of two neuroimaging datasets, respectively. 
All technical proofs are relegated to the Supplementary Material.

\section{Methodology}
\label{sec.method}
\subsection{Functional thresholding}
\label{subsec.fthreshold}
We begin by introducing some notation. 
Let $L_2(\cU)$ denotes a Hilbert space of square integrable functions defined on $\cU$ and $\mathbb S = L_2(\cU) \otimes L_2(\cU),$ where $\otimes$ is the Kronecker product. For any $Q \in {\mathbb S},$ we denote its Hilbert--Schmidt norm by $\|Q\|_{\cS} = \{\int \int Q(u,v)^2\rm{du} \rm{dv}\}^{1/2}.$ 
With the aid of Hilbert--Schmidt norm, for any regularization parameter $\lambda \geq 0,$ we first define a class of functional thresholding operators $s_\lambda:~\mathbb{S} \to \mathbb{S}$ that satisfy the following conditions:
\begin{enumerate}[(i)]
        \item 
        $\|s_\lambda(Z)\|_{\cS} \leq c \|Y\|_\cS$ for all $Z$ and $Y\in \mathbb{S}$ that satisfy $\|Z-Y\|_\cS \leq \lambda$ and some $c>0;$ 
        \item $\|s_\lambda(Z)\|_\cS = 0$ for $\|Z\|_\cS \leq \lambda;$
        \item $\|s_\lambda(Z) - Z\|_\cS \leq \lambda$ for all $Z\in \mathbb{S}.$
\end{enumerate}
Our proposed functional thresholding operators can be viewed as the functional generalization of thresholding operators \cite[]{cai2011}. Instead of a simple pointwise extension of such thresholding operators under functional domain, we advocate a global thresholding rule based on the Hilbert--Schmidt norm of functions that encourages the functional sparsity, in the sense that $s_\lambda(Z)(u,v) = 0$, for all $u,v \in \cU,$ if $\|Z\|_\cS \leq \lambda$ under condition~(ii). Condition~(iii) limits the amount of (global) functional shrinkage in the Hilbert--Schmidt norm to be no more than $\lambda.$


Conditions~(i)--(iii) are satisfied by functional versions of some commonly adopted thresholding rules, which are introduced as solutions to the following penalized quadratic loss problem with various penalties: 
\begin{equation}
\label{func_pls}
    s_{\lambda}(Z) =\underset{\theta \in {\mathbb S}}{\arg\min} \left\{\frac{1}{2} \|\theta - Z\|_{\cS}^2 + p_\lambda (\theta)\right\}
\end{equation}
with $p_{\lambda}(\theta) = \tilde p_{\lambda}(\|\theta\|_{\cS})$ being a penalty function of $\|\theta\|_{\cS}$ to enforce the functional sparsity.

The soft functional thresholding rule results from solving (\ref{func_pls}) with an $\ell_1/\ell_2$ type of penalty, $p_\lambda(\theta)=\lambda \|\theta\|_{\cS},$ and takes the form of 
$
 s_{\lambda}^\tS(Z) = Z(1-\lambda/\|Z\|_{\cS})_{+},   
$
where $(x)_{+}=\max(x,0)$ for $x \in {\mathbb R}.$ This rule can be viewed as a functional generalization of the group lasso solution under the multivariate setting \cite[]{yuan2006}.
To solve (\ref{func_pls}) with an $\ell_0/\ell_2$ type of penalty, $p_\lambda(\theta) = 2^{-1}{\lambda^2}I(\|\theta\|_\cS \neq 0),$ we obtain hard functional threhsolding rule as  
$
Z I (\|Z\|_\cS \ge \lambda),
$
where $I(\cdot)$ is an indicator function.
As a comparison, soft functional thresholding corresponds to the maximum amount of functional shrinkage allowed by condition~(iii), whereas no shrinkage results from hard functional thresholding. Taking the compromise between soft and hard functional thresholding, we next propose functional versions of SCAD \cite[]{fan2001} and adaptive lasso \cite[]{zou2006} thresholding rules.
With a SCAD penalty \cite[]{fan2001} operating on $\|\cdot\|_{\cS}$ instead of $|\cdot|$ for the univariate scalar case, SCAD functional thresholding $s_{\lambda}^\SC(Z)$ is the same as soft functional thresholding if $\|Z\|_{\cS}<2\lambda,$ and equals $Z \{(a-1) - a\lambda/\|Z\|_\cS\}/(a-2)$ for $\|Z\|_\cS \in [2\lambda, a\lambda]$ and $Z$ if $\|Z\|_{\cS} >a\lambda,$ where $a > 2.$ Analogously, adaptive lasso functional thresholding rule is
$s_{\lambda}^\AL(Z) = Z(1-\lambda^{\eta+1}/\|Z\|_{\cS}^{\eta+1})_{+}$ with $\eta \geq 0.$

Our proposed functional generalizations of soft, SCAD and adaptive lasso thresholding rules can be checked to satisfy conditions~(i)--(iii), see Section~\ref{supp.sec_con} of Supplementary Material for details. 
To present a unified theoretical analysis, we focus on functional thresholding operators $s_{\lambda}(Z)$ satisfying conditions~(i)--(iii).
It is worth noting that, although the hard functional thresholding does not satisfy condition~(i), theoretical results in Section~\ref{sec.theory} still hold for hard functional thresholding estimators under similar conditions with corresponding proofs differing slightly.

In general, conditions~(i)--(iii) are satisfied by a number of solutions to (\ref{func_pls}), where the presence of $\|\cdot\|_\cS$ in both the loss and various penalty functions leads to the solutions 
as functions of $\|Z\|_{\cS}.$
Such connection demonstrates the rationale of imposing Hilbert--Schmidt-norm based conditions~(i)--(iii).
For examples of functional data with some local spikes, 
one may suggest another class of functional thresholding operators $\tilde s_\lambda(Z)$  satisfying three supremum-norm based conditions analogous to conditions~(i)--(iii), where, for any $Q \in \mathbb{S},$ we denote its supremum norm by $\|Q\|_{\infty} = \sup_{u,v \in \cU} |Q(u,v)|.$ In this case, $\tilde s_\lambda(Z)$ can not be directly derived as the solution to (\ref{func_pls}) with $p_{\lambda}(\theta) = \tilde p_{\lambda}(\|\theta\|_{\infty}).$ 
However, by substituting $\|\cdot\|_\cS$  in $s_{\lambda}^\tS(Z),s_{\lambda}^\SC(Z)$ and $s_{\lambda}^\AL(Z)$ with $\|\cdot\|_\infty,$ the corresponding supremum-norm based functional thresholding rules can be presented and checked to satisfy three conditions for $\tilde s_\lambda(Z)$ in a similar fashion.
To study theoretical properties analogous to Theorems~\ref{thm_rate} and \ref{thm_supp} in Section~\ref{sec.theory},
the main challenge is to establish concentration bounds on some standardized processes in the supremum norm, where our tools and results in Section~\ref{supp.sec_pr} of Supplementary Material can be applied accordingly.
In this regard, the $\|\cdot\|_\cS$ that we adopt in  $s_\lambda(Z)$ is not necessarily the unique choice, but serves as the building block for the sparse covariance function estimation problem.

\subsection{Estimation}
\label{subsec.est}
We now discuss our estimation procedure based on $s_{\lambda}(Z).$ 
As the variance of $\widehat \Sigma_{jk}(u,v)$ depends on the distribution of $\{X_{ij}(u), X_{ik}(v)\}$ through higher-order moments, which is intrinsically a heteroscedastic problem, it is more desirable to use entry-dependent functional thresholds that automatically takes into account the variability of  $\widehat\Sigma_{jk}$'s. To achieve this, 
define the variance factors 
$\Theta_{jk}(u,v)  = \var\big([X_{ij}(u)-\mathbb{E}\{ X_{ij}(u)\}][X_{ik}(v)- \mathbb{E}\{X_{ik}(v)\}]\big)$ with corresponding estimators
$$\widehat \Theta_{jk}(u,v) =\frac{1}{n}\sum_{i=1}^n\Big[ \big\{X_{ij}(u)-\widebar X_j(u)\big\}\big\{X_{ik}(v)-\widebar X_k(v)\big\} - \widehat \Sigma_{jk}(u,v)\Big]^2,~~j,k=1, \dots, p.$$
Then the adaptive functional thresholding estimator $\widehat \bSigma_{\A} = \{\widehat \Sigma^\A_{jk}(\cdot,\cdot)\}_{p \times p}$ is defined by 
\begin{equation} 
\label{adp_est}
        \widehat \Sigma_{jk}^\A=  \widehat \Theta_{jk}^{1/2} \times s_{\lambda}\left(\frac{\widehat \Sigma_{jk}}{\widehat \Theta_{jk}^{1/2}}\right),
\end{equation}
which uses a single threshold level to functionally threshold standardized
entries, $\widehat \Sigma_{jk}/\widehat \Theta_{jk}^{1/2}$ for all $j,k,$ resulting in entry-dependent functional thresholds for $\widehat\Sigma_{jk}$'s.
The selection of the optimal regularization parameter $\hat \lambda$ is discussed in Section~\ref{sec.sim}.

An alternative approach to estimate $\bSigma$ is the universal functional thresholding estimator $$\widehat \bSigma_{\U} = \{\widehat \Sigma^\U_{jk}(\cdot,\cdot)\}_{p \times p}~~\text{with}~~
    \widehat \Sigma_{jk}^{\U} = s_{\lambda}\big(\widehat \Sigma_{jk}\big),
$$
where a universal threshold level is used for all entries.
In a similar spirit to \cite{rothman2009}, the consistency of $\widehat \bSigma_{\U}$ requires the assumption that marginal-covariance functions are uniformly bounded in nuclear norm, i.e., $\max_{j} \|\Sigma_{jj}\|_{\cN} \leq M,$ where $\|\Sigma_{jj}\|_{\cN} = \int_{\cU}\Sigma_{jj}(u,u) \rm{du}.$ However, intuitively, such universal method does not perform well when nuclear norms vary over a wide range, or even fails when the uniform boundedness assumption is violated. Section~\ref{sec.sim} provides some empirical evidence to support this intuition.

\subsection{Applications}
\label{subsec.app}

Many statistical problems involving multivariate functional data $\{\bX_i(\cdot)\}_{i=1}^n$ require estimating the covariance function $\bSigma.$ Under a high-dimensional regime, the functional sparsity assumption can be imposed on $\bSigma$ to facilitate its consistent sparse estimates.
Here we outline a couple of applications of our proposals for the sparse covariance function estimation.

Our first application is {\it multivariate FPCA} serving as a natural dimension reduction approach for $\bX_i(\cdot).$ With the aid of Karhunen-Lo\'eve expansion for multivariate functional data \cite[]{happ2018},
$\bX_i(\cdot) = \mathbb{E}\{\bX_i(\cdot)\} + \sum_{l=1}^{\infty} \xi_{il} \bphi_l(\cdot),$
where the principal component scores $\xi_{il} = \sum_{j=1}^p \int [X_{ij}(u)-\mathbb{E}\{X_{ij}(u)\}]\phi_{lj}(u)\rm{du}$ and eigenfunctions $\bphi_l(\cdot)=\{\phi_{l1}(\cdot) , \dots, \phi_{lp}(\cdot)\}^{\T}$ are obtained by carrying out an eigenanalysis of $\bSigma.$ When $p$ is large, we can implement our functional thresholding approach to estimate $\bSigma,$ which guarantees the consistencies of estimated eigenpairs and 
hence multivariate FPCA in high dimensions.

Our second application considers another dimension reduction framework via {\it functional factor model} \cite[]{guo2022} in the form of $\bX_i(\cdot)=\bA \bbf_i(\cdot) + \bvarepsilon_i(\cdot),$
where the common components are driven by $r$ 
functional factors $\bbf_i(\cdot)=\{f_{i1}(\cdot), \dots, f_{ir}(\cdot)\}^{\T},$ $\bA \in {\mathbb R}^{p\times r}$ is the factor loading matrix and the idiosyncratic components are $\bvarepsilon_i(\cdot).$ Denote the covariance functions of $\bX_i(\cdot),$ $\bbf_i(\cdot)$ and $\bvarepsilon_i(\cdot)$ by $\bSigma_X,$ $\bSigma_f$ and $\bSigma_{\varepsilon},$ respectively. It follows from $\int\int \bSigma_{X}(u,v) {\rm d} u {\rm d}v = \bA \int\int \bSigma_f(u,v) {\rm d} u {\rm d}v \bA^{\T} + \int\int \bSigma_{\varepsilon}(u,v) {\rm d} u {\rm d}v$
that, under certain identifiable conditions, $\bA$ can be recovered by performing eigenanalysis of $\int\int\bSigma_{X}(u,v) {\rm d} u {\rm d}v.$ To provide a parsimonious model and enhance interpretability for near-zero loadings, we can impose subspace sparsity conditions \cite[]{vu2013} on $\bA$ that results in a functional sparse $\bSigma_X$ and hence the proposed functional thresholding estimators become applicable.

Our third application explores dimension reduction under a {\it functional independent component analysis} framework, which admits the latent segmentation structure $\bX_i(\cdot)=\bB \bZ_i(\cdot)$ under the orthogonality constraint for $\bB=(\bB_1, \dots, \bB_q)\in {\mathbb R}^{p \times p}$ such that the transformed $p$-vector of curves $\bZ_i(\cdot)=\bB^{\T}\bX_i(\cdot)$ can be divided into $q$ ($q \leq p$) uncorrelated groups $\bB_1^{\T}\bX_i(\cdot), \dots, \bB_q^{\T}\bX_i(\cdot).$
It then follows from 
similar arguments in \cite{chang2017} that the columns of $\bB$ can be recovered by a permutation of $p$ eigenvectors of $\int\int\bSigma(u,v) \rm{du} \rm{dv}.$ With the enforced sparsity assumption on $\bB$ when $p$ is large, $\bSigma$ becomes functional sparse and hence our functional thresholding approach can be applied to $\widehat\bSigma$ directly.

The fourth interesting application considers estimating {\it functional graphical models} targeting at identifying the conditional dependence structure among components in $\bX_i(\cdot).$ \cite{qiao2019a} proposed to estimate a block sparse inverse covariance matrix by treating dimensions of $X_{ij}(\cdot)$'s as approaching infinity. However, to deal with truly infinite-dimensional objects, it is desirable to avoid the estimation of the unbounded inverse of $\bSigma$. For Gaussian graphical models, an innovative transformation \cite[]{fan2016} converts the problem of estimating sparse inverse covariance matrix to that of sparse covariance matrix estimation. 
It is interesting to generalize this transformation strategy to the functional domain and hence our sparse covariance function estimation approach can be adopted.


\section{Theoretical properties}
\label{sec.theory}
We begin with some notation. For a random variable $W,$ define 
$
\|W\|_{\psi} = \inf\big\{c >0: \mathbb{E}[\psi(|W|/c)] \le 1\big\},
$
where $\psi: [0,\infty) \to [0,\infty)$ is a nondecreasing, nonzero convex function with $\psi(0) = 0$ and the norm takes the value $\infty$ if no finite $c$ exists for which $\mathbb{E}[\psi(|W|/c)] \le 1.$ Denote $\psi_k(x) = \exp(x^k)-1$ for $k \ge 1$. Let the packing number $D(\epsilon,d)$ be the maximal number of points that can fit in the compact interval $\cU$ while maintaining a distance greater than $\epsilon$ between all points with respect to the semimetric $d$. We refer to Chapter 8 of \cite{Kosorok2008} for further explanations. 
For $\{X_{ij}(u): u\in\cU, i=1, \dots, n, j=1, \dots, p\},$ define the standardized processes by $ Y_{ij}(u) = [X_{ij}(u) - \mathbb{E}\{ X_{ij}(u)\}]/{\sigma_{j}(u)^{1/2}},$ where $\sigma_j(u) = \Sigma_{jj}(u,u)$. 

To present the main theorems, we need the following regularity conditions.
\begin{condition} \label{con_psi2_process}
(i) For each $i$ and $j,$ $Y_{ij}(\cdot)$ is a separable stochastic process with the semimetric $d_j(u,v) = \|Y_{1j}(u)-Y_{1j}(v)\|_{\psi_2}$ for $u,v \in \cU;$ (ii) For some $u_0 \in \cU,$ $\max_{1\leq j \leq p}\|Y_{1j}(u_0)\|_{\psi_2}$ is bounded.  
\end{condition}
\begin{condition} \label{con_finiteDistance}
The packing numbers $D(\epsilon, d_j)$'s satisfy 
$
{\max_{1 \leq j \leq p}} D(\epsilon,d_j) \le C \epsilon^{-r}
$
for some constants $C,r >0$ and $\epsilon \in (0,1].$
\end{condition}

\begin{condition} \label{con_variance_function}
There exists some constant $\tau>0$ s.t.
$
    \min_{j,k} \inf_{u,v\in \cU}\var\big\{Y_{1j}(u)Y_{1k}(v)\big\} \geq \tau. 
$
\end{condition}

\begin{condition} \label{con_rate}
 The pair $(n, p)$ satisfies $\log p /n^{1/4} \rightarrow 0$ as $n$ and $p\to \infty.$ 
 \end{condition}
 

Conditions~\ref{con_psi2_process} and \ref{con_finiteDistance} 
are standard to characterize the modulus of continuity of sub-Gaussian processes $Y_{ij}(\cdot)$'s, see Chapter 8 of \cite{Kosorok2008}. These conditions also imply that there exist some positive constants $C_0$ and $\eta$ such that 
$
\mathbb{E}[\exp( t\|Y_{1j}\|^2)] \le C_0
$
for all $|t| \le \eta$ and $j$
with $\|Y_{1j}\| = \{\int_{\cU} Y_{1j}(u)^2\rm{du}\}^{1/2},$  
which plays a crucial role in our proof when applying concentration inequalities within Hilbert space. 
Condition~\ref{con_variance_function} restricts the variances of $Y_{ij}(u)Y_{ik}(v)$'s to be uniformly bounded away from zero so that they can be well estimated. It also facilitates the development of some standardized concentration results. 
This condition precludes the case of a Brownian motion $X_{ij}(\cdot)$ starting at 0 for some $j$. However, replacing $X_{ij}(\cdot)$ with a contaminated process $X_{ij}(\cdot) + \xi_{ij},$ 
where $\xi_{ij}$'s are independent from a normal distribution with zero mean and a small variance and are independent of $X_{ij}(\cdot)$'s, Condition~\ref{con_variance_function} is fulfilled while the cross-covariance structure in $\bSigma$ remains the same
in the sense of  $\cov\{X_{ij}(u) + \xi_{ij}, X_{ik}(v)\}=\cov\{X_{ij}(u), X_{ik}(v)\}$ for $k \neq j$ and $u,v \in \cU.$
Condition~\ref{con_rate} allows the high-dimensional case, where $p$ can diverge at some exponential rate as $n$ increases.

We next establish the convergence rate of the adaptive functional thresholding estimator $\widehat \bSigma_\A$ over a large class of ``approximately sparse" covariance functions defined by
\begin{eqnarray*}
    \cC(q,s_0(p),\epsilon_0; \cU) &=& \Big\{\bSigma: \bSigma \succeq 0, \max_{1 \leq j \leq p}\sum_{k=1}^p  \|{\sigma_{j}}\|_\infty^{(1-q)/2} \|{\sigma_{k}}\|_\infty ^{(1-q)/2}\|\Sigma_{jk}\|_\cS^q \leq s_0(p) ,\Big.  \\
    && \hskip 1cm  \Big. \max_{j} \|\sigma_j^{-1}\|_\infty\|\sigma_j\|_\infty \le \epsilon_0^{-1}<\infty \Big\}
\end{eqnarray*}
for some $0 \le q <  1,$ where $\|\sigma_j\|_{\infty} = \sup_{u \in \cU} \sigma_j(u)$ and $\bSigma \succeq 0$ means that $\bSigma=\{\Sigma_{jk}(\cdot,\cdot)\}_{p \times p}$ is positive semidefinite, i.e., 
$\sum_{j,k}\int\int\Sigma_{jk}(u,v)a_j(u)a_k(v)\rm{du} \rm{dv} \geq 0$ for any  $a_j(\cdot) \in L^2(\cU)$ and $j=1, \dots, p.$ See \cite{cai2011} for a similar class of covariance matrices for non-functional data. Compared with the class $$\cC^*(q, s_0(p), M;\cU)=\big\{\bSigma: \bSigma \succeq 0, \max_j \|\sigma_j\|_{\cN} \leq M, \max_{j}\sum_{k=1}^p \|\Sigma_{jk}\|_\cS^q \\ \leq s_0(p)\big\},$$ over which the universal functional thresholding estimator $\wbSigma_\U$ can be shown to be consistent, the columns of a covariance function in $\cC(q, s_0(p),\epsilon_0;\cU)$ are required to be within a weighted $\ell_q/\ell_2$ ball instead of a standard $\ell_q/\ell_2$ ball, where the weights are determined by $\|\sigma_j\|_{\infty}$'s. Unlike $\cC^*(q, s_0(p), M;\cU),$ $\cC(q, s_0(p),\epsilon_0;\cU)$ no longer requires the uniform boundedness assumption on $\|\sigma_{j}\|_{\cN}$'s and allows $\max_j \|\sigma_{j}\|_{\cN}\rightarrow \infty.$ In the special case $q =0$, $\cC(q,s_0(p),\epsilon_0; \cU)$ corresponds to a class of truly sparse covariance functions. Notably, $s_0(p)$ can depend on $p$ and be regarded implicitly as the restriction on functional sparsity.

\begin{theorem}\label{thm_rate}
Suppose that Conditions~\ref{con_psi2_process}-\ref{con_rate} hold. Then there exists some constant $\delta >0$ such that, uniformly on $\cC(q, s_0(p),\epsilon_0; \cU),$ if $\lambda = \delta ({{\log p}/{n}})^{1/2},$ 
\begin{equation} \label{rate_est}
    \|\widehat \bSigma_\A - \bSigma\|_1  = \max_{1 \leq k \leq p}\sum_{j=1}^p\|\widehat\Sigma^{\A}_{jk} -\Sigma_{jk}\|_{\cS}  = O_P\left\{s_0(p) \Big(\frac{\log p}{n}\right)^{\frac{1-q}{2}}\Big\}.
\end{equation}
\end{theorem}

Theorem~\ref{thm_rate} presents the convergence result in the functional version of matrix $\ell_1$ norm. The rate in (\ref{rate_est}) is consistent to those of sparse covariance matrix estimates in \cite{rothman2009,cai2011}.

We finally turn to investigate the support recovery consistency of $\widehat\bSigma_\A$ over the parameter space of truly sparse covariance functions defined by 
$$
\cC_0( s_0(p); \cU) = \Big\{\bSigma: \bSigma \succeq 0, \max_{1 \le j \le p}\sum_{k=1}^p  I (\|\Sigma_{jk}\|_\cS \neq 0) \leq s_0(p)\Big\},
$$
which assumes that $\{\Sigma_{jk}(\cdot, \cdot)\}_{p \times p}$ has at most $s_0(p)$ non-zero entries on each row. The following theorem shows that, with the choice of $\lambda = \delta ({{\log p}/{n}})^{1/2}$ for some constant $\delta>0,$ $\widehat \bSigma_\A$ exactly recovers the support of $\bSigma,$ $\text{supp}(\bSigma) = \{(j,k): \|\Sigma_{jk}\|_{\cS} \neq 0\},$ with probability approaching one. 

\begin{theorem}
\label{thm_supp}
Suppose that Conditions~\ref{con_psi2_process}-\ref{con_rate} hold and 
$
\big\|\Sigma_{jk}/{\Theta_{jk}^{1/2}}\big\|_\cS > (2\delta+\gamma) (\log p/n)^{1/2}
$
for all $(j,k)\in\text{supp}(\bSigma)$ and some $\gamma>0,$
where $\delta$ is stated in Theorem~\ref{thm_rate}.  Then we have that
$$\inf_{\Sigma \in \cC_0}P\big\{\text{supp}(\widehat \bSigma_\A) = \text{supp}(\bSigma) \big \} \to 1 \text{  as  }n \to \infty.$$
\end{theorem}

Theorem~\ref{thm_supp} ensures that $\widehat \bSigma_\A$ achieves the exact recovery of functional sparsity structure in $\bSigma,$ i.e., the graph support in functional connectivity analysis, with probability tending to 1. This theorem holds 
under the condition that the Hilbert-Schmidt norms of non-zero standardized functional entries exceed a certain threshold, which ensures that non-zero components are correctly retained. See an analogous minimum signal strength condition for sparse covariance matrices in \cite{cai2011}.


\section{Partially observed functional data}
\label{sec.partial}

In this section we consider a practical scenario where each $X_{ij}(\cdot)$ is partially observed, with errors, at random measurement locations $U_{ij1}, \dots, U_{ijL_{ij}}\in \cU.$ Let $Z_{ijl}$ be the observed value of $X_{ij}(U_{ijl}).$ Then 
\begin{equation}
\label{model.partial}
Z_{ijl} = X_{ij}(U_{ijl}) + \varepsilon_{ijl},~~l = 1,\dots, L_{ij},
\end{equation}
where $\varepsilon_{ijl}$'s are i.i.d. errors with $\mathbb{E}(\varepsilon_{ijl})=0$ and $\var(\varepsilon_{ijl})=\sigma^2,$ independent of $X_{ij}(\cdot).$ For dense measurement designs all $L_{ij}$'s are larger than some order of $n,$ while for sparse designs all $L_{ij}$'s are bounded \citep{Zhang2016,qiao2020}.

\subsection{Estimation procedure}
\label{subsec.est2}
Based on the observed data, $\{(U_{ijl},Z_{ijl})\}_{1 \leq i \leq  n, 1 \leq j \leq p, 1 \leq  l \leq L_{ij}},$ we next present a unified estimation procedure that handles both densely and sparsely sampled functional data.

We first develop a nonparametric smoothing approach to estimate $\Sigma_{jk}(u,v)$'s. Without loss of generality, we assume that $\bX_i(\cdot)$ has been centered to have mean zero. 
Denote $K_h(\cdot) = h^{-1} K(\cdot/h)$ for a  univariate kernel function $K$ with a bandwidth $h>0.$ A local linear surface smoother (LLS) is employed to estimate cross-covariance functions $\Sigma_{jk}(u,v)$ ($j\neq k$) by minimizing
\begin{equation} \label{eq.kernel.cross.cov}
    \sum_{i = 1}^n \sum_{l = 1}^{L_{ij}} \sum_{m = 1}^{L_{ik}} \Big\{Z_{ijl}Z_{ikm} - \alpha_0-\alpha_1(U_{ijl}-u) -\alpha_2(U_{ikm}-v)\Big\}^2 K_{h_C}({U_{ijl}-u}) K_{h_C}({U_{ikm}-v}),
\end{equation}
with respect to $(\alpha_0, \alpha_1,\alpha_2).$ Let the minimizer of (\ref{eq.kernel.cross.cov}) be $(\hat\alpha_0, \hat\alpha_1,\hat\alpha_2)$ and the resulting estimator is 
$\widetilde \Sigma_{jk}(u,v) = \hat \alpha_0.$ To estimate marginal-covariance functions $\Sigma_{jj}(u,v)$'s, we observe that 
$\cov(Z_{ijl},Z_{ijm})= \Sigma_{jj}(U_{ijl}, U_{ijm}) + \sigma^2 I(l=m),$ and hence apply a LLS to the off-diagonals of the raw covariances $(Z_{ijl}Z_{ijm})_{1 \leq l \leq m \leq L_{ij}}.$ We consider minimizing
\begin{equation}\nonumber
    \sum_{i = 1}^n \sum_{1\leq l \neq m \leq L_{ij} } \Big\{Z_{ijl}Z_{ijm} - \beta_0-\beta_1(U_{ijl}-u) -\beta_2(U_{ikm}-v)\Big\}^2K_{h_M}({U_{ijl}-u}) K_{h_M}({U_{ikm}-v})
\end{equation}
with respect to $(\beta_0, \beta_1,\beta_2),$ thus obtaining the estimate $ \widetilde\Sigma_{jj}(u,v) = \hat \beta_0.$ 
Note that we drop subscripts $j,k$ of $h_{C,jk}$ and $j$ of $h_{M,j}$ to simplify our notation in this section. However, we select different bandwidths $h_{C,jk}$ and $h_{M,j}$ across $j,k=1, \dots, p$ in our empirical studies.  

To construct the corresponding adaptive functional thresholding estimator, a standard approach is to incorporate the variance effect of each $\widetilde \Sigma_{jk}(u,v)$ into functional thresholding. However, the estimation of $\var\{\widetilde \Sigma_{jk}(u,v)\}$'s involves estimating multiple complicated fourth moment terms \cite[]{Zhang2016}, which results in high computational burden especially for large $p.$
Since our focus is on characterizing the main variability of $\widetilde \Sigma_{jk}(u,v)$ rather than estimating its variance precisely, 
we next develop a computationally simple yet effective approach to estimate the main terms in the asymptotic variance of $\widetilde \Sigma_{jk}(u,v).$ 
For $a,b=0, 1, 2,$ let 
\begin{equation}
    \label{T.est}
T_{ab,ijk}(u,v) =  \sum_{l = 1}^{L_{ij}} \sum_{m = 1}^{L_{ik}} g_{ab}\{h_C,(u,v),(U_{ijl},U_{ikm})\}Z_{ijl}Z_{ikm},
\end{equation}
 where $ g_{ab}\big\{h,(u,v),(U_{ijl},U_{ikm})\big\} =   K_{h}({U_{ijl}-u}) K_{h}({U_{ikm}-v})(U_{ijl}-u)^a(U_{ikm}-v)^b.$
According to Section~\ref{supp.post} of Supplementary Material, minimizing (\ref{eq.kernel.cross.cov}) yields the resulting estimator
\begin{equation}\label{eq.kernel.cross.cov.solution}
     \widetilde \Sigma_{jk} = \sum_{i=1}^n\big(W_{1,jk} T_{00,ijk}+ W_{2,jk} T_{10,ijk}+W_{3,jk}T_{01,ijk}\big),
\end{equation}
where $W_{1,jk}, W_{2,jk}, W_{3,jk}$ can be represented via (\ref{eq.kernel.w}) in terms of 
\begin{equation}
\label{S.est}
S_{ab,jk}(u,v) =  \sum_{i = 1}^n \sum_{l = 1}^{L_{ij}} \sum_{m = 1}^{L_{ik}} g_{ab}\big\{h_C,(u,v),(U_{ijl},U_{ikm})\big\},~~a,b=0,1,2.   \end{equation}
It is notable that the estimator $\widetilde \Sigma_{jk}$ in (\ref{eq.kernel.cross.cov.solution}) 
is expressed as the sum of $n$ independent terms. Ignoring the cross-covariances among observations within the subject that are dominated by the corresponding variances, we propose a surrogate estimator for the asymptotic variance of $\widetilde \Sigma_{jk}$ by 
\begin{equation}\label{eq.kernel.var}
    \widetilde \Psi_{jk}   =  I_{jk}\sum_{i=1}^n\big( W_{1,jk}V_{00,ijk} +  W_{2,jk}V_{10,ijk}+  W_{3,jk}V_{01,ijk}\big)^2,
\end{equation}
where 
\begin{equation}
\label{In}
I_{jk} = \Big(\sum_{i=1}^n L_{ij} L_{ik}\Big)^{2}\Big\{ \sum_{i=1}^n \big (L_{ij} L_{ik} h_C^{-2} +L_{ij}^2 L_{ik} h_C^{-1} + L_{ij} L_{ik}^2 h_C^{-1} +L_{ij}^2 L_{ik}^2 \big)\Big\}^{-1}
\end{equation}
and
\begin{equation}
    \label{V.est}
V_{ab,ijk}(u,v)  = \sum_{l = 1}^{L_{ij}} \sum_{m = 1}^{L_{ik}}g_{ab}\big\{h_C,(u,v),(U_{ijl},U_{ikm})\big\}\big\{Z_{ijl}Z_{ikm} - \widetilde \Sigma_{jk}(u,v)\big\}. 
\end{equation}
The rationale of multiplying the rate $I_{jk}$ in (\ref{eq.kernel.var}) is to ensure that $ \widetilde \Psi_{jk}(u,v)$ converges to some finite function when $n \to \infty$ and $h_C\to 0$ as justified in Section~\ref{supp.verify.In} of Supplementary Material.  In particular, the rate $I_{jk}$ can be simplified to $\sum_{i=1}^n L_{ij} L_{ik} h_C^2$ for the sparse or moderately dense case and to $(\sum_{i=1}^n L_{ij} L_{ik})^{2}(\sum_{i=1}^n L_{ij}^2L_{ik}^2)^{-1}$ for the very dense case. Note that $I_{jk}$ is imposed in (\ref{eq.kernel.var}) mainly for the theoretical purpose and hence will not place a practical constraint on our method. 


In a similar procedure as above, the estimated variance factor $\widetilde \Psi_{jj}$ of $\widetilde \Sigma_{jj}$ for each $j$ can be obtained by operating on $\{Z_{ijl}Z_{ijm}\}_{1 \leq i\leq n, 1 \leq l \neq m \leq L_{ij}}$ instead of $\{Z_{ijl}Z_{ikm}\}_{1 \leq i\leq n, 1 \leq l\leq L_{ij}, 1 \leq m \leq L_{ik}}$ for $j\neq k.$ 
Substituting $\widehat \Theta_{jk}$ in (\ref{adp_est}) by $\widetilde \Psi_{jk},$ we obtain the smoothed adaptive functional thresholding estimator
\begin{equation} 
\label{adp_est2}
\widetilde \bSigma_{\A} = (\widetilde \Sigma^\A_{jk})_{p \times p} ~~\text{with}~~
        \widetilde \Sigma_{jk}^\A= \widetilde \Psi_{jk}  ^{1/2} \times s_{\lambda}\left(\frac{\widetilde \Sigma_{jk}}{\widetilde \Psi_{jk}  ^{1/2}}\right).
\end{equation}
For comparison, we also define the smoothed universal functional thresholding estimator as $\widetilde \bSigma_{\U} = (\widetilde \Sigma^\U_{jk})_{p \times p}$ with
    $\widetilde \Sigma_{jk}^{\U} = s_{\lambda}\big(\widetilde \Sigma_{jk}\big).
$


A natural alternative to the proposed LLS-based smoothing procedure considers pre-smoothing each individual data. For densely sampled functional data, the observations $Z_{ij1}, \dots, Z_{ijL_{ij}}$ for each $i$ and $j$ can be pre-smoothed through the local linear smoother to eliminate the contaminated noise, thus producing reconstructed random curves $\widehat X_{ij}(\cdot)$'s before subsequent analysis \cite[]{Zhang2007}. See detailed implementation of pre-smoothing 
in Section~\ref{supp.pre} of Supplementary Material. For sparsely sampled functional data, such pre-smoothing step is not 
viable, while our smoothing proposal builds strength across functions by incorporating information from all the observations, and hence is still applicable. See also Section~\ref{sec.sim.partial} for the numerical comparison between pre-smoothing and our smoothing approach under different measurement designs.

\subsection{Theoretical properties}
\label{subsec.theory2}
In this section, we investigate the theoretical properties of $\widetilde\bSigma_{\A}$ for partially observed functional data. We begin by introducing some notation. For two positive sequences $\{a_n\}$ and $\{b_n\},$ we write $a_n \lesssim b_n$ if there exits a positive constant $c_0$ such that $a_n/b_n \leq c_0.$ We write $a_n \asymp b_n$ if and only if $a_n \lesssim b_n$ and $b_n \lesssim a_n$ hold simultaneously. Before presenting the theory, we impose the following regularity conditions.

\begin{condition}
	\label{cond_T_bd}
	(i) Let $\big\{U_{ijl}: i=1, \dots, n, j \in 1, \dots,p, l=1,  \dots, L_{ij}\big\}$ be i.i.d. copies of a random variable $U$ with density $f_U(\cdot)$ defined on the compact set $\cU,$ with the $L_{ij}$'s fixed. There exist some constants $m_f$ and $M_f$ such that $0<m_f \leq {\inf}_{\cU}f_U(u)\leq {\sup}_{\cU}f_U(u) \leq M_f<\infty;$ (ii) $X_{ij},$ $\varepsilon_{ijl}$ and $U_{ijl}$ are independent for each $i,j,l.$
\end{condition}

\begin{condition}
\label{cond_design}
(i) Under the sparse measurement design, $L_{ij} \le L_0 < \infty$ for all $i, j$ and, under the dense design, $L_{ij} = L \to \infty$ 
as $n \to \infty$ with $U_{ijl}$'s independent of $i;$
(ii) The bandwidth parameters $h_C \asymp h_M \asymp h \to 0$ as $n \to \infty.$ 
\end{condition}

Condition~\ref{cond_T_bd} is standard in functional data analysis literature \cite[]{Zhang2016}. Condition~\ref{cond_design} (i) treats the number of measurement locations $L_{ij}$ as bounded and diverging under sparse and dense measurement designs, respectively. To simplify notation, we assume that $L_{ij}=L$ for the dense case and $h_C$ is of the same order as $h_M$ in Condition~\ref{cond_design} (ii).

\begin{condition}
\label{cond_sparse_rate1}
There exists some constant $\gamma_1 \in (0, 1/2]$ such that
\begin{equation}
\label{coneq_L2_sparse1}
\max_{1 \le j,k \le p }\Big\|\widetilde{\Sigma}_{jk} - \Sigma_{jk}\Big\|_{\cS} \lesssim \sqrt{\frac{\log p}{n^{2\gamma_1}}} +h^2
\end{equation}
with probability approaching one.
\end{condition}

\begin{condition}
\label{cond_sparse_rate2}
	There exist some positive constants $c_1,$ $\gamma_2 \in (0,1/2]$ and some deterministic functions $\Psi_{jk}(u,v)$'s with $\min_{j,k}\inf_{u,v \in \cU}\Psi_{jk}(u,v) \ge c_1$ such that
	\begin{equation}
\label{coneq_uniform_sparse2}
\max_{1 \le j,k \le p }\underset{u,v \in \cU}{\sup}\Big|\widetilde{\Psi}_{jk}(u,v) - \Psi_{jk}(u,v)\Big| \lesssim \sqrt{\frac{\log p}{n^{2\gamma_2}}} + h^2
\end{equation}
with probability approaching one.
\end{condition}
\begin{condition}
\label{cond_sparse_p_rate}
The pair $(n,p)$ satisfies $\log p/n^{\min(\gamma_1,\gamma_2)} \to 0$ and $\log p \ge c_2 n^{2\gamma_1} h^4$ for some positive constant $c_2$ as $n$ and $p \to \infty.$
\end{condition}

We follow \cite{qiao2020} to impose Condition~\ref{cond_sparse_rate1}, in which the parameter $\gamma_1$ depends on $h$ and possibly $L$ under the dense design. This condition is satisfied if there exist some positive constants $c_3,c_4,c_5$ such that for each $j,k = 1, \dots, p$ and $t \in (0,1],$ 
\begin{equation}
\label{coneq.L2}
P\big(\|\widetilde{\Sigma}_{jk} - \Sigma_{jk}\|_{\cS} \ge t + c_5h^2\big) \leq c_4 \exp( - {c_3n^{2\gamma_1}t^2}).    
\end{equation}
The presence of $h^2$ comes from the standard results for bias terms under the boundedness condition for the second-order partial derivatives of $\Sigma_{jk}(u,v)$ over $\cU^2$ \cite[]{yao2005a,Zhang2016}. This concentration result is fulfilled under different measurement schedules ranging from sparse to dense designs as $\gamma_1$ increases. For sparsely sampled functional data, Lemma~4 of \cite{qiao2020} established $L_2$ concentration inequality for $\widetilde \Sigma_{jk}$ for $j=k,$ which not only results in the same $L_2$ rate as that in the sparse case \cite[]{Zhang2016} but also ensures (\ref{coneq.L2}) with the choice of $\gamma_1=1/2-a$ and $h \asymp n^{-a}$ for some positive constant $a<1/2.$ Following the same proof procedure, the same concentration inequality also applies for $j \neq k$ and hence Condition~\ref{cond_sparse_rate1} is satisfied. This condition is also satisfied by densely sampled functional data, since it follows from Lemma~5 of \cite{qiao2020} that (\ref{coneq.L2}) holds for $j=k$ and, with more efforts, also for $j \neq k$ by choosing $\gamma_1=\min(1/2, 1/3+b/6-\epsilon'/2-2a/3)$ for some small constant $\epsilon'>0$ when $h \asymp n^{-a}$ and $L \asymp n^{b}$ for some constants $a, b>0.$ As $L$ grows sufficiently large, $\gamma_1=1/2,$ thus leading to the same rate as that in the ultra-dense case \cite[]{Zhang2016}. 
Condition~\ref{cond_sparse_rate2} gives the uniform convergence rate for $\widetilde \Psi_{jk}(u,v)$ in the same form as (\ref{coneq_L2_sparse1}) but with different parameter $\gamma_2.$ A denser measurement design corresponds to a larger value of $\gamma_2$ and a faster rate in (\ref{coneq_uniform_sparse2}).
See the heuristic verification of Condition~\ref{cond_sparse_rate2} in Section~\ref{supp.verify.In} of Supplementary Material.
Condition~\ref{cond_sparse_p_rate} indicates that $p$ can grow exponentially fast relative to $n.$

We next present the convergence rate of the smoothed adaptive functional thresholding estimator $\widetilde{\bSigma}_{\tA}$ over a class of ``approximate sparse" covariance functions defined by 
\begin{eqnarray*}
    \widetilde{\cC}(q,\tilde s_0(p), \epsilon_0; \cU) &=& \Big\{\bSigma: \bSigma \succeq 0, \max_{1 \leq j \leq p}\sum_{k=1}^p  \|\Psi_{jk}\|_\infty^{(1-q)/2} \left\|\Sigma_{jk}\right\|_\cS^q \leq \tilde{s}_0(p), \\
     && \hskip 1cm  \Big. \max_{j,k} \|\Psi_{jk}^{-1}\|_\infty\|\Psi_{jk}\|_\infty \le \epsilon_0^{-1}<\infty\Big\},
\end{eqnarray*}
for some $0 \le q <1.$ 

\begin{theorem}
\label{thm_sparse_2}
Suppose that Conditions~\ref{cond_T_bd}--\ref{cond_sparse_p_rate} hold. Then there exists some constants $\tilde \delta >0$ such that, uniformly on $\widetilde{\cC}(q, \tilde s_0(p),\epsilon_0; \cU),$ if $\lambda = \tilde \delta (\log p/n^{2\gamma_1})^{1/2},$
\begin{equation}
\label{sparse_rate_est}
    \|\widetilde \bSigma_\A - \bSigma\|_1  = \max_{1 \leq k \leq p}\sum_{j=1}^p\|\widetilde\Sigma^{\A}_{jk} -\Sigma_{jk}\|_{\cS}  = O_P\left\{\tilde{s}_0(p) \Big(\frac{\log p}{n^{2\gamma_1}}\Big)^{\frac{1-q}{2}}\right\}.
\end{equation}
\end{theorem}

The convergence rate of $\widetilde \bSigma_{\A}$ in (\ref{sparse_rate_est}) is governed by internal parameters $(\gamma_1,q)$ and other dimensionality parameters. Larger values of $\gamma_1$ correspond to a more frequent measurement schedule with larger $L$ and result in a faster rate. The convergence result implicitly reveals interesting phase transition phenomena depending on the relative order of $L$ to $n.$ As $L$ grows fast enough, $\gamma_1=1/2$ and the rate is consistent to that for fully observed functional data in (\ref{rate_est}), presenting that the theory for very densely sampled functional data falls in the parametric paradigm. As $L$ grows moderately fast, $\gamma_1<1/2$ and the rate is faster than that for sparsely sampled functional data but slower than the parametric rate. 

We finally present Theorem~\ref{thm_supp_2} that guarantees the support recovery consistency of $\widetilde \bSigma_{\A}.$

\begin{theorem}
\label{thm_supp_2}
Suppose that Conditions~\ref{cond_T_bd}--\ref{cond_sparse_p_rate} hold and 
$
\big\|\Sigma_{jk}/{\Psi_{jk}^{1/2}}\big\|_\cS > (2\tilde \delta + \tilde\gamma) (\log p/n^{2\gamma_1})^{1/2}
$
for all $(j,k)\in\text{supp}(\bSigma)$ and some $\tilde \gamma>0,$
where $\tilde \delta$ is stated in Theorem~\ref{thm_sparse_2}.  Then we have that
$$\inf_{\Sigma \in \cC_0}P\big\{\text{supp}(\widetilde \bSigma_A) = \text{supp}(\bSigma) \big \} \to 1 \text{  as  } n \to \infty.$$
\end{theorem}

\subsection{Fast computation}
\label{subsec.comp}
Consider a common situation in practice, where, for each $i=1, \dots, n,$ we observe the noisy versions of $X_{i1}(\cdot), \dots, X_{ip}(\cdot)$ at the same set of points, $U_{i1}, \dots, U_{iL_i} \in {\cal U},$ across $j=1, \dots, p.$  Then the original model in (\ref{model.partial}) is simplified to
\begin{equation} 
\label{eq.kernel.simple}
    Z_{ijl} = X_{ij}(U_{il}) + \varepsilon_{ijl}, ~~l=1, \dots, L_i,
\end{equation}
under which the proposed estimation procedure in Section~\ref{subsec.est2} can still be applied.
Suppose that the estimated covariance function is evaluated at a grid of
$R \times R$ locations, $\{(u_{r_1}, u_{r_2}) \in \cU^2: r_1,r_2=1, \dots, R\}.$
To serve the estimation of $p(p+1)/2$ marginal- and cross-covariance functions and the corresponding variance factors, LLSs under the simplified model in (\ref{eq.kernel.simple}) reduce the number of kernel evaluations from $O( \sum_{i= 1}^n \sum_{j = 1}^p L_{ij}R)$ to $O(\sum_{i=1}^n L_i R),$ which substantially accelerate the computation under a high-dimensional regime. 

Apparently, such nonparametric smoothing approach is conceptually simple but suffers from high computational cost in kernel evaluations. To further reduce the computational burden,
we consider fast implementations of LLSs by adopting a simple approximation technique, known as linear binning \cite[]{fan1994}, 
to the covariance function estimation.  
The key idea of the binning method is to greatly reduce the number of kernel evaluations through the fact that many of these evaluations are nearly the same.
We start by dividing $\cU$ into an equally-spaced grid of $R$ points, $u_1< \cdots < u_R \in \cU,$ with binwidth  $\Delta = u_{2}-u_{1}.$  
Denote by $w_{r}(U_{il}) = \max(1-\Delta^{-1}|U_{il} - u_{r}|,0)$ the linear weight that $U_{il}$ assigns to the grid point $u_r$ for $r=1, \dots, R.$ For the $i$-th subject, we define its ``binned weighted counts'' and ``binned weighted averages'' as
$$\varpi_{r,i} = \sum_{l = 1}^{L_{i}} w_{r}(U_{il}) \quad \text{and} \quad  \cD_{r,ij} = \sum_{l = 1}^{L_{i}} w_{r}(U_{il})Z_{ijl},$$ respectively. 
The binned implementation of smoothed adaptive functional thresholding can then be done using this modified dataset $\{(\varpi_{r,i}, \cD_{r,ij})\}_{1 \leq i \leq  n, 1 \leq j \leq p, 1 \leq r \leq R}$ 
and related kernel functions $g_{ab}\{h,(u,v),(u_{r_1},u_{r_2})\}$ for $r_1, r_2 = 1, \dots, R.$ It is notable that, with the help of such binned implementation, the number of kernel evaluations required in the covariance function estimation is further reduced from $O(\sum_{i= 1}^nL_iR)$ to $O(R),$ while only $O(\sum_{i = 1}^n L_i)$ additional operations are involved for each $j$ in the binning step \citep{fan1994}.

We next illustrate the binned implementation of LLS, denoted as BinLLS, using the example of smoothed estimates $\widetilde \Sigma_{jk}$ for $j \neq k$ in (\ref{eq.kernel.cross.cov.solution}).
Under Model~(\ref{eq.kernel.simple}), we drop subscripts $j,k$ in $W_{1,jk},$ $W_{2,jk},$ $W_{3,jk}$ and $S_{ab,jk}$ due to the same set of points $\{U_{i1}, \dots, U_{iL_i}\}$ across $j, k.$
Denote the binned approximations of $ T_{ab,ijk}$ and $S_{ab}$ 
by  $\widecheck T_{ab,ijk}$ and $\widecheck S_{ab},$ respectively. 
It follows from (\ref{T.est}) and (\ref{S.est}) that
$$
    \widecheck T_{ab,ijk}(u,v) = \sum_{r_1 = 1}^R\sum_{r_2 = 1}^R g_{ab}\{h_C,(u,v),(u_{r_1},u_{r_2})\} \cD_{r_1,ij}\cD_{r_2,ik}, 
$$
$$
    \widecheck S_{ab}(u,v)= \sum_{i = 1}^n\sum_{r_1 = 1}^R\sum_{r_2 = 1}^R g_{ab}\{h_C,(u,v),(u_{r_1},u_{r_2})\}  \varpi_{r_1,i} \varpi_{r_2,i}, 
$$
both of which together with (\ref{eq.kernel.cross.cov.solution}) yield the binned approximation of $\widetilde \Sigma_{jk}$ as
\begin{equation}\nonumber
     \widecheck \Sigma_{jk} =\sum_{i=1}^n\big( \widecheck W_{1} \widecheck T_{00,ijk}+ \widecheck W_{2} \widecheck T_{10,ijk}+\widecheck W_{3}\widecheck T_{01,ijk} \big), 
\end{equation}
where $\widecheck W_{1},\widecheck W_{2}$ and $\widecheck W_{3}$ are the binned approximations of $ W_{1}, W_{2}$ and $W_{3},$ computed by replacing the related $S_{ab}$'s in (\ref{eq.kernel.w}) of Supplementary Material with the $\widecheck S_{ab}$'s. 
It is worth noting that, for each pair $(j,k),$ the above binned implementation reduces the number of operations (i.e., additions and multiplications) from $O(R^2 \sum_{i=1}^n L_{i}^2)$ to $O(nR^2 + R^4),$ since the kernel evaluations in $g_{ab}\{h_C,(u,v),(u_{r_1},u_{r_2})\}$ no longer depend on individual observations. Table~\ref{comp.table} presents the computational complexity analysis of LLS and BinLLS under Models~(\ref{model.partial}) and (\ref{eq.kernel.simple}). It reveals that the binned implementation  dramatically improves computational speeds for both densely and sparsely sampled functional data, which is also supported by the empirical evidence in Section~\ref{sec.sim.partial}.

\begin{table}[tbp]
\caption{\label{comp.table} The computational complexity analysis of LLS and BinLLS under Models~(\ref{model.partial}) and (\ref{eq.kernel.simple}) when evaluating the corresponding smoothed covariance function estimates at a grid of $R \times R$ points.}
	\begin{center}
		\vspace{-0.1cm}
		\resizebox{5.1in}{!}{
			\begin{tabular}{cccc}
    \hline
        Method & Model & \begin{tabular}[c]{@{}c@{}}Number of  \\ kernel evaluations\end{tabular}  & \begin{tabular}[c]{@{}c@{}}Number of operations  \\(additions and multiplications)\end{tabular}  \\ \hline
        LLS & (\ref{model.partial}) & $O( \sum_{i= 1}^n \sum_{j = 1}^p L_{ij}R)$& $O(R^2\sum_{i=1}^n\sum_{j,k = 1}^pL_{ij}L_{ik})$  \\ 
        LLS  & (\ref{eq.kernel.simple}) & $O( \sum_{i= 1}^n  L_{i}R)$& $O(p^2R^2\sum_{i=1}^nL_{i}^2)$ \\ 
        BinLLS & (\ref{eq.kernel.simple})& $O(R)$ & $O(np^2R^2+p^2R^4+p\sum_{i = 1}^n L_{i})$ \\ \hline
    \end{tabular}
		}	
	\end{center}
	\vspace{-0.1cm}
\end{table}

To aid the binned implementation of the smoothed adaptive functional thresholding estimator, we then derive the binned approximation of the variance factor $\widetilde \Psi_{jk},$ denoted by  $\widecheck \Psi_{jk}.$ It follows from (\ref{V.est}) that $V_{ab,ijk}$ can be approximated by
$$
\widecheck V_{ab,ijk}(u,v)  = \sum_{r_1 = 1}^R \sum_{r_2 = 1}^{R}g_{ab}\big(h_C,(u,v),(u_{r_1},u_{r_2})\big)\big\{\cD_{r_1,ij}\cD_{r_2,ik} - \widecheck \Sigma_{jk}(u,v)\varpi_{r_1,i} \varpi_{r_2,i}\big\}.    $$
Substituting each term in (\ref{eq.kernel.var}) with its binned approximation, we obtain that
\begin{equation*}
    \widecheck \Psi_{jk}   =  I_{jk} \sum_{i=1}^n\big( \widecheck W_{1}\widecheck V_{00,ijk} +  \widecheck W_{2}\widecheck V_{10,ijk}+  \widecheck W_{3}\widecheck V_{01,ijk}\big)^2.
\end{equation*}
It is worth mentioning that, when $j=k,$
the binned approximations of $\widetilde \Sigma_{jj}$ and $\widetilde \Psi_{jj}$ can be computed in a similar fashion except that the terms corresponding to $r_1=r_2$ should be excluded from all double summations over $\{1, \dots, R\}^2.$
Finally, we obtain the binned adaptive functional thresholding estimator
$\widecheck \bSigma_{\A} = (\widecheck \Sigma^\A_{jk})_{p \times p}$ with
$ \widecheck \Sigma_{jk}^\A= \widecheck \Psi_{jk}  ^{1/2} \times s_{\lambda}\big({\widecheck \Sigma_{jk}}/{\widecheck \Psi_{jk}  ^{1/2}}\big)$ and the corresponding universal thresholding estimator $\widecheck \bSigma_{\U} = (\widecheck \Sigma^\U_{jk})_{p \times p}$ with $\widecheck \Sigma_{jk}^{\U} = s_{\lambda}\big(\widecheck \Sigma_{jk}\big).$

\section{Simulations}
\label{sec.sim}
\subsection{Setup} \label{sec.sim.setup}
We conduct a number of simulations to compare adaptive functional thresholding estimators to universal functional thresholding estimators. Sections \ref{sec.sim.full} and \ref{sec.sim.partial} consider scenarios where random functions are fully and partially observed, respectively.

In each scenario, to mimic the infinite-dimensionality of random curves, we generate functional variables by $X_{ij}(u)=\bs(u)^{\T}\btheta_{ij}$ for $i = 1, \dots, n, j=1, \dots, p$ and $u \in \cU=[0,1],$ where $\bs(u)$ is a $50$-dimensional Fourier basis function and $\btheta_{i}=(\btheta_{i1}^{\T},\dots,\btheta_{ip}^{\T} )^{\T} \in \eR^{50p}$ is generated from a mean zero multivariate Gaussian distribution with block covariance matrix
$\bOmega \in \eR^{50p \times 50p},$ whose $(j,k)$-th block is $\bOmega_{jk} \in \eR^{50 \times 50}$ for $j,k=1, \dots, p.$
The functional sparsity pattern in $\bSigma= \{\Sigma_{jk}(\cdot,\cdot)\}_{p \times p}$ with its $(j,k)$th entry $\Sigma_{jk}(u,v) = \bs(u)^{\T}\bOmega_{jk}\bs(v)$ can be characterized by the block sparsity structure in $\bOmega.$ Define $\bOmega_{jk} = \omega_{jk} \bD$ with $\bD = \text{diag}(1^{-2},\dots, 50^{-2})$ and hence $\cov(\theta_{ijk},\theta_{ijk'}) \sim k^{-2}I(k=k')$ for $k,k'=1, \dots, 50.$ Then we generate $\bOmega$ with different block sparsity patterns as follows.

\begin{itemize}
    \item Model~1 (block banded). For $j,k = 1, \dots, p/2$, $\omega_{jk} = (1-{|j-k|}/{10})_{+}$. For $ j,k = p/2+1, \dots, p$, $\omega_{jk} = 4I(j=k).$
    \item Model~2 (block sparse without any special structure). For $ j,k = p/2+1, \dots, p$, $\omega_{jk} = 4I(j=k).$ For $j,k = 1, \dots, p/2$, we generate $\bomega = (\omega_{jk})_{p/2 \times p/2} = \bB + \delta' \bI_{p/2},$ where elements of $\bB$ are sampled independently from $\text{Uniform}[0.3,0.8]$ with probability $0.2$ or $0$ with probability $0.8,$ and $\delta' = \{-\lambda_{\min}(\bB),0\} + 0.01$ to guarantee the positive definiteness of $\bOmega.$
\end{itemize}

We implement a cross-validation approach \cite[]{bickel2008} for choosing the optimal thresholding parameter $\hat \lambda$ in $\widehat\bSigma_{\A}$. Specifically, we randomly divide the sample $\{\bX_i: i=1, \dots, n\}$ into two subsamples of size $n_1$ and $n_2,$ where $n_1=n(1-1/\log n)$ and $n_2=n/\log n$ and repeat this $N$ times.
Let $\widehat\bSigma_{\A,1}^{(\nu)}(\lambda)$ and $\widehat\bSigma_{\tS,2}^{(\nu)}$ be the adaptive functional thresholding estimator as a function of $\lambda$ and the sample covariance function based on $n_1$ and $n_2$ observations, respectively, from the $\nu$th split. We select the optimal $\hat \lambda$ by minimizing $$\widehat \err(\lambda) = N^{-1}\sum_{\nu=1}^N \|\widehat\bSigma_{\A,1}^{(\nu)}(\lambda) - \widehat\bSigma_{\tS,2}^{(\nu)}\|_{\tF}^2,$$ where $\|\cdot\|_{\tF}$ denotes the functional version of Frobenius norm, i.e., for any $\boldsymbol{ Q} = \{Q_{jk}(\cdot,\cdot)\}_{p \times p}$ with each $Q_{jk} \in {\mathbb S},$ $\|\boldsymbol{ Q} \|_{\tF} =(\sum_{j,k}\|Q_{jk}\|_{\cS}^2)^{1/2}.$
The optimal thresholding parameters in $\widehat\bSigma_{\U},$ $\widetilde\bSigma_{\A}, \widetilde \bSigma_{\U}, \widecheck\bSigma_{\A},\widecheck\bSigma_{\U}$ can be selected in a similar fashion.

\subsection{Fully observed functional data}
\label{sec.sim.full}
We compare the adaptive functional thresholding estimator 
$\widehat\bSigma_{\A}$ to the universal functional thresholding estimator 
$\widehat\bSigma_{\U}$ under hard, soft, SCAD (with $a=3.7$) and adaptive lasso (with $\eta=3$) functional thresholding rules, where the corresponding $\hat \lambda$'s are selected by the cross-validation with $N=5.$ 
We generate $n=100$ observations for $p=50, 100, 150$ and replicate each simulation 100 times. We examine the performance of all competing approaches by estimation and support recovery accuracies. In terms of the estimation accuracy, Table~\ref{err.table} reports numerical summaries of losses measured by functional versions of Frobenius and matrix $\ell_1$ norms. 
To assess the support recovery consistency, we present in Table~\ref{TPR/FPR.table} the average of true positive rates (TPRs) and false positive rates (FPRs), defined as $\text{TPR} = \#\{(j,k): \|\widehat \Sigma_{jk}\|_\cS \neq 0 ~\text{and}~\| \Sigma_{jk}\|_\cS\neq 0 \}/\# \{(j,k): \| \Sigma_{jk}\|_\cS\neq 0\}$ and $\text{FPR} = \#\{(j,k): \|\widehat \Sigma_{jk}\|_\cS \neq 0 ~\text{and}~\| \Sigma_{jk}\|_\cS= 0 \}/\# \{(j,k): \| \Sigma_{jk}\|_\cS= 0\}.$

\begin{table}[tbp]
	\caption{\label{err.table}The average (standard error) functional matrix losses over 100 simulation runs.}
	\begin{center}
		\resizebox{6.2in}{!}{
			\begin{tabular}{ccrrrrrr}
			\hline															
	&		&	\multicolumn{2}{c}{$p = 50$}			&	\multicolumn{2}{c}{$p = 100$}			&	\multicolumn{2}{c}{$p = 150$}			\\	
{Model} 	&	{Method}	&	\multicolumn{1}{c}{$\widehat \bSigma_\A$}	&	\multicolumn{1}{c}{$\widehat \bSigma_\U$}	&	\multicolumn{1}{c}{$\widehat \bSigma_\A$}	&	\multicolumn{1}{c}{$\widehat \bSigma_\U$}	&	\multicolumn{1}{c}{$\widehat \bSigma_\A$}	&	\multicolumn{1}{c}{$\widehat \bSigma_\U$}	\\	\hline
\multirow{12}{*}{1}	&	 \multicolumn{7}{c}{Functional Frobenius norm}												\\		
	&	Hard	&	5.40(0.04)	&	11.90(0.02)	&	7.91(0.03)	&	17.27(0.01)	&	9.94(0.04)	&	21.36(0.01)	\\	
	&	Soft	&	6.28(0.05)	&	10.40(0.08)	&	9.41(0.05)	&	16.53(0.07)	&	11.85(0.06)	&	21.16(0.04)	\\	
	&	SCAD	&	5.68(0.05)	&	10.56(0.08)	&	8.53(0.05)	&	16.59(0.07)	&	10.80(0.06)	&	21.19(0.04)	\\	
	&	Adap. lasso	&	5.28(0.04)	&	11.42(0.07)	&	7.76(0.04)	&	17.26(0.01)	&	9.72(0.04)	&	21.36(0.01)		\\
	&	Sample	&	\multicolumn{2}{c}{	19.82(0.04)	}	&	\multicolumn{2}{c}{	39.54(0.05)	}	&	\multicolumn{2}{c}{	59.28(0.06)	}		
	\\	\cline{2-8}														
	&	 \multicolumn{7}{c}{Functional matrix $\ell_1$ norm}													\\	
	&	Hard	&	3.96(0.06)	&	9.23(0.01)	&	4.49(0.05)	&	9.31(0.01)	&	4.78(0.05)	&	9.34(0.01)	\\	
	&	Soft	&	5.04(0.07)	&	8.14(0.08)	&	5.88(0.05)	&	9.15(0.02)	&	6.21(0.04)	&	9.31(0.01)	\\	
	&	SCAD	&	4.40(0.08)	&	8.32(0.07)	&	5.35(0.06)	&	9.18(0.02)	&	5.75(0.05)	&	9.31(0.01)	\\	
	&	Adap.lasso	&	3.85(0.06)	&	8.91(0.07)	&	4.52(0.05)	&	9.30(0.01)	&	4.83(0.06)	&	9.34(0.01)	\\	
	&	Sample	&	\multicolumn{2}{c}{	26.60(0.13)	}	&	\multicolumn{2}{c}{	52.65(0.18)	}	&	\multicolumn{2}{c}{	78.69(0.22)	}		
	\\	\hline														
\multirow{12}{*}{2}	&	 \multicolumn{7}{c}{Functional Frobenius norm}												\\		
	&	Hard	&	5.67(0.03)	&	9.39(0.02)	&	9.48(0.04)	&	15.79(0.01)	&	14.00(0.05)	&	22.26(0.01)	\\	
	&	Soft	&	6.14(0.03)	&	8.55(0.04)	&	10.28(0.05)	&	15.00(0.05)	&	14.8(0.05)	&	21.89(0.04)	\\	
	&	SCAD	&	5.94(0.03)	&	8.59(0.04)	&	9.96(0.05)	&	15.02(0.05)	&	14.49(0.06)	&	21.91(0.04)	\\	
	&	Adap. lasso	&	5.44(0.03)	&	9.10(0.04)	&	8.99(0.04)	&	15.73(0.02)	&	13.02(0.05)	&	22.25(0.01)	\\	
	&	Sample	&	\multicolumn{2}{c}{	21.80(0.04)	}	&	\multicolumn{2}{c}{	43.51(0.06)	}	&	\multicolumn{2}{c}{	65.22(0.07)	}		
\\	\cline{2-8}															
	&	 \multicolumn{7}{c}{Functional matrix $\ell_1$ norm}													\\	
	&	Hard	&	2.85(0.03)	&	4.74(0.01)	&	4.77(0.05)	&	7.11(0.01)	&	7.65(0.07)	&	10.31(0.01)	\\	
	&	Soft	&	3.31(0.03)	&	4.51(0.04)	&	5.37(0.04)	&	6.90(0.02)	&	8.21(0.05)	&	10.21(0.01)	\\	
	&	SCAD	&	3.22(0.03)	&	4.48(0.03)	&	5.29(0.04)	&	6.91(0.02)	&	8.14(0.05)	&	10.21(0.01)	\\	
	&	Adap. lasso	&	2.75(0.03)	&	4.66(0.02)	&	4.62(0.05)	&	7.08(0.01)	&	7.35(0.07)	&	10.30(0.01)	\\	
	&	Sample	&	\multicolumn{2}{c}{	28.06(0.12)	}	&	\multicolumn{2}{c}{	56.01(0.19)	}	&	\multicolumn{2}{c}{	84.13(0.23)	}		
	\\	\hline

			\end{tabular}
		}	
	\end{center}
\end{table}

\begin{table}[tbp]
	\caption{\label{TPR/FPR.table} The average TPRs/ FPRs over 100 simulation runs.}
	\begin{center}
		\resizebox{6.2in}{!}{
			\begin{tabular}{cccccccc}
			\hline
	&			&	\multicolumn{2}{c}{$p = 50$}			&	\multicolumn{2}{c}{$p = 100$}			&	\multicolumn{2}{c}{$p = 150$}			\\
Model	&	Method	&	$\widehat \bSigma_\A$	&	$\wbSigma_\U$	&	$\widehat \bSigma_\A$	&	$\wbSigma_\U$	&	$\widehat \bSigma_\A$	&	$\wbSigma_\U$	\\	\hline	
																	
\multirow{4}{*}{1}	&	Hard	&	0.71/0.00	&	0.00/0.00	&	0.66/0.00	&	0.00/0.00	&	0.64/0.00	&	0.00/0.00	\\		
	&	Soft	&	0.89/0.08	&	0.47/0.17	&	0.85/0.04	&	0.22/0.05	&	0.84/0.03	&	0.06/0.01	\\		
	&	SCAD	&	0.89/0.07	&	0.42/0.13	&	0.85/0.04	&	0.20/0.04	&	0.84/0.03	&	0.05/0.01	\\		
	&	Adap. lasso	&	0.78/0.00	&	0.11/0.02	&	0.74/0.00	&	0.00/0.00	&	0.73/0.00	&	0.00/0.00	\\	\hline	
\multirow{4}{*}{2}	&	Hard	&	0.77/0.00	&	0.00/0.00	&	0.68/0.00	&	0.00/0.00	&	0.63/0.00	&	0.00/0.00	\\		
	&	Soft	&	0.99/0.06	&	0.50/0.07	&	0.97/0.04	&	0.30/0.04	&	0.96/0.04	&	0.11/0.02	\\		
	&	SCAD	&	0.99/0.06	&	0.47/0.06	&	0.98/0.05	&	0.29/0.04	&	0.97/0.05	&	0.10/0.01	\\		
	&	Adap. lasso	&	0.91/0.00	&	0.10/0.01	&	0.86/0.00	&	0.01/0.00	&	0.83/0.00	&	0.00/0.00	\\	\hline	

			\end{tabular}
		}	
	\end{center}
\end{table}

Several conclusions can be drawn from Tables~\ref{err.table} and \ref{TPR/FPR.table}. First, in all scenarios, $\widehat \bSigma_{\A}$ provides substantially improved accuracy over $\widehat \bSigma_{\U}$ regardless of the thresholding rule or the loss  used. 
We also obtain the sample covariance function $\widehat\bSigma_{\tS},$ the results of which deteriorate severely compared with $\widehat \bSigma_{\A}$ and $\widehat \bSigma_{\U}.$
Second, for support recovery, again $\widehat \bSigma_{\A}$ uniformly outperforms $\widehat \bSigma_{\U}$, which fails to recover the functional sparsity pattern especially when $p$ is large.
Third, the adaptive functional thresholding approach using the hard and the adaptive lasso functional thresholding rules tends to have lower losses and lower TPRs/FPRs than that using the soft and the SCAD functional thresholding rules.

\subsection{Partially observed functional data} 
\label{sec.sim.partial}
In this section, we assess the finite-sample performance of  LLS and  BinLLS methods to handle partially observed functional data.  
We first generate random functions $X_{ij}(\cdot)$ for $i=1, \dots, n, j=1, \dots, p$ by the same procedure as in Section~\ref{sec.sim.setup} with either non-sparse or sparse $\bSigma$ depending on $p.$ We then generate the observed values $Z_{ijl}$ from equation (\ref{eq.kernel.simple}), where the measurement locations $U_{il}$ and errors $\varepsilon_{ijl}$ are sampled independently from  $\text{Uniform[0,1]}$ and ${\cal N}(0,0.5^2),$ respectively.
We consider settings of $n=100$ and $L_i = 11, 21, 51, 101,$ changing from sparse to moderately dense to very dense measurement schedules. We use the Gaussian kernel with the optimal bandwidths proportional to $n^{-1/6},$ $(nL_i^2)^{-1/6}$ and $n^{-1/4},$ respectively, as suggested in \cite{Zhang2016}, so for the empirical work in this paper we choose the proportionality constants in the range $(0,1],$ which gives good results in all settings we consider.

To compare BinLLS with LLS in terms of the computational speed and estimation accuracy, we first consider a low-dimensional example $p = 6$ with non-sparse $\bSigma$ generated by modifying Model~1 with $\omega_{jk} = (1-{|j-k|}/{10})_{+}$ for $j,k = 1, \dots, 6.$ In addition to our proposed smoothing methods, we also implement local-linear-smoother-based pre-smoothing and its binned implementation, denoted as LLS-P and BinLLS-P, respectively. Table \ref{table.kernel.6} reports numerical summaries of estimation errors evaluated at $R=21$ equally-spaced points in $[0,1]$ and the corresponding CPU time on the processor Intel(R) Xeon(R) CPU E5-2690 v3 @ 2.60GHz. The results for the sample covariance function $\widehat \bSigma_{\tS}$ based on fully observed $\bX_1(\cdot), \dots, \bX_n(\cdot)$ are also provided as the baseline for comparison. Note that, LLS is too slow to implement for the case $L_i=101,$ so we do not report its result here.


A few trends are observable from Table~\ref{table.kernel.6}.
First, the binned implementations (BinLLS and BinLLS-P) attain similar or even lower estimation errors compared with their direct implementations (LLS and LLS-P) under all scenarios, while resulting in considerably faster computational speeds especially under dense designs. For example, BinLLS runs over $400$ times faster than LLS when $L_i=51.$ 
Second, all methods provide higher estimation accuracies as $L_i$ increases, and enjoy similar performance when functions are very densely observed, e.g., $L_i=51$ and $101,$ compared with the fully observed functional case.
However, the performance of LLS-P and BinLLS-P deteriorates severely under sparse designs, e.g., $L_i=11$ and $21,$ since limited information is available from a small number of observations per subject.
Among all competitors, we conclude that BinLLS is overall a unified approach that can handle both sparsely and densely sampled functional data well with increased computational efficiency and guaranteed estimation accuracy.



\begin{table}[tbp]
	\caption{ The average (standard error) functional matrix losses and average CPU time for $p = 6$ over 100 simulation runs. \label{table.kernel.6}}
	\begin{center}
		\resizebox{6.5in}{!}{
\begin{tabular}{ccccccccc}
\hline
\multicolumn{1}{c|}{$L_i$}                         & Method                  & \begin{tabular}[c]{@{}c@{}}Functional \\ Frobenius norm\end{tabular} & \begin{tabular}[c]{@{}c@{}}Functional \\ matrix $\ell_1$ norm\end{tabular} & \multicolumn{1}{c|}{\begin{tabular}[c]{@{}c@{}}Elapsed time \\ (sec)\end{tabular}} & Method   & \begin{tabular}[c]{@{}c@{}}Functional \\ Frobenius norm\end{tabular} & \begin{tabular}[c]{@{}c@{}}Functional \\ matrix $\ell_1$ norm\end{tabular} & \begin{tabular}[c]{@{}c@{}}Elapsed time \\ (sec)\end{tabular} \\ \hline
\multicolumn{1}{c|}{\multirow{2}{*}{11}}  & BinLLS                  & 1.57(0.02)                                                         & 1.72(0.03)                                                               & \multicolumn{1}{c|}{2.06}                                                         & BinLLS-P & 4.14(0.03)                                                         & 4.36(0.04)                                                               & 0.18                                                         \\
\multicolumn{1}{c|}{}                         & LLS                     & 1.62(0.02)                                                         & 1.76(0.03)                                                               & \multicolumn{1}{c|}{50.52}                                                        & LLS-P     & 4.23(0.04)                                                         & 4.47(0.05)                                                               & 0.22                                                         \\
\multicolumn{1}{c|}{\multirow{2}{*}{21}}  & BinLLS                  & 1.28(0.02)                                                         & 1.42(0.03)                                                               & \multicolumn{1}{c|}{2.07}                                                         & BinLLS-P & 2.66(0.02)                                                         & 2.80(0.02)                                                               & 0.19                                                         \\
\multicolumn{1}{c|}{}                         & LLS                     & 1.28(0.02)                                                         & 1.42(0.03)                                                               & \multicolumn{1}{c|}{136.88}                                                        & LLS-P     & 2.67(0.02)                                                         & 2.82(0.03)                                                               & 0.29                                                         \\
\multicolumn{1}{c|}{\multirow{2}{*}{51}}  & BinLLS                  & 1.06(0.02)                                                         & 1.20(0.03)                                                               & \multicolumn{1}{c|}{2.21}                                                         & BinLLS-P & 1.12(0.03)                                                         & 1.26(0.03)                                                               & 0.20                                                         \\
\multicolumn{1}{c|}{}                         & LLS                     & 1.04(0.02)                                                         & 1.18(0.03)                                                               & \multicolumn{1}{c|}{967.75}                                                        & LLS-P     & 1.12(0.03)                                                         & 1.26(0.03)                                                               & 0.39                                                         \\
\multicolumn{1}{c|}{\multirow{2}{*}{101}} & BinLLS & 1.00(0.02)                                       & 1.14(0.03)                                              & \multicolumn{1}{c|}{2.23}                                        & BinLLS-P & 0.99(0.02)                                                         & 1.13(0.03)                                                               & 0.21                                                        \\
\multicolumn{1}{c|}{}                         &        LLS                 &                                                              -        &          -                                                                  & \multicolumn{1}{c|}{-}                                                              & LLS-P     & 0.97(0.02)                                                         & 1.11(0.03)                                                               & 0.64                                                       \\ \hline
\multicolumn{2}{c}{\multirow{2}{*}{$\widehat   \bSigma_\tS$}}      & \multicolumn{2}{c}{Functional Frobenius norm}                                                                                                                               & \multicolumn{3}{c}{Functional matrix $\ell_1$ norm}                                           & \multicolumn{2}{c}{Elapsed time (sec)}                                                                                                                                                                            \\
         \multicolumn{2}{c}{}                                           & \multicolumn{2}{c}{1.04(0.03)}                                                                                                                                            & \multicolumn{3}{c}{1.20(0.03)}                                                              & \multicolumn{2}{c}{0.11}                                                                                                                                                                                         \\ \hline
\end{tabular}
		}	
	\end{center}
\end{table}

\begin{table}[!htbp]
	\caption{\label{err.table.kernel.50} The average (standard error) functional matrix losses for partially observed functional scenarios and $p=50$ over 100 simulation runs.}
	\begin{center}
		\resizebox{6.5in}{!}{
			\begin{tabular}{ccrrrrrrrr}
			\hline																			
	&		&	\multicolumn{2}{c}{$L_i = 11$}			&	\multicolumn{2}{c}{$L_i = 21$}			&	\multicolumn{2}{c}{$L_i = 51$}			&	\multicolumn{2}{c}{$L_i = 101$}			\\	
{Model} 	&	{Method}	&	\multicolumn{1}{c}{$\widecheck \bSigma_\A$}	&	\multicolumn{1}{c}{$\widecheck \bSigma_\U$}	&	\multicolumn{1}{c}{$\widecheck \bSigma_\A$}	&	\multicolumn{1}{c}{$\widecheck \bSigma_\U$}	&	\multicolumn{1}{c}{$\widecheck \bSigma_\A$}	&	\multicolumn{1}{c}{$\widecheck \bSigma_\U$}	&	\multicolumn{1}{c}{$\widecheck \bSigma_\A$}	&	\multicolumn{1}{c}{$\widecheck \bSigma_\U$}	\\	\hline
\multirow{12}{*}{1}	&	 \multicolumn{9}{c}{Functional Frobenius norm}																\\		
	&	Hard	&	7.78(0.03)	&	12.65(0.01)	&	6.61(0.04)	&	12.26(0.01)	&	5.83(0.04)	&	12.04(0.02)	&	5.57(0.04)	&	11.89(0.04)	\\	
	&	Soft	&	8.69(0.04)	&	12.63(0.01)	&	7.64(0.05)	&	11.75(0.06)	&	6.94(0.05)	&	10.51(0.07)	&	6.71(0.05)	&	10.05(0.07)	\\	
	&	SCAD	&	8.36(0.05)	&	12.63(0.01)	&	7.13(0.05)	&	11.80(0.06)	&	6.28(0.05)	&	10.67(0.07)	&	5.99(0.05)	&	10.27(0.07)	\\	
	&	Adap. lasso	&	7.69(0.04)	&	12.64(0.01)	&	6.57(0.04)	&	12.21(0.02)	&	5.83(0.04)	&	11.54(0.08)	&	5.57(0.04)	&	11.05(0.10)	\\	
\cline{2-10}																		
	&	 \multicolumn{9}{c}{Functional matrix $\ell_1$ norm}																	\\	
	&	Hard	&	5.35(0.05)	&	9.36(0.01)	&	4.68(0.06)	&	9.30(0.01)	&	4.09(0.06)	&	9.24(0.02)	&	3.87(0.06)	&	9.13(0.05)	\\	
	&	Soft	&	6.38(0.06)	&	9.35(0.01)	&	5.86(0.07)	&	8.94(0.05)	&	5.43(0.07)	&	8.13(0.08)	&	5.29(0.07)	&	7.84(0.08)	\\	
	&	SCAD	&	6.12(0.07)	&	9.35(0.01)	&	5.40(0.08)	&	8.99(0.05)	&	4.78(0.08)	&	8.32(0.07)	&	4.56(0.08)	&	8.09(0.07)	\\	
	&	Adap.lasso	&	5.31(0.07)	&	9.36(0.01)	&	4.71(0.07)	&	9.28(0.02)	&	4.15(0.07)	&	8.89(0.07)	&	3.98(0.07)	&	8.59(0.09)	\\	
\hline								\multirow{12}{*}{2}	&	 \multicolumn{9}{c}{Functional Frobenius norm}																\\	
	&	Hard	&	8.12(0.03)	&	10.41(0.02)	&	6.85(0.04)	&	9.89(0.01)	&	6.06(0.04)	&	9.60(0.02)	&	5.75(0.04)	&	9.51(0.02)	\\
	&	Soft	&	8.35(0.03)	&	10.37(0.01)	&	7.35(0.03)	&	9.60(0.03)	&	6.72(0.03)	&	8.86(0.04)	&	6.48(0.03)	&	8.56(0.04)	\\
	&	SCAD	&	8.32(0.03)	&	10.37(0.01)	&	7.23(0.04)	&	9.60(0.03)	&	6.50(0.04)	&	8.89(0.04)	&	6.23(0.04)	&	8.61(0.04)	\\
	&	Adap. lasso	&	7.83(0.03)	&	10.39(0.01)	&	6.69(0.04)	&	9.84(0.02)	&	5.97(0.04)	&	9.40(0.04)	&	5.71(0.04)	&	9.16(0.04)	\\
\cline{2-10}																		
	&	 \multicolumn{9}{c}{Functional matrix $\ell_1$ norm}																	\\
	&	Hard	&	3.82(0.04)	&	4.91(0.01)	&	3.36(0.04)	&	4.82(0.01)	&	3.00(0.05)	&	4.78(0.01)	&	2.85(0.05)	&	4.77(0.01)	\\
	&	Soft	&	3.96(0.02)	&	4.88(0.01)	&	3.71(0.03)	&	4.72(0.02)	&	3.50(0.03)	&	4.55(0.03)	&	3.44(0.03)	&	4.47(0.03)	\\
	&	SCAD	&	3.96(0.02)	&	4.88(0.01)	&	3.67(0.03)	&	4.72(0.02)	&	3.41(0.03)	&	4.55(0.02)	&	3.32(0.03)	&	4.48(0.02)	\\
	&	Adap. lasso	&	3.65(0.04)	&	4.90(0.01)	&	3.28(0.04)	&	4.80(0.01)	&	2.96(0.04)	&	4.73(0.01)	&	2.88(0.04)	&	4.69(0.02)	\\
\hline

			\end{tabular}
		}	
	\end{center}
\end{table}

\begin{table}[!htbp]
	\caption{\label{TPR/FPR.table.kernel.50} The average TPRs/ FPRs for partially observed functional scenarios and $p=50$ over 100 simulation runs.}
	\begin{center}
		\resizebox{6.5in}{!}{
			\begin{tabular}{cccccccccc}
	\hline																				
	&			&	\multicolumn{2}{c}{$L_i = 11$}			&	\multicolumn{2}{c}{$L_i = 21$}			&	\multicolumn{2}{c}{$L_i=51$}			&	\multicolumn{2}{c}{$L_i=101$}			\\
Model	&	Method	&	$\widecheck \bSigma_\A$	&	$\widecheck \bSigma_\U$	&	$\widecheck \bSigma_\A$	&	$\widecheck \bSigma_\U$	&	$\widecheck \bSigma_\A$	&	$\widecheck \bSigma_\U$	&	$\widecheck \bSigma_\A$	&	$\widecheck \bSigma_\U$	\\	\hline
																				
\multirow{4}{*}{1}	&	Hard	&	0.63/0.00	&	0.00/0.00	&	0.66/0.00	&	0.00/0.00	&	0.69/0.00	&	0.01/0.00	&	0.71/0.00	&	0.03/0.00	\\	
	&	Soft	&	0.85/0.05	&	0.01/0.00	&	0.87/0.07	&	0.22/0.09	&	0.89/0.08	&	0.5/0.17	&	0.89/0.08	&	0.57/0.18	\\	
	&	SCAD	&	0.86/0.06	&	0.01/0.00	&	0.87/0.07	&	0.2/0.07	&	0.88/0.07	&	0.45/0.14	&	0.89/0.07	&	0.51/0.14	\\	
	&	Adap. lasso	&	0.72/0.00	&	0.00/0.00	&	0.75/0.00	&	0.01/0.00	&	0.77/0.00	&	0.12/0.02	&	0.78/0.00	&	0.20/0.03	\\	\hline
\multirow{4}{*}{2}	&	Hard	&	0.58/0.00	&	0.00/0.00	&	0.69/0.00	&	0.00/0.00	&	0.75/0.00	&	0.01/0.00	&	0.79/0.00	&	0.01/0.00	\\	
	&	Soft	&	0.95/0.04	&	0.03/0.01	&	0.97/0.05	&	0.22/0.03	&	0.99/0.06	&	0.48/0.06	&	0.99/0.06	&	0.58/0.07	\\	
	&	SCAD	&	0.95/0.04	&	0.03/0.01	&	0.97/0.06	&	0.22/0.03	&	0.99/0.07	&	0.46/0.06	&	0.99/0.07	&	0.54/0.06	\\	
	&	Adap. lasso	&	0.80/0.00	&	0.00/0.00	&	0.86/0.00	&	0.02/0.00	&	0.90/0.00	&	0.08/0.00	&	0.91/0.00	&	0.15/0.01	\\	\hline

			\end{tabular}
		}	
	\end{center}
\end{table}

We next examine the performance of BinLLS-based adaptive and universal functional thresholding estimators in terms of estimation accuracy and support recovery consistency using the same performance measures as in
Tables~\ref{err.table}--\ref{TPR/FPR.table}.
Tables~\ref{err.table.kernel.50}--\ref{TPR/FPR.table.kernel.50} and Tables~\ref{err.table.kernel.100}--\ref{TPR/FPR.table.kernel.100} of Supplementary Material report numerical results for settings of $p=50$ and $100,$ respectively, satisfying Models~1 and 2 under different measurement schedules. We observe a few apparent patterns.
First, $\widecheck \bSigma_\A$ substantially outperforms $\widecheck\bSigma_\U$ with significantly lower estimation errors in all settings.
Second, $\widecheck \bSigma_\A$ works consistently well in recovering the functional sparsity structures especially under the soft and SCAD functional thresholding rules, while $\widecheck \bSigma_\U$ fails to identify such patterns.
Third, the estimation and support recovery consistencies of $\widecheck \bSigma_\A$ and $\widecheck \bSigma_\U$ are improved as $L_i$ increases. When curves are very densely observed, e.g., $L_i=101,$ we observe that both estimators enjoy similar performance with $\widehat\bSigma_\A$ and $\widehat\bSigma_\U$ in Tables~\ref{err.table} and \ref{TPR/FPR.table}. Such observation provides empirical evidence to support our remark for Theorem~\ref{thm_sparse_2} about the same convergence rate between very densely observed and fully observed functional scenarios.

\section{Real Data}
\label{sec.real}

\subsection{ADHD dataset}
\label{sec.real.ADHD}
In this section, we illustrate our adaptive functional thresholding estimation using the ADHD-200 Sample, collected by New York University Medical Center.
This dataset consists of  resting-state fMRI scans with Blood Oxygenation Level-Dependent (BOLD) signals recorded every 2 seconds in the whole brain with $L=172$ locations in total,
for $n_\ADHD= 90$ patients diagnosed with attention-deficit/hyperactivity disorder (ADHD) and $n_\TDC = 87$ typically-developing controls (TDC). 
The preprocessing of the raw fMRI data is performed by Neuro Bureau using the Athena pipeline \cite[]{bellec2017}. See Figure~\ref{data_ADHD} of Supplementary Material for plots of pre-smoothed BOLD signals at a selection of regions of interest (ROIs). Following \cite{li2018} based on the same dataset, we treat the signals at different ROIs as multivariate functional data. 
Our goal is to construct resting state functional connectivity networks among $p = 116$ ROIs \cite[]{tzourio2002}, with the first 90 ROIs from the cerebrum and the last 26 ROIs from the cerebellum, for ADHD and TDC groups, respectively. 
To this end, we implement adaptive and universal functional thresholding methods to discover the networks for two groups.

\begin{figure}[tbp]

\centering
\begin{subfigure}{0.3\linewidth}
  \includegraphics[width=4.5cm,height=4.5cm]{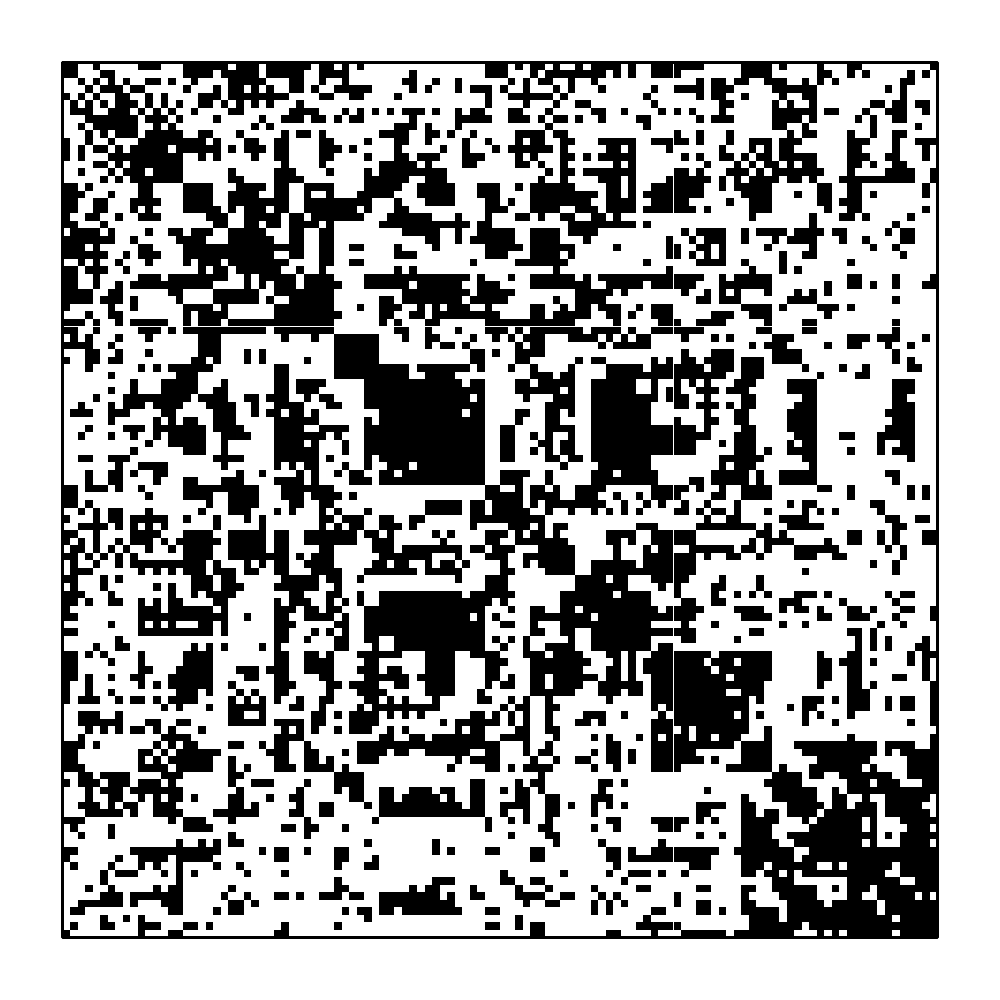}
  \caption{ADHD: $\wbSigma_\A$ ($57.50\%$ zeros)} \label{fig:1a}
  \par\medskip 
  \includegraphics[width=4.5cm,height=4.5cm]{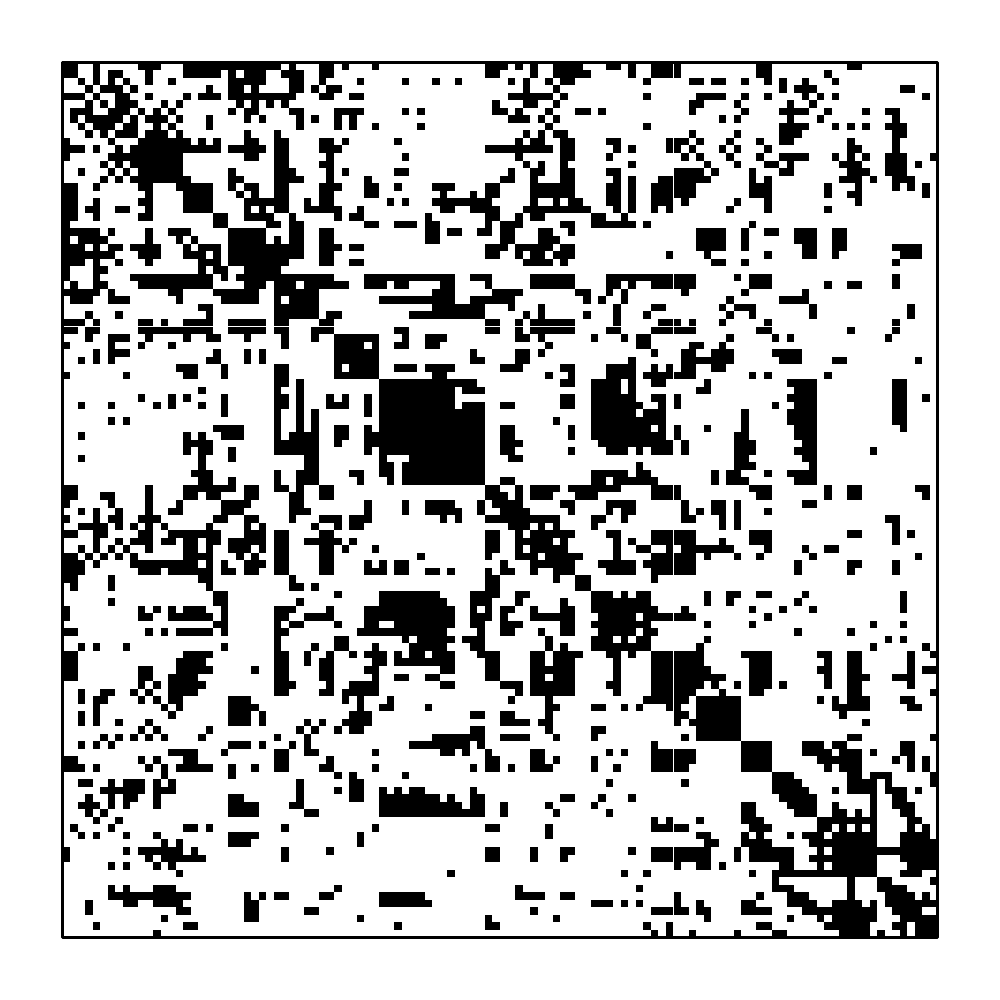}
  \caption{TDC: $\wbSigma_\A$ ($71.24\%$ zeros)} \label{fig:1b}
\end{subfigure}
\begin{subfigure}{0.3\linewidth}
  \includegraphics[width=4.5cm,height=4.5cm]{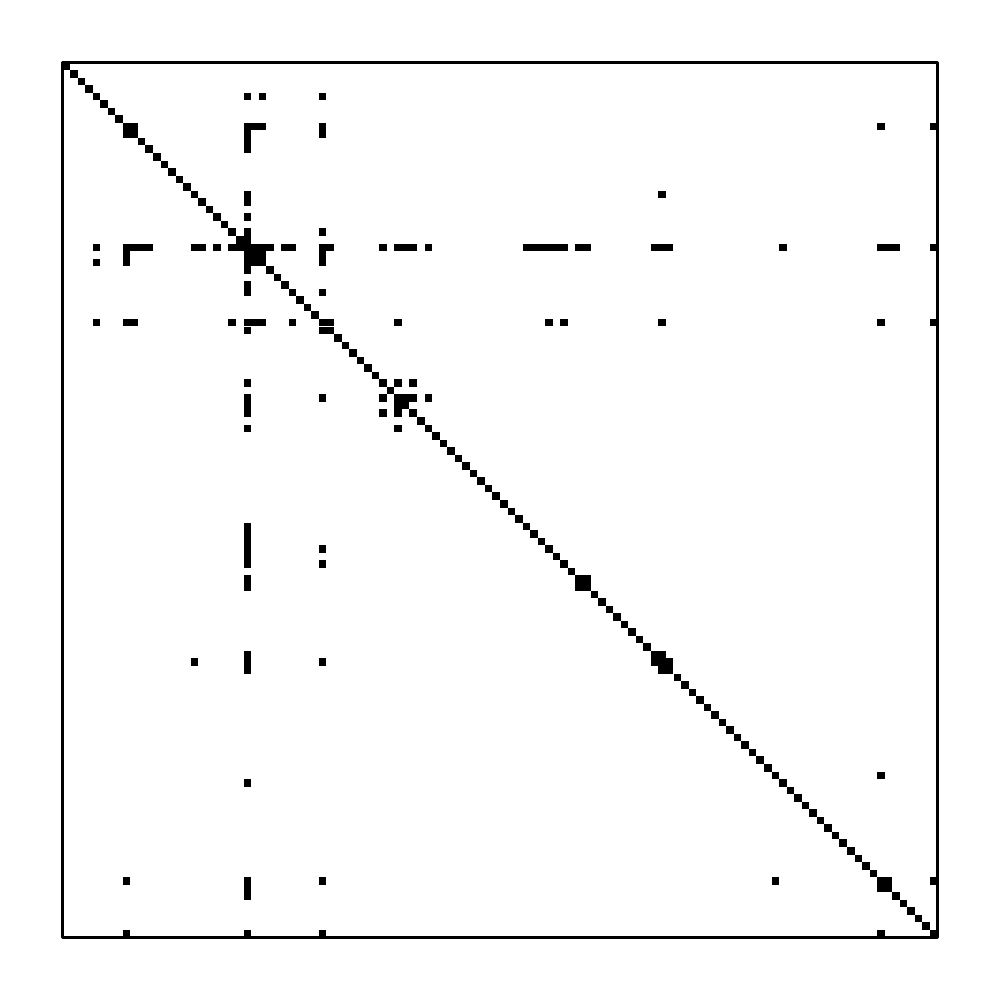}
  \caption{ADHD: $\wbSigma_\U$ ($98.94\%$ zeros)} \label{fig:1c}
  \par\medskip 
  \includegraphics[width=4.5cm,height=4.5cm]{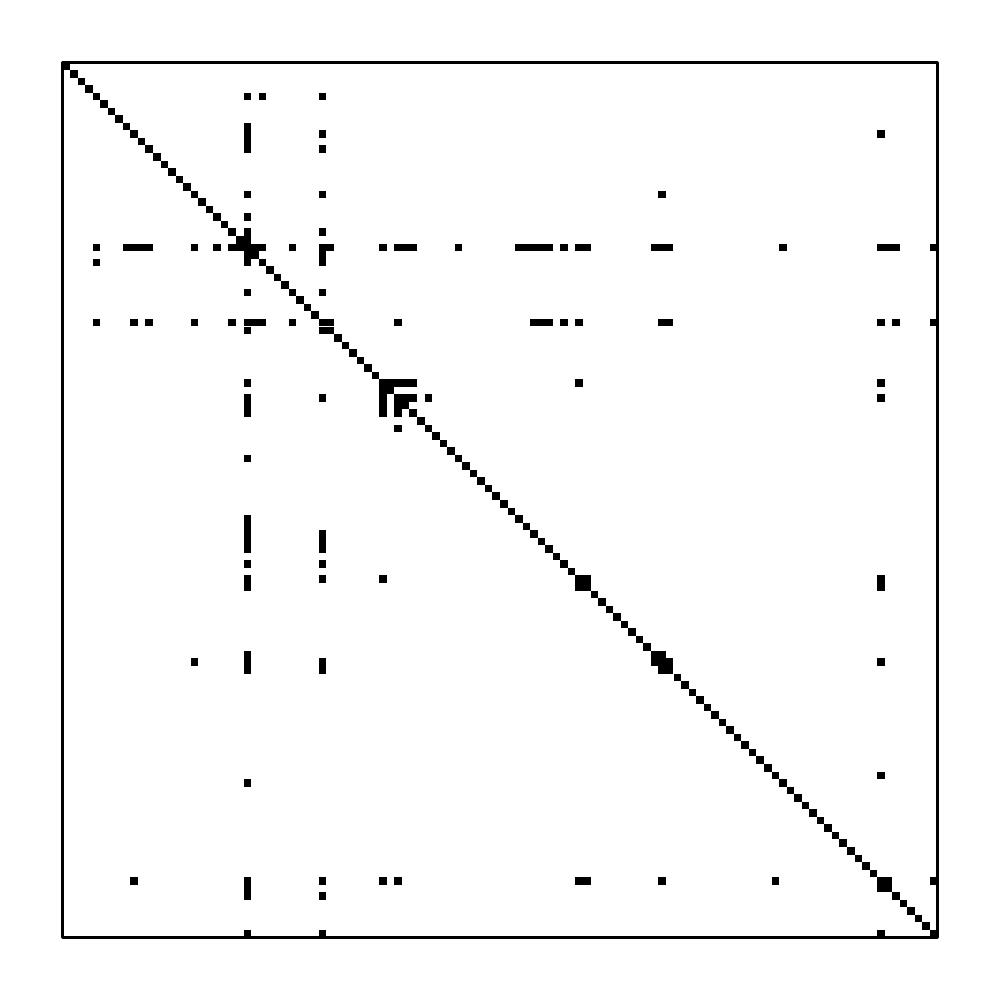}
  \caption{TDC: $\wbSigma_\U$ ($98.85\%$ zeros)} \label{fig:1d}
\end{subfigure}
\begin{subfigure}{0.3\linewidth}
  \includegraphics[width=4.5cm,height=4.5cm]{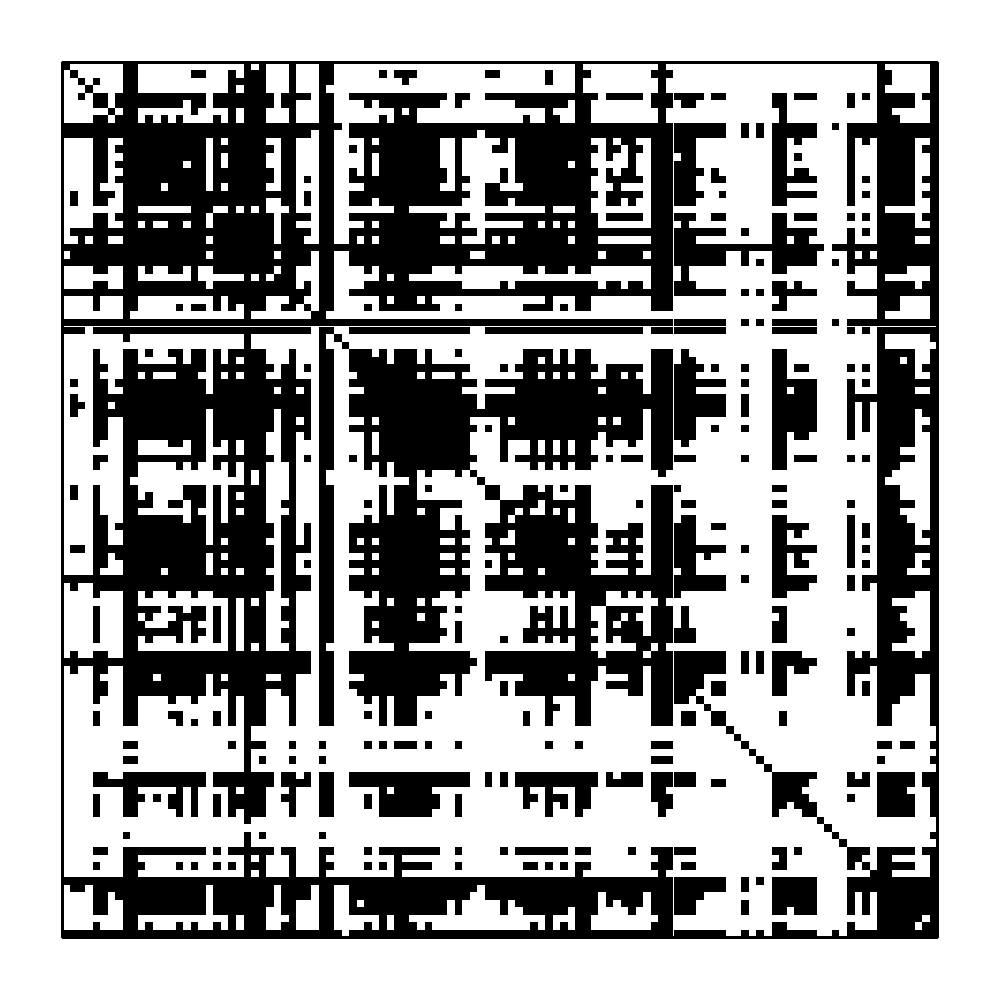}
  \caption{ADHD: $\wbSigma_\U$ ($57.50\%$ zeros)} \label{fig:1e}
  \par\medskip 
  \includegraphics[width=4.5cm,height=4.5cm]{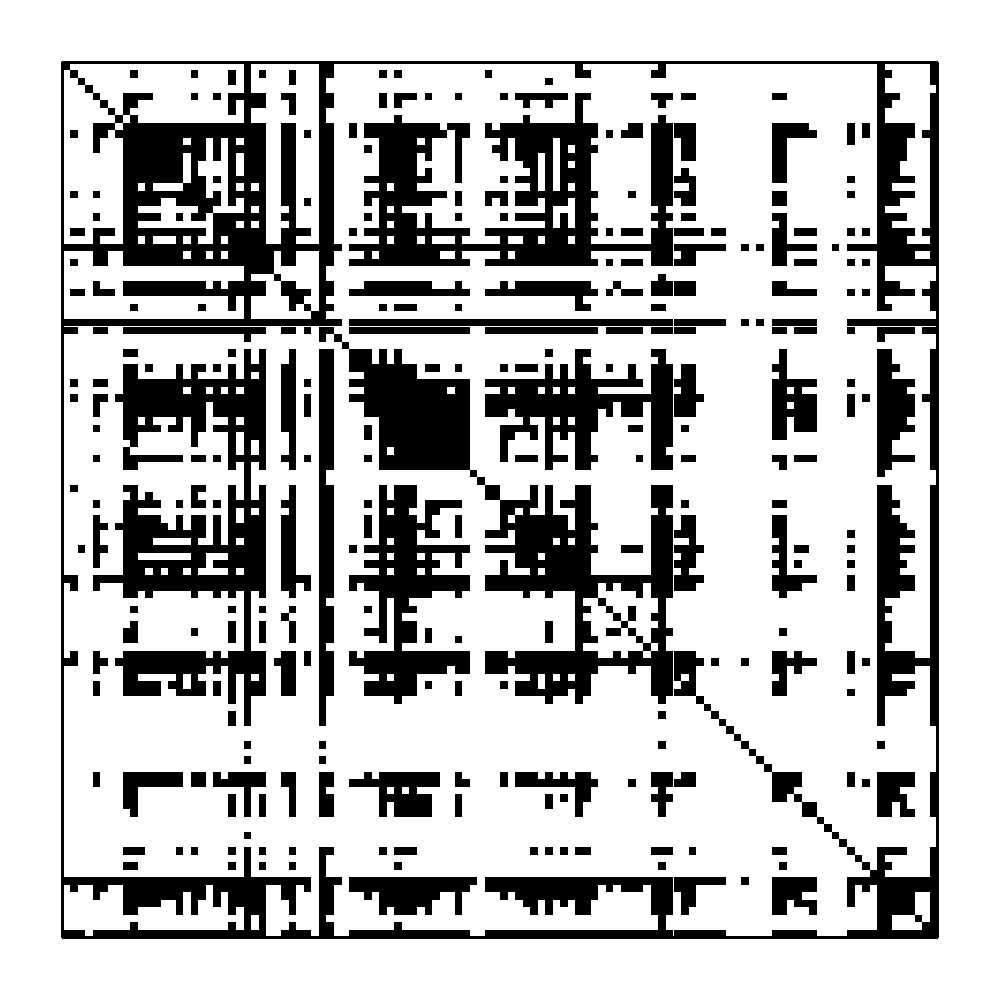}
  \caption{TDC: $\wbSigma_\U$ ($71.24\%$ zeros)} \label{fig:1f}
\end{subfigure}
\centering
\caption{\label{hm_ADHD}{The sparsity structures in $\widehat \bSigma_\A$ and $\widehat \bSigma_\U$ for ADHD and TDC groups: (a)--(d) with the corresponding $\hat\lambda$ selected by fivefold cross-validation using soft functional threosholding rule;  (e)--(f) with the same sparsity levels as those in (a)--(b). Black corresponds to non-zero entries of $\widehat \bSigma_\A$ and $\widehat \bSigma_\U$ (identified edges connecting a subset of ROIs).
}}
\end{figure}

Figure~\ref{hm_ADHD} plots the sparsity patterns in estimated covariance functions corresponding to identified functional connectivity networks. We observe several interesting patterns. 
First, with $\hat \lambda$ selected by the cross-validation, $\widehat \bSigma_\A$ in Fig.~\ref{hm_ADHD}(a)--(b) reveal clear blockwise connectivity structures with two blocks coinciding with the regions of the cerebrum and the cerebellum, while $\widehat \bSigma_\U$ in Fig.~\ref{hm_ADHD}(c)--(d) result in very sparse networks. 
Second, under the same sparsity levels as those of $\wbSigma_\A$ in Fig.~\ref{hm_ADHD}(a)--(b), $\wbSigma_\U$ in Fig.~\ref{hm_ADHD}(e)--(f)
only retain edges related to large marginal-covariance functions but fail to identify some essential within-network connections, e.g., those of the cerebellar region \cite[]{dobromyslin2012} on the bottom right corner.
Third, the ADHD group has increased connections relative to the TDC group, which is in line with the finding in \cite{konrad2010} that
ADHD patients tend to
exhibit abnormal spontaneous functional connectivity patterns.


\subsection{HCP dataset}
\label{sec.hcp}
In this section, we aim to investigate the association between the brain functional connectivity and  fluid intelligence (\textit{gF}), the capacity to solve problems independently of  acquired knowledge
\cite[]{cattell1987}. 
The dataset contains subjects of resting-state fMRI scans and the corresponding \textit{gF} scores, measured by the 24-item Raven’s Progressive Matrices, 
from the Human Connectome Project (HCP). 
We follow many recent proposals based on HCP by modelling signals as multivariate random functions with each ROI representing one random function  \cite[]{zapata2019,lee2021,miao2022}.
We focus our analysis on $n_\low = 73$ subjects with intelligence scores $\textit{gF} \leq 8 $ and $n_\high = 85$ subjects with $\textit{gF} \geq 23 $, and consider $p = 83$ ROIs
of three generally acknowledged modules in neuroscience study \cite[]{finn2015}: the medial frontal (29 ROIs), frontoparietal (34 ROIs) and default mode modules (20 ROIs). For each subject, the BOLD signals at each ROI are collected every 0.72 seconds for a total of $L=1200$ measurement locations (14.4 minutes). 
We first implement the ICA-FIX preprocessed pipeline \cite[]{glasser2013} and a standard band-pass filter at $[0.01,0.08]$ Hz to exclude frequency bands not implicated in resting state functional connectivity \cite[]{biswal1995}. 
Figure~\ref{data_HCP} of Supplementary Material displays examplified trajectories of pre-smoothed data.
The adaptive functional thresholding method is then adopted to estimate the sparse covariance function and therefore the brain networks.

\begin{figure}[!htbp]
\centering
\begin{subfigure}{0.3\linewidth}
  \includegraphics[width=4.5cm,height=4.5cm]{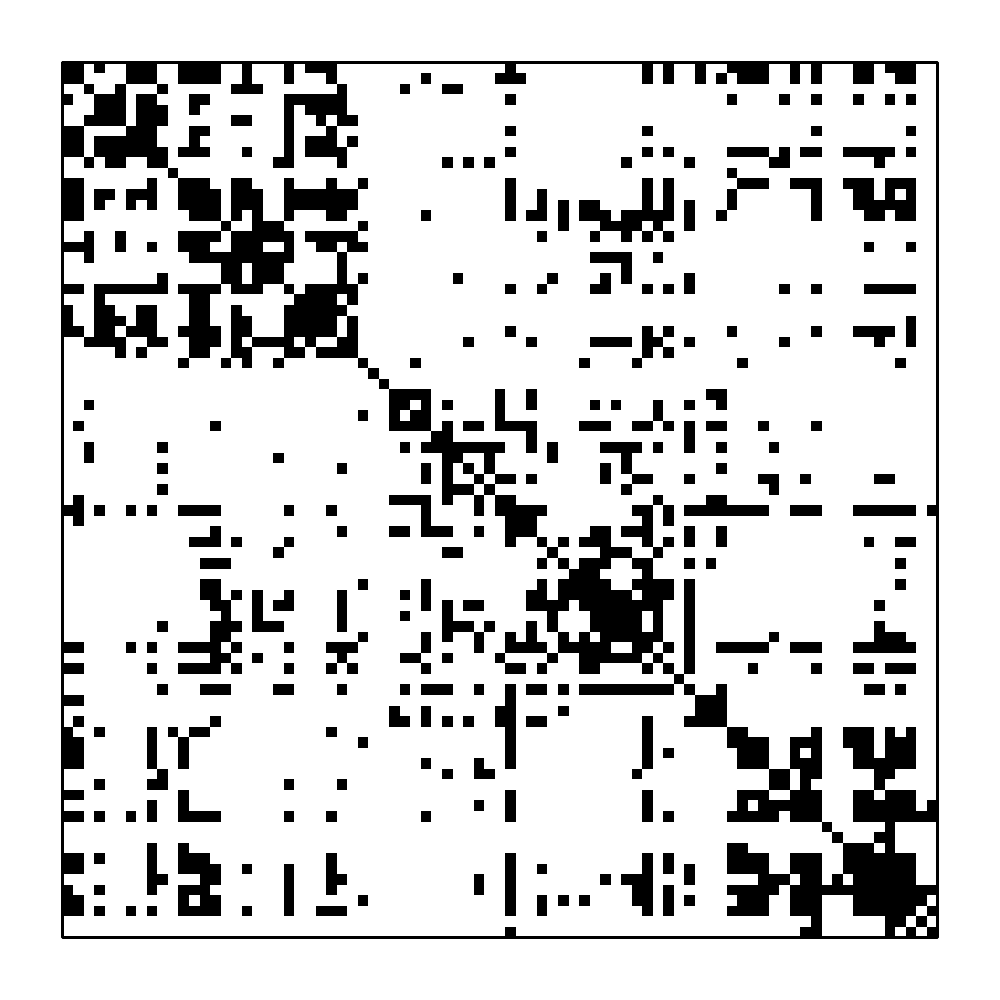}
  \caption{$\textit{gF} \leq 8$ : $\wbSigma_\A$ ($80.42\%$ zeros)} 
  \par\medskip 
  \includegraphics[width=4.5cm,height=4.5cm]{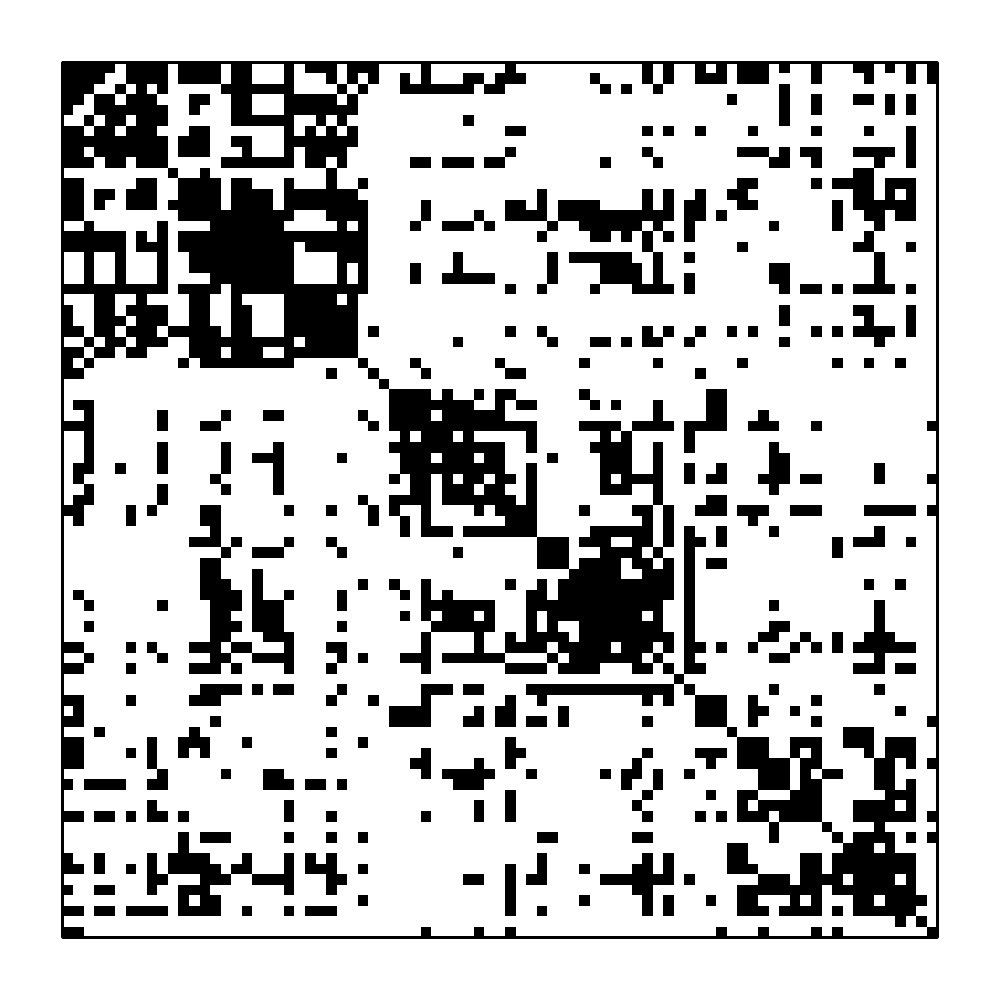}
  \caption{$\textit{gF} \geq 23$: $\wbSigma_\A$ ($72.93\%$ zeros)} 
\end{subfigure}
\begin{subfigure}{0.3\linewidth}
  \includegraphics[width=4.5cm,height=4.5cm]{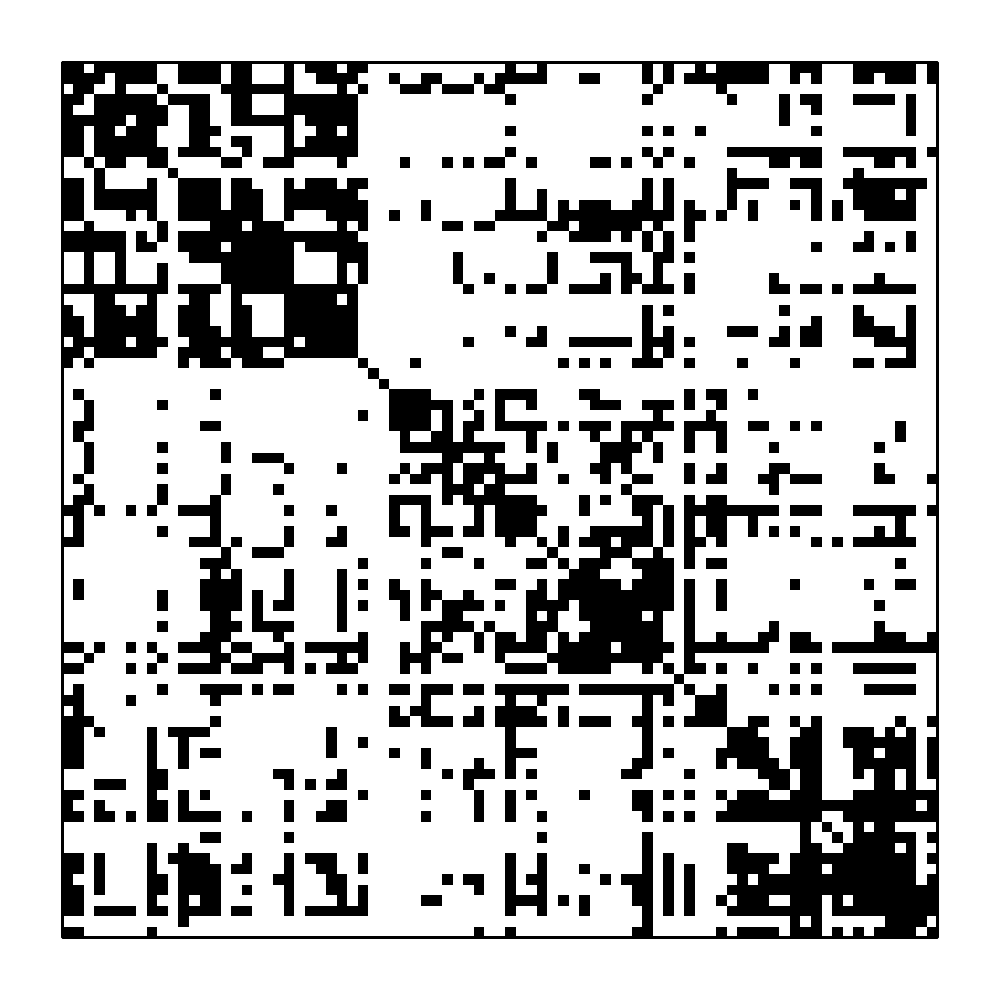}
  \caption{$\textit{gF} \leq 8$: $\wbSigma_\A$ ($70\%$ zeros)} 
  \par\medskip 
  \includegraphics[width=4.5cm,height=4.5cm]{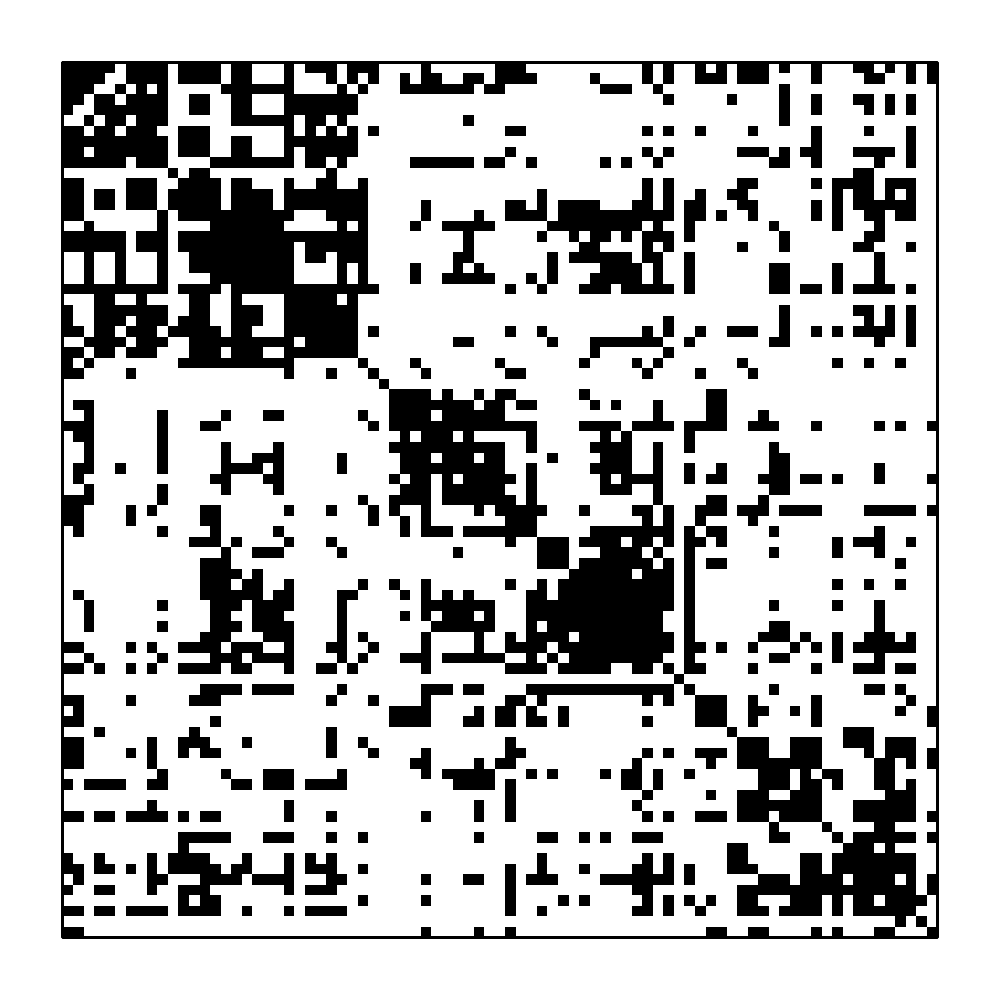}
  \caption{$\textit{gF} \geq 23$: $\wbSigma_\A$ ($70\%$ zeros)} 
\end{subfigure}
\begin{subfigure}{0.3\linewidth}
  \includegraphics[width=4.5cm,height=4.5cm]{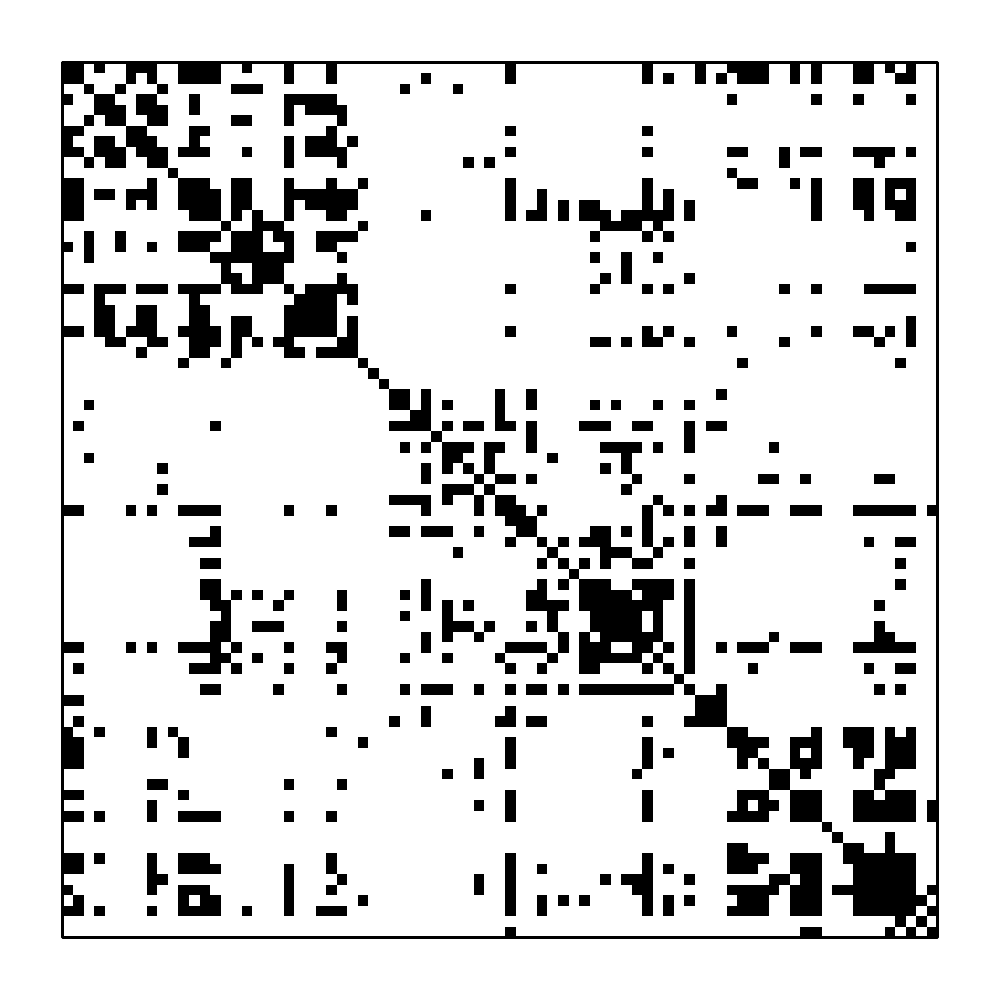}
  \caption{$\textit{gF} \leq 8$: $\wbSigma_\A$ ($85\%$ zeros)} 
  \par\medskip 
  \includegraphics[width=4.5cm,height=4.5cm]{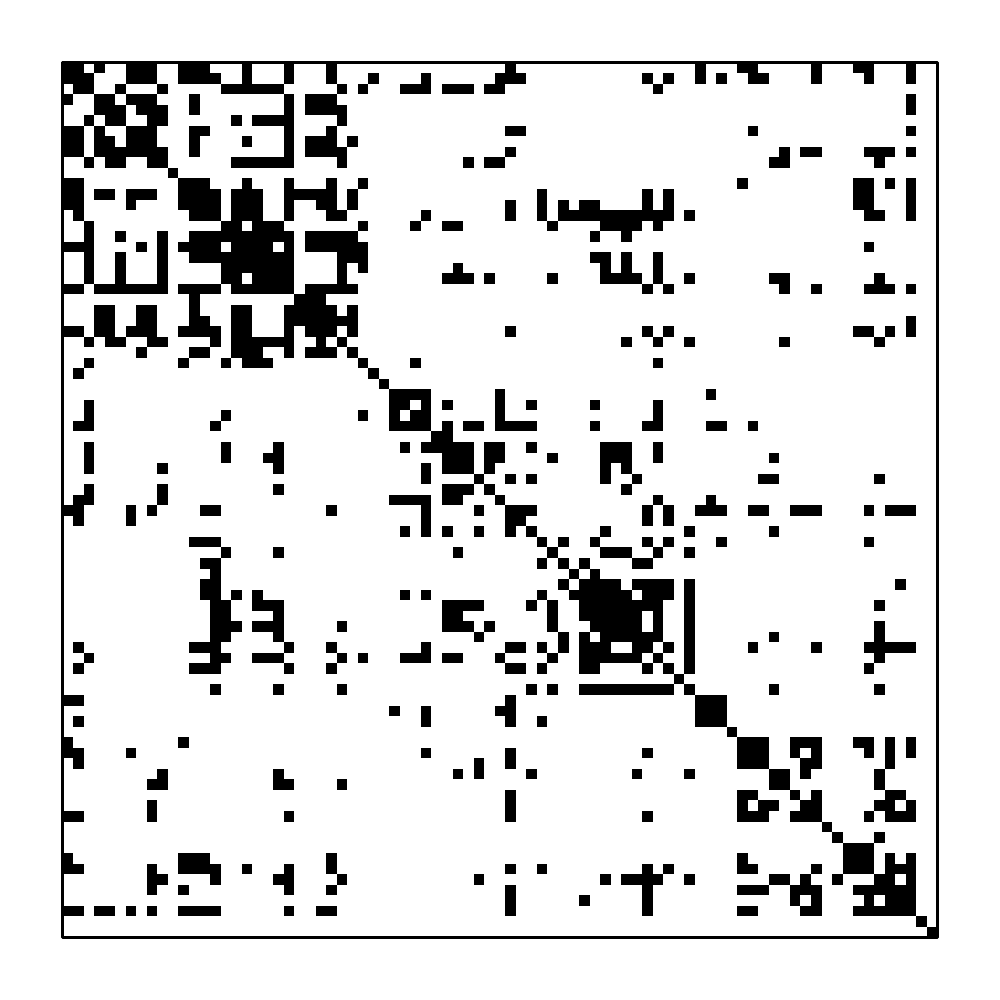}
  \caption{$\textit{gF} \geq 23$: $\wbSigma_\A$ ($85\%$ zeros)} 
\end{subfigure}

\centering
\caption{\label{hm_hcp}{Estimated sparsity structures in $\wbSigma_\A$ using soft functional thresholding rule at fluid intelligence $\textit{gF} \leq 8$ and $\textit{gF} \geq 23$: (a)--(b) with the corresponding $\hat \lambda$ selected by fivefold cross-validation;   (c)--(f)  with the estimated functional sparsity levels set at $70\%$ and $85\%$.
}}
\end{figure}

\begin{figure}[!htbp]
\centering
\begin{subfigure}{0.49\linewidth}
  \includegraphics[width=7cm]{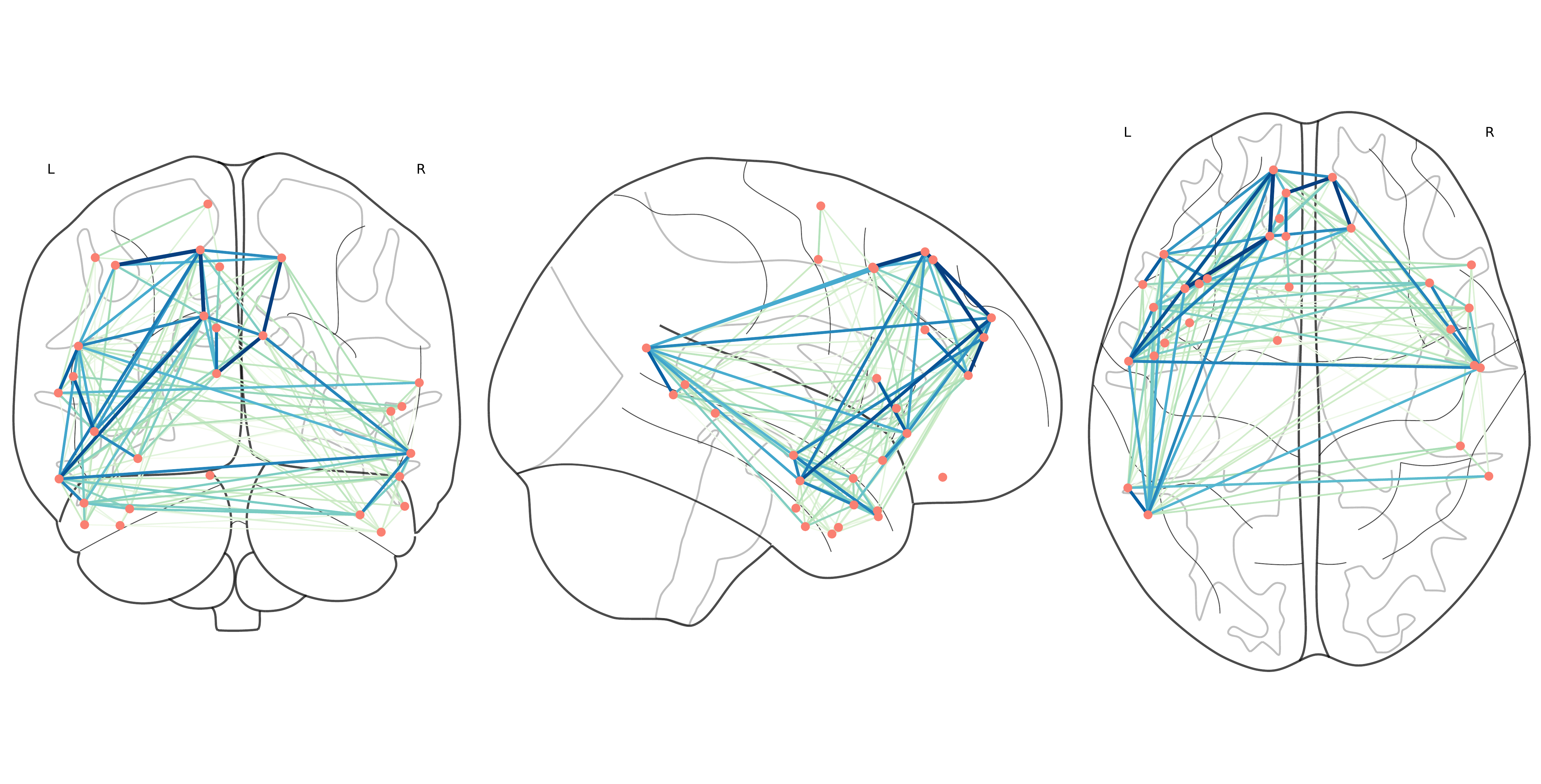}
  \caption{$\textit{gF}  \leq 8$: the medial frontal module in Fig.~\ref{hm_hcp}(e)} 
  \par\medskip 
  \includegraphics[width=7cm]{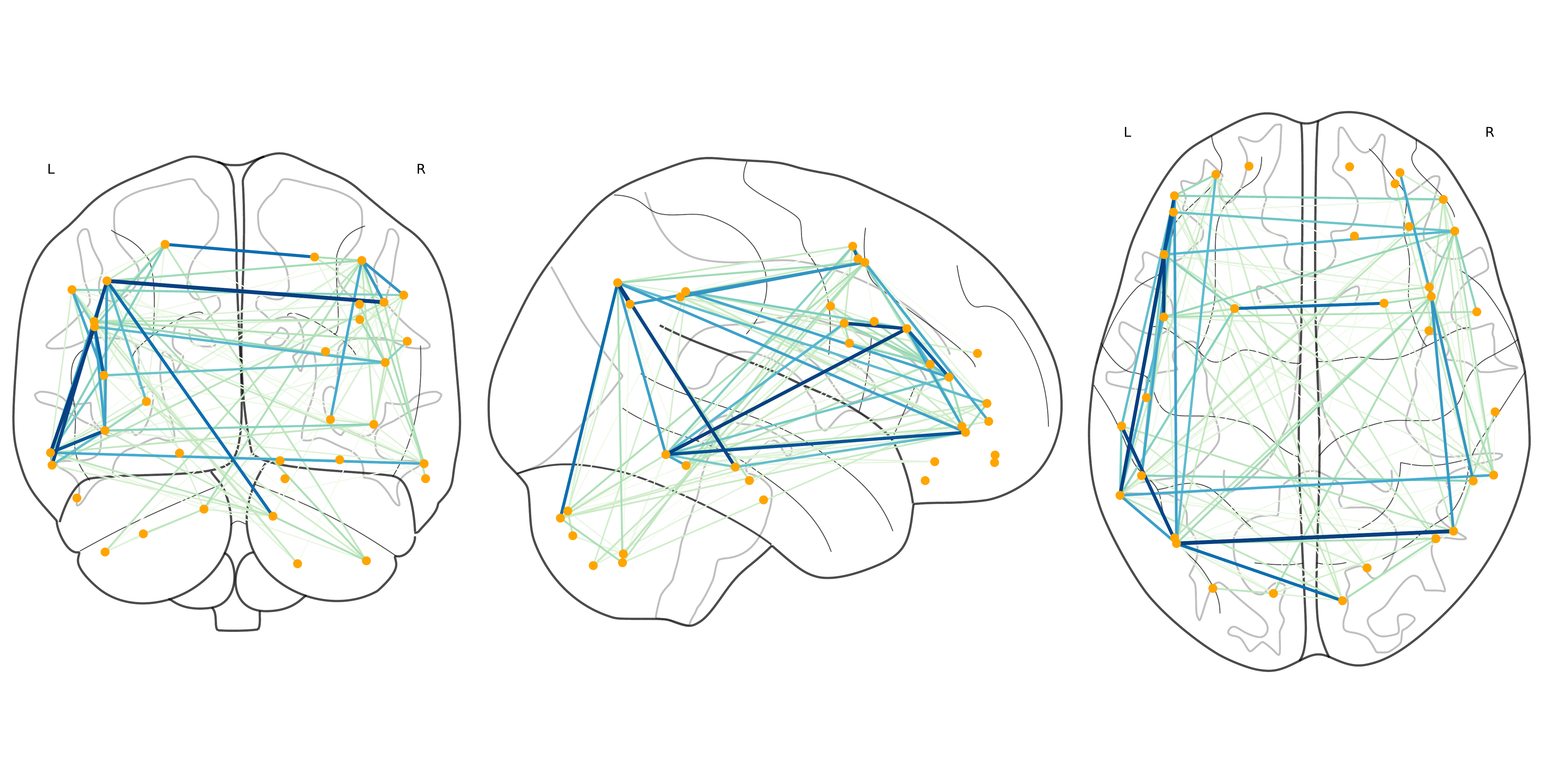}
  \caption{$\textit{gF}  \leq 8$: the frontoparietal module in Fig.~\ref{hm_hcp}(e)} 
  \par\medskip 
  \includegraphics[width=7cm]{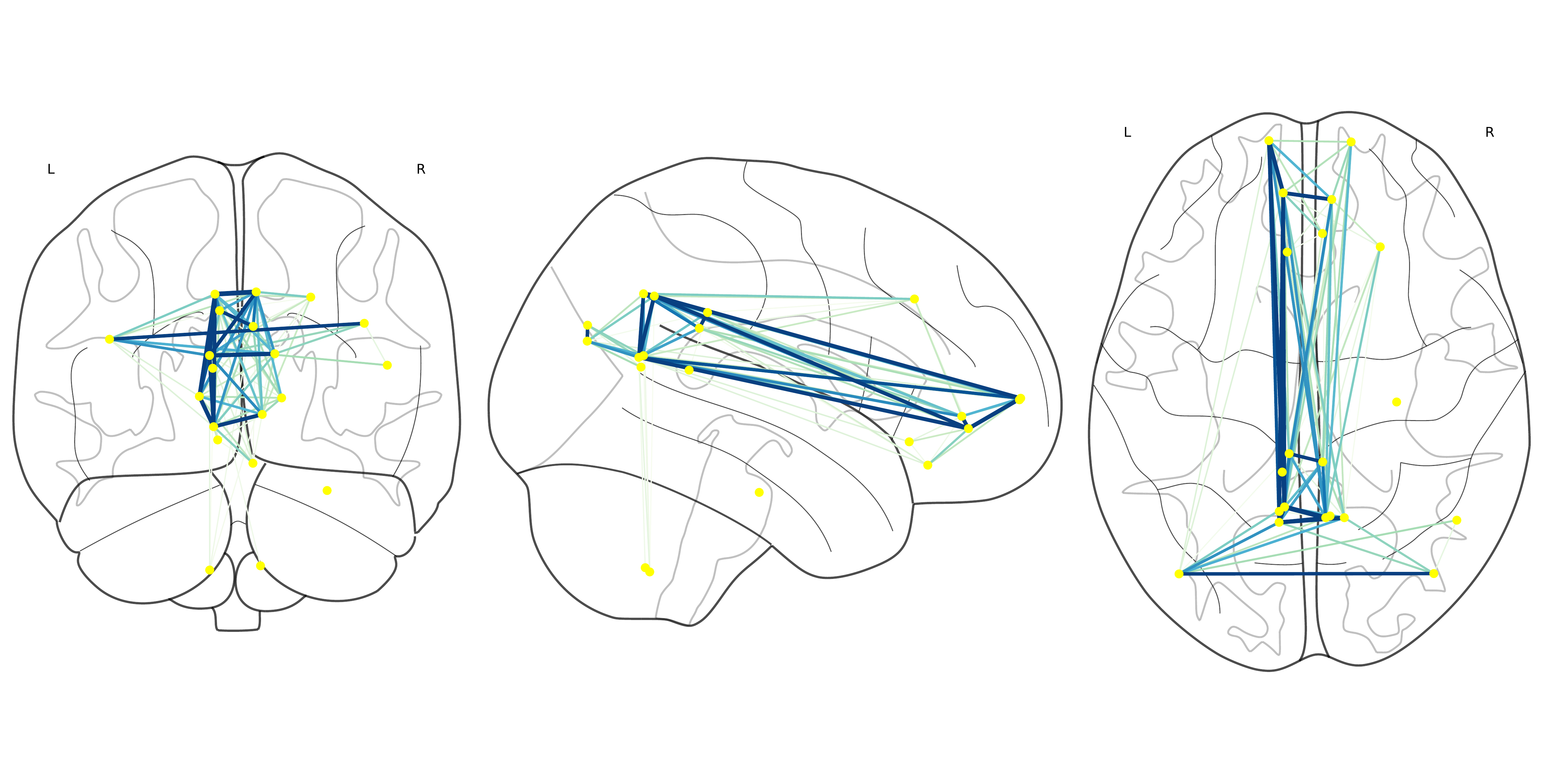}
  \caption{$\textit{gF}  \leq 8$: the default mode module in Fig.~\ref{hm_hcp}(e)} 
\end{subfigure}
\centering
\begin{subfigure}{0.49\linewidth}
  \includegraphics[width=7cm]{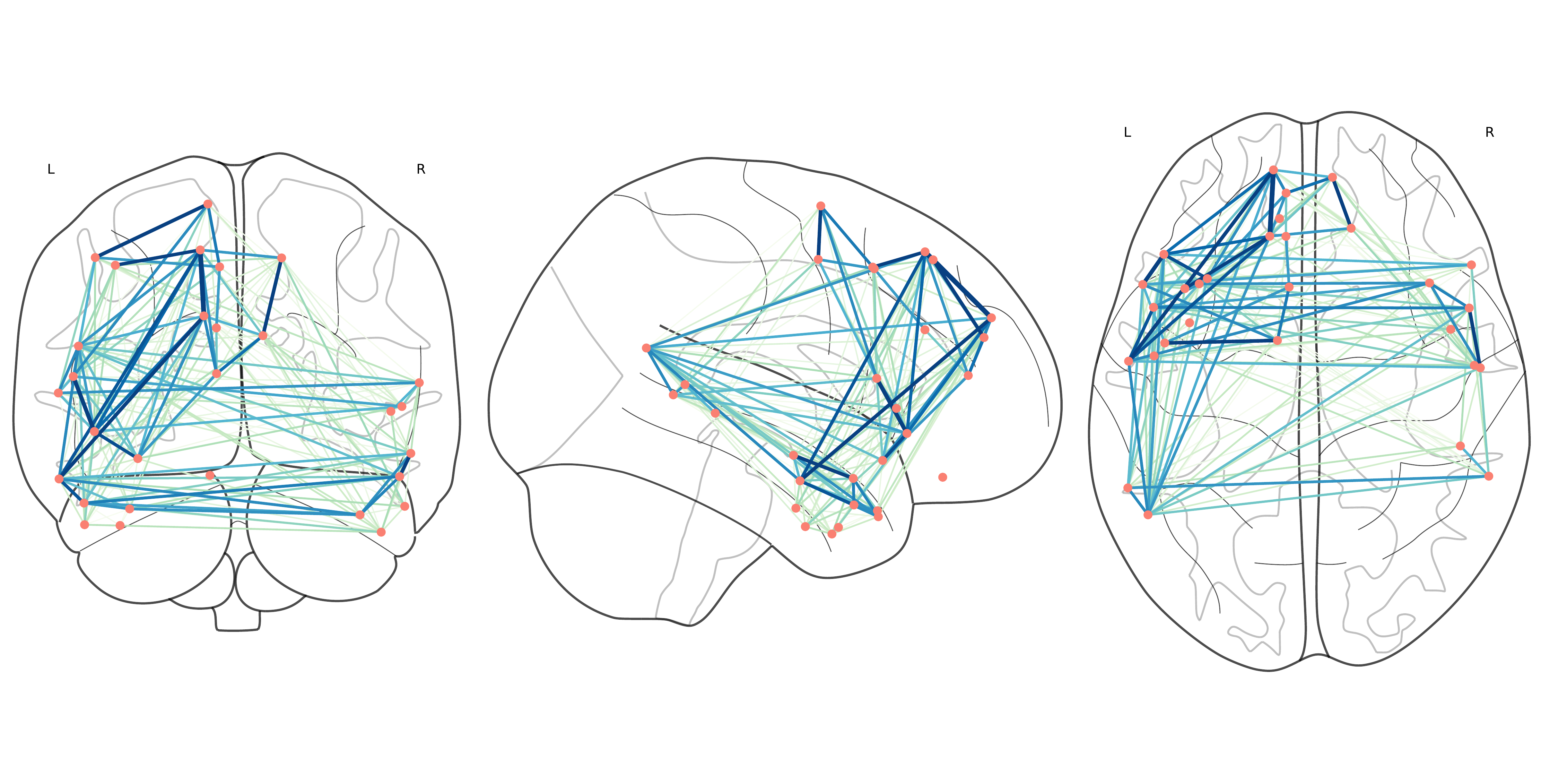}
  \caption{$\textit{gF}   \geq 23$: the medial frontal module in Fig.~\ref{hm_hcp}(f)} 
  \par\medskip 
  \includegraphics[width=7cm]{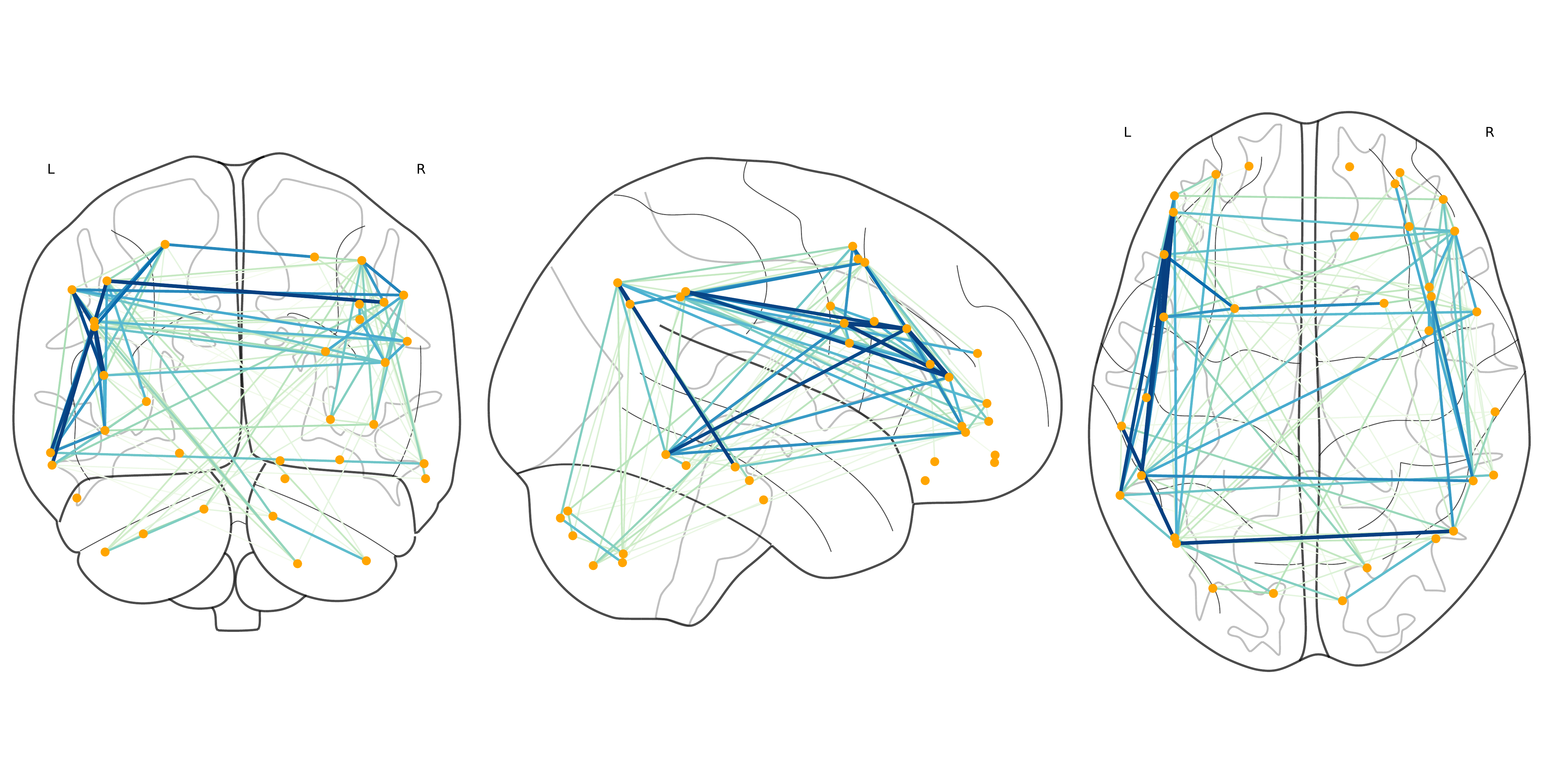}
  \caption{$\textit{gF}   \geq 23$: the frontoparietal module in Fig.~\ref{hm_hcp}(f)} 
  \par\medskip 
  \includegraphics[width=7cm]{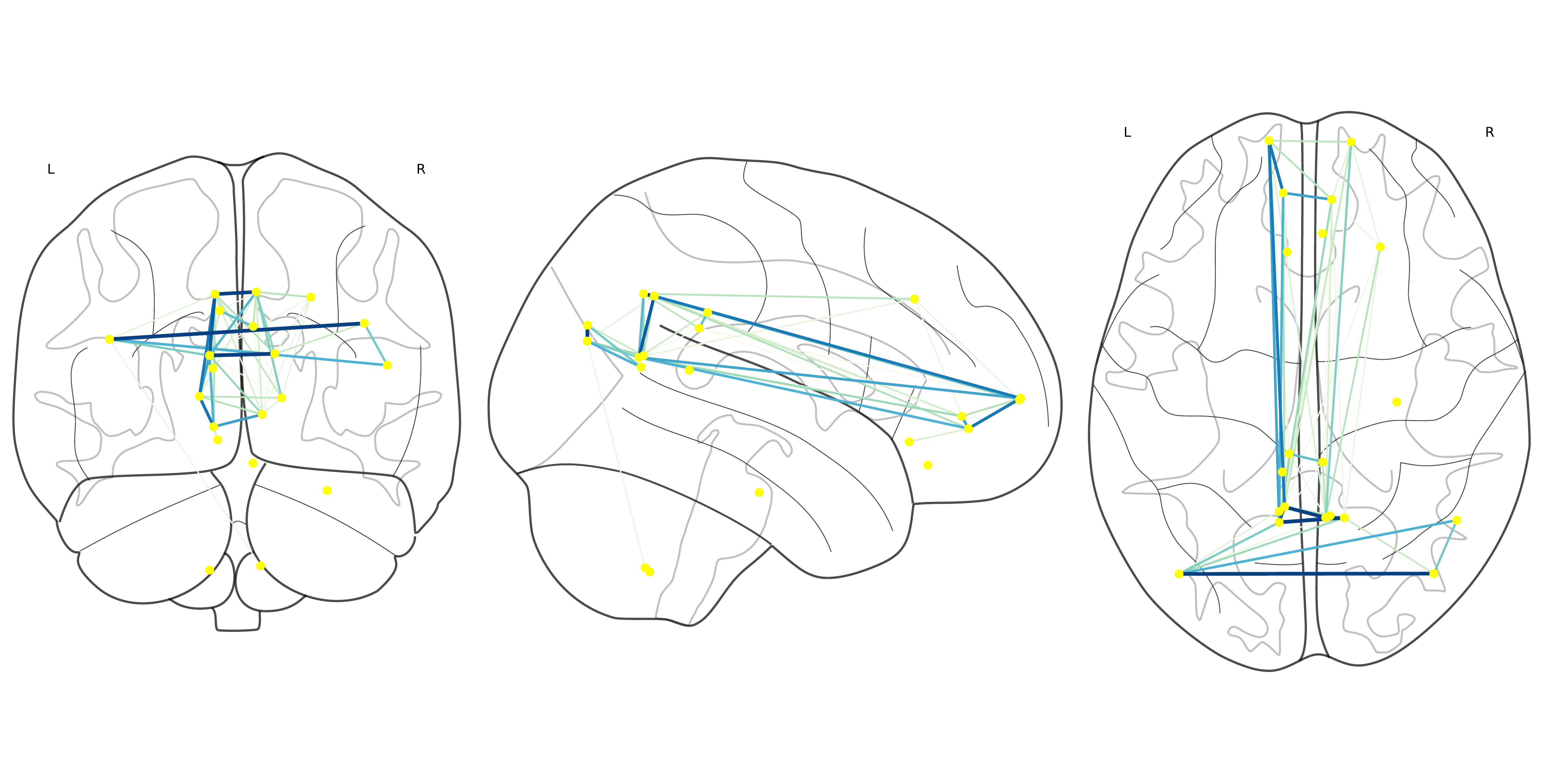}
  \caption{$\textit{gF}   \geq 23$: the default mode module in Fig.~\ref{hm_hcp}(f)} 
\end{subfigure}

\centering
\caption{\label{hcp_network15}{ The connectivity strengths in Fig.~\ref{hm_hcp}(e)--(f) at fluid intelligence $\textit{gF} \leq 8$ and $\textit{gF} \geq 23$.
Salmon, orange and yellow nodes represent the ROIs in the medial frontal, frontoparietal and  default mode modules, respectively. 
The edge color from cyan to blue corresponds to the value of $\|\widehat \Sigma_{jk}^\A\|_\cS/\{\|\widehat \Sigma_{jj}^\A\|_\cS\|\widehat \Sigma_{kk}^\A\|_\cS\}^{1/2}$ from small to large.}}
\end{figure}

The sparsity structures in $\hat \bSigma_\A$ for both groups are displayed in Figure~\ref{hm_hcp}. With $\widehat \lambda$ selected by the cross-validation, 
the network associated with $\wbSigma_\A$ for subjects with $\textit{gF} \geq 23$ is more densely connected than that with $\textit{gF} \leq 8$, as evident from Fig.~\ref{hm_hcp}(a)--(b).
We further set the sparsity level to $70\%$ and $85\%,$ and present the corresponding sparsity patterns in Fig.~\ref{hm_hcp}(c)--(f). The results clearly indicate the existence of three diagonal blocks under all sparsity levels, complying with the identification of the medial frontal, frontoparietal and default mode modules in \cite{finn2015}.
We also implement the universal functional thresholding method. However, compared with $\wbSigma_\A,$ the results of $\wbSigma_\U$ suffer from the heteroscedasticity, as demonstrated in Sections~\ref{sec.sim} and \ref{sec.real.ADHD}, and fail to detect any noticeable block structure, hence we choose not to report them here.
To explore the impact of \textit{gF} on the functional connectivity, we compute the connectivity strength 
using the standardized form $\|\widehat \Sigma_{jk}^\A\|_\cS/\{\|\widehat \Sigma_{jj}^\A\|_\cS\|\widehat \Sigma_{kk}^\A\|_\cS\}^{1/2}$ for $ j, k = 1 \dots, p.$ Interestingly, we observe from Figure~\ref{hcp_network15} that subjects with $\textit{gF} \geq 23$ tend to have enhanced brain connectivity in the medial frontal and frontoparietal modules, while the connectivity strength in the default mode module declines. 
This agrees with existing neuroscience literature reporting a strong positive association between  intelligence score and the medial frontal/frontoparietal functional connectivity in the resting state \citep{van2009,finn2015}, 
and lends support to the conclusion that lower default mode module activity is associated with better cognitive performance \citep{anticevic2012}.


\appendix

\spacingset{1}
\bibliography{paperbib}

@article{cai2011,
 author = {Tony Cai and Weidong Liu},
title = {Adaptive Thresholding for Sparse Covariance Matrix Estimation},
journal = {Journal of the American Statistical Association},
volume = {106},
pages = {672-684},
year  = {2011},
 }

@article{guo2022,
	author               = {Guo, S. and Qiao, X. and Wang, Q.},
	journal              = {arXiv:2112.13651v2},
	pages                = {},
	title                = {Factor modelling for high-dimensional functional time series},
	volume               = {},
	year                 = {2022},
}

@article{vu2013,
  author  = {Vu, Vincent Q. and Lei, Jing},
  journal = {The Annals of Statistics},
  pages   = {2905--2947},
  title   = {Minimax sparse principal subspace estimation in high dimensions},
  volume  = {41},
  year    = {2013}
}

@article{chang2017,
	author               = {Chang, J. and Guo, B. and Yao, Q.},
	journal              = {The Annals of Statistics},
	pages                = {2094-2124},
	title                = {Principal component analysis for second-order stationary vector time series},
	volume               = {46},
	year                 = {2018},
}

@article{qiao2019a,
	author               = {Qiao, X. and Guo, S. and James, G.},
	journal              = {Journal of the American Statistical Association},
	title                = {Functional graphical models},
	volume={114},

  pages={211--222},
  year={2019},
}

@article {qiao2020,
    AUTHOR = {Qiao, Xinghao and Qian, Cheng and James, Gareth M. and Guo,
              Shaojun},
     TITLE = {Doubly functional graphical models in high dimensions},
   JOURNAL = {Biometrika},
  FJOURNAL = {Biometrika},
    VOLUME = {107},
      YEAR = {2020},
     PAGES = {415--431},
      ISSN = {0006-3444},
   MRCLASS = {62R10 (62H22 92B15)},
  MRNUMBER = {4108937},
       DOI = {10.1093/biomet/asz072},
}

@article {bickel2008,
    AUTHOR = {Bickel, Peter J. and Levina, Elizaveta},
     TITLE = {Covariance regularization by thresholding},
  JOURNAL = {The Annals of Statistics},
    VOLUME = {36},
      YEAR = {2008},
     PAGES = {2577--2604},
      ISSN = {0090-5364},
   MRCLASS = {62H12 (62F12 62G09)},
  MRNUMBER = {2485008},
MRREVIEWER = {M. Hu\v{s}kov\'{a}},
       DOI = {10.1214/08-AOS600},
}

@article {rothman2009,
    AUTHOR = {Rothman, Adam J. and Levina, Elizaveta and Zhu, Ji},
     TITLE = {Generalized thresholding of large covariance matrices},
  JOURNAL = {Journal of the American Statistical Association},
    VOLUME = {104},
      YEAR = {2009},
     PAGES = {177--186},
      ISSN = {0162-1459},
   MRCLASS = {62H20 (62J07)},
  MRNUMBER = {2504372},
MRREVIEWER = {Bao Xue Zhang},
       DOI = {10.1198/jasa.2009.0101},
}

@article{konrad2010,
  title={{Is the ADHD brain wired differently? A review on structural and functional connectivity in attention deficit hyperactivity disorder}},
  author={Konrad, Kerstin and Eickhoff, Simon B},
  journal={Human Brain Mapping},
  volume={31},
  pages={904--916},
  year={2010},
  publisher={Wiley Online Library}
}

@article{bellec2017,
  title={{The neuro bureau ADHD-200 preprocessed repository}},
  author={Bellec, Pierre and Chu, Carlton and Chouinard-Decorte, Francois and Benhajali, Yassine and Margulies, Daniel S and Craddock, R Cameron},
  journal={Neuroimage},
  volume={144},
  pages={275--286},
  year={2017},
  publisher={Elsevier}
}

@article{tzourio2002,
  title={{Automated anatomical labeling of activations in SPM using a macroscopic anatomical parcellation of the MNI MRI single-subject brain}},
  author={Tzourio-Mazoyer, Nathalie and Landeau, Brigitte and Papathanassiou, Dimitri and Crivello, Fabrice and Etard, Olivier and Delcroix, Nicolas and Mazoyer, Bernard and Joliot, Marc},
  journal={Neuroimage},
  volume={15},
  pages={273--289},
  year={2002},
  publisher={Academic Press}
}

@article{finn2015,
  title={Functional connectome fingerprinting: identifying individuals using patterns of brain connectivity},
  author={Finn, Emily S and Shen, Xilin and Scheinost, Dustin and Rosenberg, Monica D and Huang, Jessica and Chun, Marvin M and Papademetris, Xenophon and Constable, R Todd},
  journal={Nature Neuroscience},
  volume={18},
  pages={1664--1671},
  year={2015},
  publisher={Nature Publishing Group}
}

@article{glasser2013,
  title={The minimal preprocessing pipelines for the Human Connectome Project},
  author={Glasser, Matthew F and Sotiropoulos, Stamatios N and Wilson, J Anthony and Coalson, Timothy S and Fischl, Bruce and Andersson, Jesper L and Xu, Junqian and Jbabdi, Saad and Webster, Matthew and Polimeni, Jonathan R and others},
  journal={Neuroimage},
  volume={80},
  pages={105--124},
  year={2013},
  publisher={Elsevier}
}

@article{fan2001,
  title={Variable selection via nonconcave penalized likelihood and its oracle properties},
  author={Fan, Jianqing and Li, Runze},
  journal={Journal of the American Statistical Association},
  volume={96},
  pages={1348--1360},
  year={2001},
  publisher={Taylor \& Francis}
}

@article{zou2006,
  title={The adaptive lasso and its oracle properties},
  author={Zou, Hui},
  journal={Journal of the American Statistical Association},
  volume={101},
  pages={1418--1429},
  year={2006},
  publisher={Taylor \& Francis}
}

@book {Kosorok2008,
    AUTHOR = {Kosorok, Michael R.},
     TITLE = {Introduction to empirical processes and semiparametric
              inference},
    SERIES = {Springer Series in Statistics},
 PUBLISHER = {Springer, New York},
      YEAR = {2008},
     PAGES = {xiv+483},
      ISBN = {978-0-387-74977-8},
   MRCLASS = {62-02 (60B05 60B12 60F17 62Gxx)},
  MRNUMBER = {2724368},
MRREVIEWER = {Gutti J. Babu},
       DOI = {10.1007/978-0-387-74978-5},
}

@article{li2018,
author = {Bing Li and Eftychia Solea},
title = {A Nonparametric Graphical Model for Functional Data With Application to Brain Networks Based on {fMRI}},
journal = {Journal of the American Statistical Association},
volume = {113},
pages = {1637-1655},
year  = {2018},
publisher = {Taylor & Francis},
}

@article{park2021,
author = {Park, Juhyun and Ahn, Jeongyoun and Jeon, Yongho},
    title = "{Sparse functional linear discriminant analysis}",
    journal = {Biometrika},
    volume = {109},
    pages = {209-226},
    year = {2021},
}

@article{happ2018,
    author = {Clara Happ and Sonja Greven},
    title = {Multivariate Functional Principal Component Analysis for
      Data Observed on Different (Dimensional) Domains},
    journal = {Journal of the American Statistical Association},
    year = {2018},
    volume = {113},
    issue = {522},
    pages = {649-659},
  }

@article{chen2016,
author = {Ziqi Chen and Chenlei Leng},
title = {Dynamic Covariance Models},
journal = {Journal of the American Statistical Association},
volume = {111},
pages = {1196-1207},
year  = {2016},
}

@article{yuan2006,
 author               = {Yuan, M. and Lin, Y.},
 journal              = {Journal of the Royal Statistical Society: Series B},
 pages                = {49-67},
 title                = {Model selection and estimation in regression with grouped variables},
 volume               = {68},
 year                 = {2006},
 }

@book{cattell1987,
  title={Intelligence: Its structure, growth and action},
  author={Cattell, Raymond Bernard},
  year={1987},
  publisher={Elsevier}
}

@article{dobromyslin2012,
  title={Distinct functional networks within the cerebellum and their relation to cortical systems assessed with independent component analysis},
  author={Dobromyslin, Vitaly I and Salat, David H and Fortier, Catherine B and Leritz, Elizabeth C and Beckmann, Christian F and Milberg, William P and McGlinchey, Regina E},
  journal={Neuroimage},
  volume={60},
  pages={2073--2085},
  year={2012},
  publisher={Elsevier}
}

@article{avella2018,
  title={Robust estimation of high-dimensional covariance and precision matrices},
  author={Avella-Medina, Marco and Battey, Heather S and Fan, Jianqing and Li, Quefeng},
  journal={Biometrika},
  volume={105},
  pages={271--284},
  year={2018},
  publisher={Oxford University Press}
}

@article{biswal1995,
  title={Functional connectivity in the motor cortex of resting human brain using echo-planar {MRI}},
  author={Biswal, Bharat and Zerrin Yetkin, F and Haughton, Victor M and Hyde, James S},
  journal={Magnetic resonance in medicine},
  volume={34},
  pages={537--541},
  year={1995},
  publisher={Wiley Online Library}
}

@article{chang2010,
  title={Time--frequency dynamics of resting-state brain connectivity measured with {fMRI}},
  author={Chang, Catie and Glover, Gary H},
  journal={Neuroimage},
  volume={50},
  pages={81--98},
  year={2010},
  publisher={Elsevier}
}

@article{rogers2007,
  title={Assessing functional connectivity in the human brain by {fMRI}},
  author={Rogers, Baxter P and Morgan, Victoria L and Newton, Allen T and Gore, John C},
  journal={Magnetic resonance imaging},
  volume={25},
  pages={1347--1357},
  year={2007},
  publisher={Elsevier}
}

@article{lotte2018,
  title={A review of classification algorithms for EEG-based brain--computer interfaces: a 10 year update},
  author={Lotte, Fabien and Bougrain, Laurent and Cichocki, Andrzej and Clerc, Maureen and Congedo, Marco and Rakotomamonjy, Alain and Yger, Florian},
  journal={Journal of Neural Engineering},
  volume={15},
  pages={031005},
  year={2018},
  publisher={IOP Publishing}
}

@article{wang2020,
  title = {Nonparametric estimation of large covariance matrices with conditional sparsity},
journal = {Journal of Econometrics},
volume = {223},
pages = {53-72},
year = {2021},
author = {Hanchao Wang and Bin Peng and Degui Li and Chenlei Leng},
}

@article {fan2016,
    AUTHOR = {Fan, Yingying and Lv, Jinchi},
     TITLE = {Innovated scalable efficient estimation in ultra-large
              {G}aussian graphical models},
JOURNAL = {The Annals of Statistics},
    VOLUME = {44},
      YEAR = {2016},
     PAGES = {2098--2126},
      ISSN = {0090-5364},
   MRCLASS = {62H12 (62F12 62J05)},
  MRNUMBER = {3546445},
       DOI = {10.1214/15-AOS1416},
}

@article{storey2005,
  title={Significance analysis of time course microarray experiments},
  author={Storey, John D and Xiao, Wenzhong and Leek, Jeffrey T and Tompkins, Ronald G and Davis, Ronald W},
  journal={Proceedings of the National Academy of Sciences},
  volume={102},
  pages={12837--12842},
  year={2005},
  publisher={National Acad Sciences}
}

@article{kong2016,
  title={Partially functional linear regression in high dimensions},
  author={Kong, Dehan and Xue, Kaijie and Yao, Fang and Zhang, Hao H},
  journal={Biometrika},
  volume={103},
  pages={147--159},
  year={2016},
  publisher={Oxford University Press}
}

@article {Zhang2016,
    AUTHOR = {Zhang, Xiaoke and Wang, Jane-Ling},
     TITLE = {From sparse to dense functional data and beyond},
  
 JOURNAL = {The Annals of Statistics},
    VOLUME = {44},
      YEAR = {2016},
     PAGES = {2281--2321},
      ISSN = {0090-5364},
   MRCLASS = {62G20 (62G05 62G08)},
  MRNUMBER = {3546451},
       DOI = {10.1214/16-AOS1446},
}

@article{fan1994,
  title={Fast implementations of nonparametric curve estimators},
  author={Fan, Jianqing and Marron, James S},
  journal={Journal of Computational and Graphical Statistics},
  volume={3},

  pages={35--56},
  year={1994},
  publisher={Taylor \& Francis}
}

@article {Zhang2007,
    AUTHOR = {Zhang, Jin-Ting and Chen, Jianwei},
     TITLE = {Statistical inferences for functional data},

  JOURNAL = {The Annals of Statistics},
    VOLUME = {35},
      YEAR = {2007},

     PAGES = {1052--1079},
      ISSN = {0090-5364},
   MRCLASS = {62G07 (62G10 62J12)},
  MRNUMBER = {2341698},
       DOI = {10.1214/009053606000001505},
}

@article {Yao2005a,
    AUTHOR = {Yao, Fang and M\"{u}ller, Hans-Georg and Wang, Jane-Ling},
     TITLE = {Functional data analysis for sparse longitudinal data},
   JOURNAL = {Journal of the American Statistical Association},
    VOLUME = {100},
      YEAR = {2005},

     PAGES = {577--590},
      ISSN = {0162-1459},
   MRCLASS = {62H25 (62G05)},
  MRNUMBER = {2160561},
MRREVIEWER = {M. Riedel},
       DOI = {10.1198/016214504000001745},
}

@article{miao2022,
  title={A Wavelet-Based Independence Test for Functional Data with an Application to {MEG} Functional Connectivity},
  author={Miao, Rui and Zhang, Xiaoke and Wong, Raymond KW},
  journal={Journal of the American Statistical Association, in press},
  pages={},
  year={2022},
  publisher={Taylor \& Francis}
}

@article{lee2021,
  title={Conditional functional graphical models},
  author={Lee, Kuang-Yao and Ji, Dingjue and Li, Lexin and Constable, Todd and Zhao, Hongyu},
  journal={Journal of the American Statistical Association, in press},
  pages={},
  year={2021},
  publisher={Taylor \& Francis}
}

@article{zapata2019,
  author = {Zapata, J and Oh, S Y and Petersen, A},
    title = "{Partial separability and functional graphical models for multivariate Gaussian processes}",
    journal = {Biometrika, in press},
    year = {2021},
    volume = {}
}

@article{van2009,
  title={Efficiency of functional brain networks and intellectual performance},
  author={Van Den Heuvel, Martijn P and Stam, Cornelis J and Kahn, Ren{\'e} S and Pol, Hilleke E Hulshoff},
  journal={Journal of Neuroscience},
  volume={29},

  pages={7619--7624},
  year={2009},
  publisher={Soc Neuroscience}
}

@article{anticevic2012,
  title={The role of default network deactivation in cognition and disease},
  author={Anticevic, Alan and Cole, Michael W and Murray, John D and Corlett, Philip R and Wang, Xiao-Jing and Krystal, John H},
  journal={Trends in cognitive sciences},
  volume={16},

  pages={584--592},
  year={2012},
  publisher={Elsevier}
}
\bibliographystyle{dcu}

\newpage
\spacingset{1.7}
\begin{center}
	{\noindent \bf \large Supplementary material to ``Adaptive functional thresholding for sparse covariance function estimation in high dimensions"}\\
\end{center}
 \begin{center}
 	{\noindent Qin Fang, Shaojun Guo and Xinghao Qiao}
 \end{center}
\bigskip

\setcounter{page}{1}
\setcounter{section}{0}
\renewcommand\thesection{\Alph{section}}
\setcounter{lemma}{0}
\renewcommand{\thelemma}{\Alph{section}\arabic{lemma}}
\setcounter{equation}{0}
\renewcommand{\theequation}{S.\arabic{equation}}


This supplementary material contains  the technical proofs for the fully observed functional scenario in Section~\ref{supp.sec_pr}, derivations of functional thresholding rules in Section~\ref{supp.sec_con}, additional methodological details and technical proofs for the partially observed functional scenario in Section~\ref{supp.sec_partial} and additional empirical results in Section~\ref{supp.sec_emp}.

\section{Technical proofs}
\label{supp.sec_pr}



 
Before stating the regularity conditions, we make some notation. For a function $Z\in \mathbb{S},$ define 
$\|Z\|_{\infty} = {\sup}_{u,v\in \cU} |Z(u,v)|.$ 
For two sequences of real processes $\{a_n(u), u \in \cU\}$
and $\{b_n(u), u \in \cU\}$,  we write $a_n(u) \lesssim b_n(u)$ if there exists some constant $c$ such
that $|a_n(u)| \le  c|b_n(u)|$ holds for all $n$ and $u \in \cU.$ 
Without loss of generality, in the following we assume that $\mathbb{E}\{X_{ij}(u)\} \equiv 0$ and  both estimators $\widehat{\Sigma}_{jk}(u,v)$ and $\widehat{\Theta}_{jk}(u,v)$ are defined as
\begin{equation*}
   \widehat{\Sigma}_{jk}(u,v) = \frac{1}{n}\sum_{i=1}^n X_{ij}(u) X_{ik}(v)~\mbox{and}~  \widehat{\Theta}_{jk}(u,v) = \frac{1}{n}\sum_{i=1}^n X_{ij}(u)^2 X_{ik}(v)^2 - \widehat{\Sigma}_{jk}(u,v)^2,
\end{equation*}
respectively.

\begin{lemma} \label{lm_theta}
Suppose that Conditions~\ref{con_psi2_process}--\ref{con_rate} hold. Then for any $M>0,$ there exists some constant $\rho_1 >0$  such that 
\begin{equation} \nonumber
    P\left \{ \max_{j,k} \left\| \frac{\widehat \Theta_{jk}  - \Theta_{jk}}{\Theta_{jk}}\right\|_{\infty} \geq  \rho_1 \frac{\log^2 p}{n^{1/2}} \right\} = O(p^{-M}).
\end{equation}
\end{lemma}

\textbf{Proof}.  Denote $\widetilde{\Theta}_{jk}(u,v) = \mathbb{E}\{X_{ij}(u)^2X_{ik}(v)^2\}.$ We decompose $\widehat \Theta_{jk} (u,v) - \Theta_{jk}(u,v)$ as 
\begin{equation}\nonumber
\begin{split}
  &\widehat \Theta_{jk} (u,v) - \Theta_{jk}(u,v) 
  \\=~&\Sigma_{jk}(u,v)^2 - \widehat \Sigma_{jk}(u,v)^2 + \frac{1}{n} \sum_{i = 1}^n \Big\{X_{ij}(u)^2X_{ik}(v)^2 - \widetilde{\Theta}_{jk}(u,v)\Big\}. 
  \end{split}
\end{equation}
By Condition~\ref{con_variance_function}, $\Theta_{jk}(u,v) \ge \tau \sigma_j(u) \sigma_k(v)$ for each $j,k = 1,\ldots,p.$ Hence, 
\begin{equation}\nonumber
\begin{split}
  &\left|\frac{\widehat \Theta_{jk} (u,v) - \Theta_{jk}(u,v)}{\Theta_{jk}(u,v)}\right| 
  \\ \le 
  &\left|\frac{\Sigma_{jk}(u,v)^2 - \widehat \Sigma_{jk}(u,v)^2 }{\tau \sigma_j(u) \sigma_k(v)}\right|+ 
  \left|\frac{1}{n} \sum_{i = 1}^n \frac{X_{ij}(u)^2X_{ik}(v)^2 - \widetilde{\Theta}_{jk}(u,v)}{\tau\sigma_j(u) \sigma_k(v)}\right|
  \\= &~ H_{jk}^{(1)}(u,v)  + H_{jk}^{(2)}(u,v).
\end{split}
\end{equation}

First, consider the concentration bound for $\|H_{jk}^{(1)}\|_\infty$. Denote $\widetilde{Y}_{ijk}(u,v) = Y_{ij}(u) Y_{ik}(v) - {\Sigma_{jk}(u,v)}/\{\sigma_j(u)^{1/2}\sigma_k(v)^{1/2}\}$ and let $d_{jk}((u,v),(u',v')) = d_j(u,u') + d_k(v,v').$ Applying Theorem~8.4 in \cite{Kosorok2008} under Conditions \ref{con_psi2_process} and \ref{con_finiteDistance}, we obtain that, there exists some constant $C_1>0$ such that $\big\|\sup_{u \in \cU}|Y_{1j}(u)|\big\|_{\psi_2} \leq C_1$
for all $j =1,\ldots,p.$ By the property of $\psi_1$-norm, we have that
\begin{equation*} 
\begin{split}
     &\left\|Y_{ij}(u)Y_{ik}(v)-Y_{ij}(u')Y_{ik}(v')\right\|_{\psi_1} \\ 
     &\leq \left\|Y_{ij}(u)\{Y_{ik}(v)-Y_{ik}(v')\}\right\|_{\psi_1} +\left\|\{Y_{ij}(u)-Y_{ij}(u')\}Y_{ik}(v')\right\|_{\psi_1}
      \\& \leq \left\|Y_{ij}(u)\right\|_{\psi_2}\left\|Y_{ik}(v)-Y_{ik}(v')\right\|_{\psi_2} + \left\|Y_{ik}(v')\right\|_{\psi_2}\left\|Y_{ij}(u)-Y_{ij}(u')\right\|_{\psi_2}
      \\& \lesssim  \{d_j(u,u') + d_k(v,v')\} = d_{jk}((u,v),(u',v')),
\end{split}
\end{equation*}
which implies that
\begin{equation} \label{eq_Y2distance}
     \left\|\widetilde{Y}_{ijk}(u,v)- \widetilde{Y}_{ijk}(u',v')\right\|_{\psi_1} \lesssim  d_{jk}((u,v),(u',v')).
\end{equation}
Note that
\begin{equation}\nonumber
  \bar{Z}_{jk}(u,v) = \frac{ \widehat \Sigma_{jk}(u,v) - \Sigma_{jk}(u,v)}{\sigma_j(u)^{1/2} \sigma_k(v)^{1/2}} = \frac{1}{n} \sum_{i=1}^n \left\{Y_{ij}(u) Y_{ik}(v) - \frac{\Sigma_{jk}(u,v)}{\sigma_j(u)^{1/2}\sigma_k(v)^{1/2}}\right\},
\end{equation}
and for a random variable $X$ and any integer $m\ge 1$, $\mathbb{E}\|X\|^m \le m! \|X\|_{\psi_1}^m.$ By Bernstein's inequality and Lemma 8.3 of \cite{Kosorok2008}, we have that for $u,v,u',v' \in \cU,$ 
\begin{equation} \nonumber
     \left\|{n^{1/2}}\Big\{\bar{Z}_{jk}(u,v)- \bar{Z}_{jk}(u',v')\Big\}\right\|_{\psi_1} \lesssim d_{jk}((u,v),(u',v')).
\end{equation}
For the semimetric $d_{jk},$ $D(\epsilon,d_{jk}) \le D(\epsilon/2,d_j) D(\epsilon/2,d_k) \lesssim \epsilon^{-2r}.$ Applying Theorem 8.4 in \cite{Kosorok2008} with Conditions \ref{con_psi2_process} and \ref{con_finiteDistance} again, we obtain that, there exists some constant $C_2>0$ such that 
\begin{equation*} \label{con_psi2_pointwise}
    \max_{1\le j,k \le p}\left\|\sup_{u,v \in \cU}|{n^{1/2}}\bar{Z}_{jk}(u,v)|\right\|_{\psi_1} \leq C_2.
\end{equation*}
This immediately implies that there exist some universal constant $C_3 >0$ such that for any $x>0,$
\begin{equation*}
P \left\{  \max_{j,k} \sup_{u,v \in \cU}\left|\frac{ \widehat \Sigma_{jk}(u,v) - \Sigma_{jk}(u,v)}{\sigma_j(u)^{1/2} \sigma_k(v)^{1/2}} \right| > x \right\} \lesssim p^2\exp\{ - C_3{n}^{1/2}x\}. 
\end{equation*}
As a result, for any $M >0$, there exists some constant $\tilde \rho_1 >0$ such that 
\begin{equation}
\label{ineq_sigma}
P \left\{  \max_{j,k} \sup_{u,v \in \cU}\left|\frac{ \widehat \Sigma_{jk}(u,v) - \Sigma_{jk}(u,v)}{\sigma_j(u)^{1/2} \sigma_k(v)^{1/2}} \right| > \tilde \rho_1 {\frac{\log p}{n^{1/2}}} \right\} \lesssim p^{-M}. 
\end{equation}
Observe that
\begin{equation*}
    \left|\frac{ \widehat \Sigma_{jk}(u,v)^2 - \Sigma_{jk}(u,v)^2}{\sigma_j(u) \sigma_k(v)} \right| \le \left|\frac{ \widehat \Sigma_{jk}(u,v) - \Sigma_{jk}(u,v)}{\sigma_j(u)^{1/2} \sigma_k(v)^{1/2}} \right|^2 + 2 \left|\frac{ \widehat \Sigma_{jk}(u,v) - \Sigma_{jk}(u,v)}{\sigma_j(u)^{1/2} \sigma_k(v)^{1/2}} \right|,
\end{equation*}
since $|\Sigma_{jk}(u,v)| \le \sigma_j(u)^{1/2} \sigma_k(v)^{1/2}.$ By the inequality~(\ref{ineq_sigma}), we have that 
\begin{equation}
\label{eq_I1}
P \left\{  \max_{j,k} \|H_{jk}^{(1)}\|_\infty > 2\tilde \rho_1 {\frac{\log p}{n^{1/2}}} + \tilde \rho_1^2 \frac{\log^2p}{n} \right\} \lesssim p^{-M}. 
\end{equation}

We next control the bound for $\|H_{jk}^{(2)}\|_\infty$ through the truncation technique.  Note that 
\begin{equation*}
    \frac{X_{ij}(u)^2X_{ik}(v)^2 - \widetilde{\Theta}_{jk}(u,v)}{\sigma_j(u) \sigma_k(v)} = Y_{ij}(u)^2 Y_{ik}(v)^2 - \frac{\widetilde{\Theta}_{jk}(u,v)}{\sigma_j(u) \sigma_k(v)}.
\end{equation*}
Define that $Y^*_{ij}(u) = Y_{ij}(u) I\left\{ \|Y_{ij}\|_\infty \leq C_4 \log^{1/2} (p\vee n)\right\}$ and 
$$
Z_{ijk}^*(u,v) = Y_{ij}^{*}(u)^2Y_{ik}^{*}(v)^2 - \mathbb{E} \{Y_{ij}^{*}(u)^2Y_{ik}^{*}(v)^2\}.
$$
By the property of $\psi_1$-norm and $|Y_{ij}^{*}(u)^2 - Y_{ij}^{*}(u')^2| \le 2 C_4 \log^{1/2} (p\vee n) |Y_{ij}^{*}(u) - Y_{ij}^{*}(u')|$, we have that
\begin{equation*} 
\begin{split}
     &\left\|Y_{ij}^{*}(u)^2Y_{ik}^{*}(v)^2-Y_{ij}^{*}(u')^2Y_{ik}^{*}(v')^2\right\|_{\psi_1} \\ 
     &\leq \left\|Y_{ij}^{*}(u)^2\{Y_{ik}^{*}(v)^2-Y_{ik}^{*}(v')^2\}\right\|_{\psi_1} +\left\|\{Y_{ij}^{*}(u)^2-Y_{ij}^{*}(u')^2\}Y_{ik}^{*}(v')^2\right\|_{\psi_1}
      \\& \lesssim  \log (p\vee n) \Big\{\|Y_{ij}^{*}(u)\|_{\psi_2}\left\|Y_{ik}^*(v)-Y_{ik}^*(v')\right\|_{\psi_2} + \|Y_{ik}^{*}(v')\|_{\psi_2}\left\|Y_{ij}^*(u)-Y_{ij}^*(u')\right\|_{\psi_2} \Big\}
      \\& \lesssim   \log (p\vee n) \{d_j(u,u') + d_k(v,v')\} \lesssim \log (p\vee n) d_{jk}((u,v),(u',v')),
\end{split}
\end{equation*}
which implies that, similar to (\ref{eq_Y2distance}),
\begin{equation} \nonumber
     \left\|Z^*_{ijk}(u,v)- Z^*_{ijk}(u',v')\right\|_{\psi_1} \lesssim \log (p\vee n) d_{jk}((u,v),(u',v')).
\end{equation}
Let $ \bar{Z}^*_{jk}(u,v) = n^{-1}\sum_{i=1}^n Z^*_{ijk}(u,v).$
We apply the similar technique of $\bar{Z}_{jk}$ above to  the term $\bar{Z}^*_{jk}$ and obtain that there exists some universal constant $C_5 >0$ such that for any $x>0,$
\begin{equation*}
P \left\{  \max_{j,k} \sup_{u,v \in \cU}\left|\frac{ \bar{Z}^*_{jk}(u,v)}{\log (p\vee n)} \right| > x \right\} \lesssim p^2\exp( - C_5{n^{1/2}}x). 
\end{equation*}
As a result, for any $M >0$, there exists some constant $\tilde \rho_2 >0$ such that 
\begin{equation}
\nonumber
P \left\{  \max_{j,k} \sup_{u,v \in \cU}\left|\bar{Z}^*_{jk}(u,v) \right| > \tilde \rho_2 {\frac{\log^2 (p\vee n)}{n^{1/2}}} \right\} \lesssim p^{-M}. 
\end{equation}
Now we consider the bound of the term $\|Y_{ij}\|_{\infty}$. By Conditions \ref{con_psi2_process}-\ref{con_finiteDistance} and Theorem 8.4 of \cite{Kosorok2008}, we immediately have that there exists some constant $C_6>0$
\begin{equation*}
    \max_{1 \le i\le n, 1\le j\le p}\Big\|\sup_{u \in \cU}|Y_{ij}(u)|\Big\|_{\psi_2} \le C_6,
\end{equation*}
which also implies that there exists some constant $C_7 >0$ such that for any  $x >0,$ 
\begin{equation*}
    P \left \{ \max_{1 \le i\le n, 1\le j\le p}\|Y_{ij}(u)\|_{\infty} > x\right\} \lesssim n p \exp(- C_7 x^2).
\end{equation*}
Hence we obtain that for any $M >0$, there exists some constant $C_4 >0$ such that 
\begin{equation}
\label{maximal_Y}
    P \left \{ \max_{1 \le i\le n, 1 \le  j\le p}\|Y_{ij}\|_{\infty} > C_4 \log^{1/2}(p \vee n)\right\} \lesssim (p \vee n)^{-M}.
\end{equation}
On the event 
$$
\Omega_{n0} = \Big\{\max_{1 \le i\le n, 1 \le j\le p}\|Y_{ij}\|_{\infty} \le C_4 \log^{1/2}(p \vee n)\Big\},
$$ we find that 
\begin{equation*}
\begin{split}
      Y_{ij}(u)^{2} Y_{ik}(v)^{2} - \frac{\widetilde{\Theta}_{jk}(u,v)}{\sigma_j(u)\sigma_k(v)}
     =~ & Y_{ij}^{*}(u)^2 Y_{ik}^{*}(v)^2 - \mathbb{E} \Big\{Y_{ij}^{*}(u)^2 Y_{ik}^{*}(v)^2\Big\} \\
     & + \mathbb{E} \Big\{Y_{ij}^{*}(u)^2 Y_{ik}^{*}(v)^2 - Y_{ij}(u)^{2} Y_{ik}(v)^{2}\Big\}.
     \end{split}
\end{equation*}
Note that $Y_{ij}^{*}(u)^2 - Y_{ij}(u)^2 = Y_{ij}(u)^2I\{\|Y_{ij}\|_{\infty} > C_4 \log^{1/2}(p \vee n)\}.$ 
By the inequality (\ref{maximal_Y}), we can obtain that 
$$
\left |\mathbb{E} \Big\{Y_{ij}^{*}(u)^2 Y_{ik}^{*}(v)^2 - Y_{ij}(u)^{2} Y_{ik}(v)^{2}\Big\}\right| \lesssim (p \vee n)^{-M}.
$$
Therefore, for any $M>0,$ there exist some constant $\tilde \rho_3>0$ such that 
\begin{equation}
\label{eq_I2}
\begin{split}
 P \left \{ \max_{1 \le j\le p}\|H_{jk}^{(2)}\|_{\infty} > \tilde \rho_3 {\frac{\log^2 (p\vee n)}{n^{1/2}}}\right\} \lesssim p^{-M}.
     \end{split}
\end{equation}
Combining (\ref{eq_I1}) and (\ref{eq_I2}), we obtain that for any $M>0$, there exists some constant $\rho_1 >0$ such that 
\begin{equation*}
    P\left \{ \max_{j,k} \left\| \frac{\widehat \Theta_{jk}  - \Theta_{jk}}{{\Theta_{jk}}}\right\|_{\infty} \geq \rho_1 {\frac{\log^2 (p \vee n)}{n^{1/2}}} \right\} \lesssim p^{-M}.
\end{equation*}
The proof is complete. $\square$

\begin{lemma}
\label{lm_thetasqrt}
Suppose that Conditions~\ref{con_psi2_process}--\ref{con_rate} hold. Then for any  $M >0$, there exist some  constant $\rho_2 >0$ such that
\begin{equation} \label{eq_thetasqrt}
    \max_{j,k}\left\| \frac{{\Theta_{jk}^{1/2}}-{ \widehat \Theta_{jk}^{1/2}}}{{ \widehat \Theta_{jk}^{1/2}}} \right\|_\infty \leq  \rho_2 \frac{\log^2 p}{n^{1/2}}
\end{equation}
with probability greater than $1- O(p^{-M}).$
\end{lemma}

\textbf{Proof}.  Let the event $\Omega_{n}(s) = \{\|(\widehat \Theta_{jk} - \Theta_{jk})/\Theta_{jk}\|_\infty\leq s {\log^2 p/n^{1/2}} \le 1/2\}.$ For any  $M>0,$ it follows from Lemma \ref{lm_theta} that there exists some constant $\rho_1 >0$ such that
 $    P\{\Omega_n(\rho_1)\} \ge  1 - O(p^{-M}).$ 
Since
\begin{equation} \nonumber
    \begin{split}
        \left\|\frac{\Theta_{jk} }{\widehat\Theta_{jk}}\right\|_\infty = \left\|\frac{\Theta_{jk} - \widehat \Theta_{jk}}{\widehat \Theta_{jk}} +1 \right\|_\infty \leq \left\|\frac{\Theta_{jk} - \widehat \Theta_{jk}}{\Theta_{jk}}\right\|_\infty \left\|\frac{\Theta_{jk} }{\widehat\Theta_{jk}}\right\|_\infty +1,
    \end{split}
\end{equation} 
hence, on the event $\Omega_n(\rho_1),$  
we have that $\|\Theta_{jk}/\widehat \Theta_{jk}\|_\infty \leq 2.$ As a  result, on the event $\Omega_n(\rho_1),$ it  follows that 
\begin{equation*} \nonumber
    \begin{split}
        \left\| \frac{{\Theta_{jk}^{1/2}}-{ \widehat \Theta_{jk}^{1/2}}}{{ \widehat \Theta_{jk}^{1/2}}} \right\|_\infty = \left\| \frac{{\Theta_{jk}}-{ \widehat \Theta_{jk}}}{{ \widehat \Theta_{jk}} + {\widehat \Theta_{jk}^{1/2} \Theta_{jk}^{1/2}}} \right\|_\infty 
        \leq \left\|\frac{\Theta_{jk} - \widehat \Theta_{jk}}{\Theta_{jk}}\right\|_\infty \left\|\frac{\Theta_{jk} }{\widehat\Theta_{jk}}\right\|_\infty \leq 2 \rho_1 {\frac{\log^2 p}{n^{1/2}}}.
    \end{split}
\end{equation*}
Take $\rho_2 =  2  \rho_1$ and  the proof is complete.
$\square$


\begin{lemma} \label{lm_cb_ratio}
Suppose that Conditions~\ref{con_psi2_process}--\ref{con_rate} holds. Then for any $M >0$,  there exist some  positive constant  $\rho_3>0$  such that
\begin{equation}\nonumber
\max_{j,k} \left\| \frac{\widehat \Sigma_{jk}  - \Sigma_{jk}}{{\widehat \Theta_{jk}^{1/2}}}\right\|_\cS \le  \rho_3  \left(\frac{\log p}{n}\right)^{1/2}
\end{equation}
with probability greater than $1-O(p^{-M}).$
\end{lemma}


\textbf{Proof}.  Let  $\widetilde{Y}_{ijk}(u,v) = Y_{ij}(u) Y_{ik}(v) - {\Sigma_{jk}(u,v)}/\{\sigma_j(u)^{1/2}\sigma_k(v)^{1/2}\}$ and
\begin{equation}\nonumber
  \bar{Z}_{jk}(u,v) = \frac{ \widehat \Sigma_{jk}(u,v) - \Sigma_{jk}(u,v)}{\sigma_j(u)^{1/2} \sigma_k(v)^{1/2}} = \frac{1}{n} \sum_{i=1}^n \widetilde{Y}_{ijk}(u,v). 
\end{equation}
We first derive the concentration bound of $\|\bar{Z}_{jk}\|_{\cS}.$  It  follows from the proof  of Lemma  \ref{lm_theta} that there exists some constant $C_8 >0$ such that
$$
\max_{j,k} \left\| \sup_{u,v \in \cU}\Big| \widetilde{Y}_{1jk}(u,v)\Big|\right\|_{\psi_1} \le C_8.
$$
which further implies that $\max_{j,k} \left\| \big\| \widetilde{Y}_{1jk}\big\|_{\cS}\right\|_{\psi_1} \le C_8.$
As a result, it follows from Theorem~2.5 of \textcolor{blue}{Bosq (2000)}
that there exists some universal constant $C_9>0$ such that for any $x>0$
\begin{equation} \nonumber
    P\left( \left \| \bar{Z}_{jk} \right\|_{\cS} \geq x \right) \leq 2 \exp\{-C_9n \min(x^2, x)\}.
\end{equation}
For any $M >0,$ there exists some constant $\tilde \rho >0$  that
\begin{equation} \label{eq_cb_sigma}
\left\|\bar{Z}_{jk} \right\|_{\cS} \le \tilde \rho  \left(\frac{\log p}{n}\right)^{1/2}
\end{equation}
with probability greater than $1 - O(p^{-M}).$  

Now we derive the  bound of 
$\left\| \big(\widehat \Sigma_{jk}  - \Sigma_{jk}\big)/\widehat \Theta_{jk}^{1/2}\right\|_\cS$.  Note that Condition \ref{con_variance_function} implies that $\Theta_{jk}(u,v) \geq \tau \sigma_{j}(u) \sigma_{k}(v).$  
We obtain  that
\begin{equation}\nonumber
\begin{split}
     \left\| \frac{\widehat \Sigma_{jk}  - \Sigma_{jk}}{{\widehat \Theta_{jk}^{1/2}}}\right\|_\cS  
     \leq 
     \left\| \frac{\widehat \Sigma_{jk}  - \Sigma_{jk}}{{ \Theta_{jk}^{1/2}}}\right\|_\cS \left\| \frac{{\Theta_{jk}^{1/2}}}{{ \widehat \Theta_{jk}^{1/2}}}\right\|_\infty
      \leq \left\| \tau^{-1/2} \bar{Z}_{jk}\right\|_\cS
     \left(\left\| \frac{{\Theta_{jk}^{1/2}}-{ \widehat \Theta_{jk}^{1/2}}}{{ \widehat \Theta_{jk}^{1/2}}}\right\|_\infty +1\right).
\end{split}
\end{equation}
Hence, together with (\ref{eq_cb_sigma}) and Lemma \ref{lm_thetasqrt}, the lemma follows. The proof is complete.
$\square$


\textbf{Proof of Theorem~\ref{thm_rate}}.
For easy representation, define 
\begin{equation*}
\widehat\Phi_{jk}(u,v)= \frac{\widehat \Sigma_{jk}(u,v)}{\widehat \Theta_{jk}(u,v)^{1/2}},\quad \widetilde\Phi_{jk}(u,v)= \frac{\Sigma_{jk}(u,v)}{\widehat \Theta_{jk}(u,v)^{1/2}} \quad \mbox{and} \quad \Phi_{jk}(u,v) = \frac{\Sigma_{jk}(u,v)}{\Theta_{jk}(u,v)^{1/2}}. 
\end{equation*}
Let 
$$
\Omega_{n1} = \Big\{\max_{j,k}\|\widehat\Phi_{jk}- \widetilde\Phi_{jk}\|_\cS\leq \lambda \Big\}, \Omega_{n2} = \left\{\max_{j,k}\left\|\frac{\widehat \Theta_{jk} - \Theta_{jk}}{\Theta_{jk}}\right\|_{\infty} \le  \frac{1}{2}\right\}.
$$
It is immediate to see that under the  event $\Omega_{n2},$ $2^{-1} \|\Theta_{jk}\|_{\infty} \le \|\widehat \Theta_{jk}\|_{\infty} \le 2 \|\Theta_{jk}\|_{\infty}$ for all $j$ and $k.$ 
By Conditions~\ref{con_psi2_process}--\ref{con_variance_function}, we have $\Theta_{jk}(u,v) \leq C'  \sigma_{j}(u) \sigma_{k}(v)$ and $
\Theta_{jk}(u,v) \ge \tau  \sigma_{j}(u) \sigma_{k}(v) $
Then under the event $\Omega_{n1} \cap \Omega_{n2}$ and Conditions (i)-(iii) on $S_{\lambda}(Z)$, we obtain that
\begin{equation*}
    \begin{split}
        &\sum_{k = 1}^p\|\widehat \Sigma_{jk}^\A - \Sigma_{jk}\|_\cS
        \\=& \sum_{k = 1}^p\|\widehat \Sigma_{jk}^\A - \Sigma_{jk}\|_\cS I\{\|\widehat\Phi_{jk}\|_\cS\geq \lambda\} + \sum_{k = 1}^p\| \Sigma_{jk}\|_\cS I\{\|\widehat\Phi_{jk}\|_\cS < \lambda\}
        \\ 
        \leq &\sum_{k = 1}^p \left\{ \|s_{\lambda}( \widehat\Phi_{jk}) - \widehat\Phi_{jk} \|_\cS + \|\widehat\Phi_{jk}-\widetilde\Phi_{jk}\|_\cS \right\}\big\|\widehat \Theta_{jk}^{1/2} \big\|_\infty I\{\|\widehat\Phi_{jk}\|_\cS \geq  \lambda, \|\widetilde \Phi_{jk}\|_\cS \geq  \lambda\}
        \\ 
        & +\sum_{k = 1}^p  \Big\|\big[s_{\lambda}( \widehat\Phi_{jk}) - \widetilde \Phi_{jk} \big]\widehat \Theta_{jk}^{1/2} \Big\|_\cS I\{\|\widehat\Phi_{jk}\|_\cS \geq  \lambda, \|\widetilde \Phi_{jk}\|_\cS <  \lambda\} + \sum_{k = 1}^p\| \Sigma_{jk}\|_\cS I\{\|\widetilde \Phi_{jk}\|_\cS < 2\lambda\}
        \\ 
        \leq &\sum_{k=1}^p 2\lambda \big\|\widehat \Theta_{jk}^{1/2} \big\|_\infty I\{  \| \widetilde \Phi_{jk}\|_\cS \geq  \lambda\} + \sum_{k=1}^p (1+c) {\|\widetilde \Phi_{jk}\|_\cS} \|\widehat \Theta_{jk}^{1/2}\|_{\infty}I\{  \|\widetilde \Phi_{jk}\|_\cS <  \lambda\}
        \\ 
        &+ \sum_{k=1}^p \|\widetilde \Phi_{jk}\|_\cS \big\|\widehat \Theta_{jk}^{1/2} \big\|_\infty I\{  \|\widetilde \Phi_{jk}\|_\cS <2  \lambda\}
        \\
        \lesssim &~ \lambda^{1-q} \sum_{k=1}^p \big\|\widehat \Theta_{jk}\big\|_\infty^{1/2} \|\widetilde \Phi_{jk}\|_\cS^q
        \lesssim   \lambda^{1-q}  \sum_{k=1}^p \big\|\sigma_{j} \big\|_\infty^{(1-q)/2} \big\|\sigma_{k}\big\|_\infty^{(1-q)/2}\|\Sigma_{jk}\|_\cS^q
      \lesssim  s_0(p)\left(\frac{\log p}{n}\right)^{\frac{1-q}{2}}.
    \end{split}
 \end{equation*}
 Since there exists some constant $\delta >0$ such that $P\{\Omega_{n1}^{C}\} + P\{\Omega_{n2}^C\}  \lesssim p^{-M},$ the theorem follows. 
$\square$

\textbf{Proof of Theorem~\ref{thm_supp}}. 
We consider two sets: $S_{n1} = \{(j,k): \|\widehat \Sigma_{jk}^\A\|_\cS \neq 0 ~\text{and}~\| \Sigma_{jk}\|_\cS= 0 \}$ and $S_{n2} = \{(j,k): \|\widehat \Sigma_{jk}^\A\|_\cS = 0 ~\text{and}~\| \Sigma_{jk}\|_\cS \neq 0 \}.$ It suffices to prove that 
$$
P\big(|S_{n1}| >0 \big) + P\big(|S_{n2}| >0 \big) \to 0,
$$
as $n,p \to \infty.$
By Conditions (i)-(iii) on $S_{\lambda}(Z),$
$$
S_{n1} =  \left\{(j,k): \left\|\frac{\widehat \Sigma_{jk}}{\widehat \Theta_{jk}^{1/2}}\right\|_\cS > \lambda ~\text{and}~\| \Sigma_{jk}\|_\cS= 0 \right\} \subset \left\{(j,k): \left\|\frac{\widehat \Sigma_{jk} - \Sigma_{jk}}{\widehat \Theta_{jk}^{1/2}}\right\|_\cS > \lambda \right\}
$$
Therefore, with the choice $\lambda= \delta ({\log p}/{n})^{1/2},$ we obtain
\begin{equation} \label{eq_S1}
    \begin{split}
        P(|S_{n1}|>0) \leq P\left\{\max_{j,k} \left\|\frac{\widehat \Sigma_{jk} - \Sigma_{jk}}{\widehat \Theta_{jk}^{1/2}}\right\|_\cS > \lambda\right\} \lesssim p^{-M}.
    \end{split}
\end{equation}
for some prespecified $M >0.$
Similarly, we have 
\begin{equation} \nonumber
    \begin{split}
        S_{n2} &= \left\{(j,k): \left\|\frac{\widehat \Sigma_{jk}}{\widehat \Theta_{jk}^{1/2}}\right\|_\cS \leq \lambda ~\text{and}~\| \Sigma_{jk}\|_\cS \neq 0 \right\}. 
    \end{split}
\end{equation}
Note that $\|\Sigma_{jk}\|_\cS \neq 0$ implies that 
\begin{equation} \label{eq_S2_con}
    \begin{split}
        (2\delta +\gamma)\left(\frac{\log p}{n}\right)^{1/2}< \left\|\frac{  \Sigma_{jk}}{\Theta_{jk}^{1/2}}\right\|_\cS \leq \left[ \left\|\frac{\Sigma_{jk} - \widehat \Sigma_{jk}}{\widehat \Theta_{jk}^{1/2}}\right\|_\cS + \left\|\frac{\widehat \Sigma_{jk}}{\widehat \Theta_{jk}^{1/2}}\right\|_\cS\right]\left\|\frac{  \widehat \Theta_{jk}^{1/2}}{\Theta_{jk}^{1/2}}\right\|_\infty.
    \end{split}
\end{equation}
Let $\Omega_{n3} = \Big\{\|(\widehat \Theta_{jk}^{1/2} - \Theta_{jk}^{1/2})/\widehat \Theta_{jk}^{1/2}\|_\infty\leq \epsilon \Big\}$ for  some small constant $0 < \epsilon < \gamma/(4\delta + 2\gamma).$
Conditioned on the event of $\Omega_{n3}$, the  inequality
\begin{equation} \nonumber
    \begin{split}
        \left\|\frac{\widehat \Theta_{jk}^{1/2} }{\Theta_{jk}^{1/2}}\right\|_\infty  \leq \left\|\frac{\widehat\Theta_{jk}^{1/2} - \Theta_{jk}^{1/2}}{\widehat\Theta_{jk}^{1/2}}\right\|_\infty \left\|\frac{\widehat\Theta_{jk}^{1/2} }{\Theta_{jk}^{1/2}}\right\|_\infty +1
    \end{split}
\end{equation} 
implies that  $\|\widehat\Theta_{jk}^{1/2}/\Theta_{jk}^{1/2}\|_\infty \leq 1/(1-\epsilon).$ This together with (\ref{eq_S2_con}) shows that 
\begin{equation} \nonumber
    \begin{split}
        S_{n2} \cap \Omega_{n3} 
      & \subset \left\{(j,k): \left\|\frac{\widehat \Sigma_{jk} - \Sigma_{jk}}{\widehat \Theta_{jk}^{1/2}}\right\|_\cS >  \delta \left(\frac{\log p}{n}\right)^{1/2} \right\}.
    \end{split}
\end{equation}
As a result, 
\begin{equation} \label{eq_S2}
    \begin{split}
        P(|S_{n2}|>0)  \leq P(\Omega_{n3}^{C}) + P\left\{\max_{j,k} \left\|\frac{\widehat \Sigma_{jk} - \Sigma_{jk}}{\widehat \Theta_{jk}^{1/2}}\right\|_\cS > \delta \left(\frac{\log p}{n}\right)^{1/2}\right\} 
        \lesssim p^{-M}.
    \end{split}
\end{equation}
Combining (\ref{eq_S1}) and (\ref{eq_S2}), we complete our proof.
$\square$

\section{Examples of functional thresholding operators}
\label{supp.sec_con}
In Section~\ref{supp_fr_con}, we verify that our proposed soft, SCAD and adaptive lasso functional thresholding rules satisfy conditions~(i)--(iii) in Section~\ref{sec.method}. We then present the derivations of these three functional thresholding rules in Section~\ref{supp_fr_penalty}.

\subsection{Verification of conditions~(i)--(iii)} 
\label{supp_fr_con}
It is directly implied from the thresholding rules that the soft, SCAD and adaptive lasso functional methods satisfy condition~(ii). Since the soft functional thresholding has the largest amount of functional shrinkage in the Hilbert--Schmidt norm compared with SCAD and adaptive lasso methods, it suffices to show that the soft functional thresholding satisfies condition~(iii). For $\|Z\|_\cS \leq \lambda,$  the thresholding effect leads to $\|0-Z\|_\cS \leq \lambda.$ When $\|Z\|_\cS > \lambda,$ we obtain that $\|Z \lambda / \|Z\|_\cS\|_\cS =  \lambda.$ 

We next show that the above three thresholding methods satisfy condition (i). By the triangle inequality, $\|Z-Y\|_\cS \leq \lambda$ in condition (i) implies that $\big|\| Z\|_\cS -\lambda\big| \leq \|Y\|_\cS.$
\begin{itemize}
    \item Soft functional thresholding: If $\|Z\|_\cS \leq \lambda,$  $0 \leq c\|Y\|_\cS$ directly holds for all $Y \in \mathbb{S}$ and $c>0.$ When $\|Z\|_\cS > \lambda,$ we have $\|s^\tS_\lambda(Z)\|_{\cS}=  \| Z\|_\cS -\lambda \leq \|Y\|_\cS$ with the choice of  $c = 1.$
    \item SCAD functional thresholding: When $\|Z\|_\cS \leq 2\lambda,$  $s^\SC_{\lambda}(Z)$ is the same as the soft functional thresholding rule. 
    For $\|Z\|_\cS >2\lambda,$ we have $\|s^\SC_\lambda(Z)\|_{\cS} \leq \|Z\|_\cS\leq \|Y\|_\cS + \lambda \leq \|Y\|_\cS +\|Z\|_\cS/2$ and hence $\|s^\SC_\lambda(Z)\|_{\cS} \leq \|Z\|_\cS \leq 2 \|Y\|_\cS.$ Combining the above results, we take $c = 2.$
    \item Adaptive lasso functional thresholding:
    Let $\ceil{\eta} $ denote  the smallest integer greater than or equal to $\eta.$ For $\|Z\|_\cS \leq \lambda,$ this condition holds for all $Y \in \mathbb{S}$ and $c>0.$ For $\|Z\|_\cS >\lambda,$ we have that 
$\|s^\AL_{\lambda}(Z)\|_\cS = \|Z(1-\lambda^{\eta+1}/\|Z\|_{\cS}^{\eta+1})\|_\cS = (\|Z\|_{\cS}^{\eta+1} -\lambda^{\eta+1})/\|Z\|_{\cS}^{\eta} \leq (\|Z\|_{\cS}^{\ceil{\eta}+1} -\lambda^{\ceil{\eta}+1})/\|Z\|_{\cS}^{\ceil{\eta}} = (\|Z\|_\cS - \lambda)(\|Z\|_{\cS}^{\ceil{\eta}} + \|Z\|_{\cS}^{\ceil{\eta}-1}\lambda+ \cdots + \lambda^{\ceil{\eta}})/\|Z\|_{\cS}^{\ceil{\eta}} \leq (\ceil{\eta}+1) \|Y\|_\cS.$ Hence, for any $\eta\geq 0,$ we can find $c =\ceil{\eta}+1.$ In the special case of $\eta=0,$ $s^\AL_{\lambda}(Z)$ degenerates to the soft functional thresholding rule with $c = 1,$ which is consistent with our finding for the soft functional thresholding.
\end{itemize}

\subsection{Derivations of the functional thresholding rules from various penalty functions} \label{supp_fr_penalty}
Soft functional thresholding can be obtained via
\begin{equation}
\label{func_pls_soft}
    s^\tS_{\lambda}(Z) =\underset{\theta \in {\mathbb S}}{\arg\min} \left\{\frac{1}{2} \|\theta - Z\|_{\cS}^2 + \lambda \|\theta\|_{\cS}\right\}.
\end{equation}
First, we show that if $\|Z\|_\cS \le \lambda,$ then $\|s^\tS_{\lambda}(Z)\|_\cS = 0$ and hence $s^\tS_{\lambda}(Z) = 0.$ This results from the fact that, for any $\theta,$
\begin{equation*}
\begin{split}
    \frac{1}{2} \|\theta - Z\|_{\cS}^2 + \lambda \|\theta\|_{\cS} &\ge \frac{1}{2} \big(\|\theta\|_\cS - \|Z\|_\cS\big)^2 + \lambda \|\theta\|_{\cS}\\
    & = \frac{1}{2} \|\theta\|_\cS^2 +(\lambda - \|Z\|_\cS) \|\theta\|_\cS + \frac{1}{2}\|Z\|_{\cS}^2 \ge \frac{1}{2}\|Z\|_{\cS}^2.
\end{split}
\end{equation*}
Second, we show that if $\|Z\|_\cS > \lambda,$ then $\|s^\tS_{\lambda}(Z)\|_\cS \neq 0.$ In fact, we can find $\theta_c = c Z$ with $ c= 1 - \lambda/\|Z\|_\cS > 0$ such that 
\begin{equation*}
    \frac{1}{2} \|\theta_c - Z\|_{\cS}^2 + \lambda \|\theta_c\|_{\cS}  = \frac{1}{2} (1 - c)^2\|Z\|_\cS^2 +  \lambda c\|Z\|_{\cS} < \frac{1}{2}\|Z\|_\cS^2.
\end{equation*}
As a result, we are able to take the first derivative of (\ref{func_pls_soft}) with respect to $\theta$ and set $p'_\lambda(\theta) =\theta - Z +\lambda \theta/\|\theta\|_\cS = 0.$ Thus, 
$\widehat \theta = Z {\|\widehat \theta\|_\cS}/\big(\|\widehat \theta\|_\cS + \lambda\big),$
which implies that $ \|\widehat \theta\|_\cS  = \|Z\|_\cS - \lambda $. Combining the above results, we have that 
$\widehat \theta = Z (1-\lambda/\|Z\|_\cS)_+.$
 
The SCAD and adaptive lasso functional thresholding rules can be derived in a similar fashion. Hence, we only present their penalty functions here. The functional version of SCAD penalty takes the form of 
$$p_\lambda(\theta) = \lambda \|\theta\|_\cS I(\|\theta\|_\cS \leq \lambda) + \frac{2a\lambda \|\theta\|_\cS - \|\theta\|_\cS^2 -\lambda^2}{2(a-1)} I(\lambda < \|\theta\|_\cS \leq a\lambda) + \frac{\lambda^2(a+1)}{2} I (\|\theta\|_\cS> a\lambda), $$
for $a>2.$ For the functional version of adaptive lasso penalty, we use 
$p_\lambda(\theta) = \lambda^{\eta+1} \|Z\|_\cS^{-\eta} \|\theta\|_\cS,$ for $\eta \ge 0.$ A similar adaptive lasso penalty function operating on $|\cdot|$ for the univariate scalar case can be found in \cite{rothman2009}.

\section{Partially observed functional data} 
\label{supp.sec_partial}
Section~\ref{supp.post} gives the expression of the local linear surface smoother for the cross-covariance estimation. Section~\ref{supp.pre} presents the details of pre-smoothing for densely sampled functional data. Section~\ref{supp.tp.partial} provides all technical proofs for the partially observed functional scenario. Section~\ref{supp.verify.In} presents the heuristic verification of $I_{jk}$ in (\ref{In}) and Condition~\ref{cond_sparse_rate2}.

\subsection{Local linear surface smoother}
\label{supp.post}
We use (\ref{eq.kernel.cross.cov}) to derive the expression of its minimizer. Recall $T_{ab,ijk}$ and $S_{ab, jk}$ in (\ref{T.est}) and (\ref{S.est}), respectively, for $a,b=0,1,2,$ $i=1, \dots, n$ and $j,k=1, \dots, p.$
To minimize the objective in (\ref{eq.kernel.cross.cov}), some calculations lead to the resulting estimator
\begin{equation} \label{eq.kernel.w}
\begin{split}
     \widehat \Sigma_{jk} =& \sum_{i=1}^n\frac{(S_{20}S_{02}-S_{11}^2)T_{00,ijk} - (S_{10}S_{02}-S_{01}S_{11})T_{10,ijk} +(S_{10}S_{11}-S_{01}S_{20})T_{01,ijk}}{(S_{20}S_{02}-S_{11}^2)S_{00} - (S_{10}S_{02}-S_{01}S_{11})S_{10} +(S_{10}S_{11}-S_{01}S_{20})S_{01}}
     \\ := & \sum_{i=1}^n \left( W_{1,jk} T_{00,ijk}+ W_{2,jk} T_{10,ijk}+W_{3,jk}T_{01,ijk} \right),
\end{split}
\end{equation}
where we drop subscripts $j,k$ in $S_{ab,jk}$'s to simplify the notation. Note that, under Model~(\ref{eq.kernel.simple}), $S_{ab,jk}$'s no longer depend on $j,k,$ and hence subscripts $j,k$ in $S_{ab,jk}$'s can be dropped.

\subsection{Pre-smoothing} 
\label{supp.pre}
When each random function $X_{ij}(\cdot)$ is densely observed with errors satisfying Model~(\ref{model.partial}), the commonly adopted pre-smoothing approach applies local linear smoother to estimate each $X_{ij}(\cdot)$ before subsequent analysis. The reconstructed individual function is obtained by $\widehat X_{ij}(u)=\hat a_0,$ where
\begin{equation} \nonumber
     (\hat a_0, \hat a_1) =  \argmin_{a_0, a_1}  \sum_{l = 1}^{L_{ij}}  \big\{Z_{ijl} - a_0-a_1(U_{ijl}-u)\big\}^2 K_{h_X}({U_{il}-u}).
\end{equation}
Let $T_{a,ij}(u) = \sum_{l = 1}^{L_{ij}} K_{h_X}({U_{ijl}-u}) (U_{ijl}-u)^aZ_{ijl}$ and $S_{a,ij}(u)=\sum_{l = 1}^{L_{ij}} K_{h_X}({U_{ijl}-u}) (U_{ijl}-u)^b$ for $a=0,1,2.$ Solving the minimization problem above yields that
$$
\widehat X_{ij}(u) = \frac{S_{2,ij}(u)T_{0,ij}(u) - S_{1,ij}(u)T_{1,ij}(u)}{S_{2,ij}(u)S_{0,ij}(u)-\{S_{1,ij}(u)\}^2}.
$$
Under the simplified model in (\ref{eq.kernel.simple}), we drop the subscript $j$ in $L_{ij}$ and $S_{a,ij}$ in the expression of $\widehat X_{ij}(u)$ above.
For an equally-spaced grid of $R$ points $u_1< \cdots < u_R \in \cU,$ the binned approximation of $\widehat X_{ij}(u)$ is
$$
\widecheck X_{ij}(u) = \frac{\widecheck S_{2,i}(u)\widecheck T_{0,ij}(u) - \widecheck S_{1,i}(u)\widecheck T_{1,ij}(u)}{\widecheck S_{2,i}(u)\widecheck S_{0,i}(u)-\{\widecheck S_{1,i}(u)\}^2},
$$
where $ \widecheck T_{a,ij}(u) = \sum_{r = 1}^R K_{h_X}({u_{r}-u}) (u_{r}-u)^a \cD_{r,ij}$ and $ \widecheck S_{a,i}(u) = \sum_{r = 1}^R K_{h_X}({u_{r}-u}) (u_{r}-u)^a \varpi_{r,i}.$ See also Table~\ref{comp.table.pre} for the computational complexity analysis of the pre-smoothing based on local linear smoother and its binned implementation, denoted as LLS-P and BinLLS-P respectively, under Models~(\ref{model.partial}) and (\ref{eq.kernel.simple}).  

\begin{table}[!htbp]
\caption{\label{comp.table.pre} The computational complexity analysis of LLS- and BinLLS-based pre-smoothings under Models~(\ref{model.partial}) and (\ref{eq.kernel.simple}) when evaluating the reconstructed functions at a grid of $R$ points.}
	\begin{center}
		\resizebox{5.1in}{!}{
			\begin{tabular}{cccc}
    \hline
        Method & Model & \begin{tabular}[c]{@{}c@{}}Number of  \\ kernel evaluations\end{tabular}  & \begin{tabular}[c]{@{}c@{}}Number of operations  \\(additions and multiplications)\end{tabular}  \\ \hline
        LLS-P & (\ref{model.partial}) & $O( R\sum_{i= 1}^n \sum_{j = 1}^p L_{ij})$& $O(R\sum_{i = 1}^n\sum_{j=1}^pL_{ij})$  \\ 
        LLS-P  & (\ref{eq.kernel.simple}) & $O(R \sum_{i= 1}^n  L_{i})$& $O(pR\sum_{i=1}^nL_{i})$ \\ 
        BinLLS-P & (\ref{eq.kernel.simple})& $O(R)$ & $O(npR^2+p\sum_{i = 1}^n L_{i})$ \\ \hline
    \end{tabular}
		}	
	\end{center}
\end{table}

\subsection{Technical proofs}
\label{supp.tp.partial}
\textbf{Proof of Theorem~\ref{thm_sparse_2}}.
Define
\begin{equation*}
\widetilde\Lambda_{jk}(u,v)= \frac{\widetilde \Sigma_{jk}(u,v)}{\widetilde \Psi_{jk}(u,v)^{1/2}},\quad \widecheck\Lambda_{jk}(u,v)= \frac{\Sigma_{jk}(u,v)}{\widetilde \Psi_{jk}(u,v)^{1/2}} \quad \mbox{and} \quad \Lambda_{jk}(u,v) = \frac{\Sigma_{jk}(u,v)}{\Psi_{jk}(u,v)^{1/2}}.
\end{equation*}
Let
$$
\widetilde \Omega_{n1} = \Big\{\max_{j,k}\|\widetilde\Lambda_{jk}- \widecheck\Lambda_{jk}\|_\cS\leq \lambda \Big\}, ~~\widetilde  \Omega_{n2} = \left\{\max_{j,k}\left\|\frac{\widetilde \Psi_{jk} - \Psi_{jk}}{\Psi_{jk}}\right\|_{\infty} \le  \frac{1}{2}\right\}.
$$

First, we can obtain from Condition~\ref{cond_sparse_rate2} that $P(\widetilde \Omega_{n2}^C) = o(1).$ Note that 
\begin{equation}\nonumber
\begin{split}
     \left\| \frac{\widetilde \Sigma_{jk}  - \Sigma_{jk}}{{\widetilde \Psi_{jk}^{1/2}}}\right\|_\cS
     \leq
     \left\| \frac{\widetilde \Sigma_{jk}  - \Sigma_{jk}}{{ \Psi_{jk}^{1/2}}}\right\|_\cS \left\| \frac{{\Psi_{jk}^{1/2}}}{{ \widetilde \Psi_{jk}^{1/2}}}\right\|_\infty
      \lesssim \left\| \widetilde \Sigma_{jk}  - \Sigma_{jk}\right\|_{\cS}
     \left(\left\| \frac{{\Psi_{jk}^{1/2}}-{ \widetilde \Psi_{jk}^{1/2}}}{{ \widetilde \Psi_{jk}^{1/2}}}\right\|_\infty +1\right).
\end{split}
\end{equation}
It follows from Condition~\ref{cond_sparse_rate1} that there exists some constant $\tilde \delta >0$ such that $ P\{(\widetilde\Omega_{n1})^{C}\} = o(1).$
We also can see that under the event $\widetilde\Omega_{n2},$ $2^{-1} \|\Psi_{jk}\|_{\infty} \le \|\widetilde \Psi_{jk}\|_{\infty} \le 2 \|\Psi_{jk}\|_{\infty}$ for all $j$ and $k.$
Then on the event $\widetilde\Omega_{n1} \cap \widetilde\Omega_{n2}$ and Conditions (i)-(iii) on $S_{\lambda}(Z)$, we obtain that
\begin{equation*}
    \begin{split}
        &\sum_{k = 1}^p\|\widetilde \Sigma_{jk}^\A - \Sigma_{jk}\|_\cS
        \\=& \sum_{k = 1}^p\|\widetilde \Sigma_{jk}^\A - \Sigma_{jk}\|_\cS I\{\|\widetilde\Lambda_{jk}\|_\cS\geq \lambda\} + \sum_{k = 1}^p\| \Sigma_{jk}\|_\cS I\{\|\widetilde\Lambda_{jk}\|_\cS < \lambda\}
        \\
        \leq &\sum_{k = 1}^p \left\{ \|s_{\lambda}( \widetilde\Lambda_{jk}) - \widetilde\Lambda_{jk} \|_\cS + \|\widetilde\Lambda_{jk}-\widecheck\Lambda_{jk}\|_\cS \right\}\big\|\widetilde \Psi_{jk}^{1/2} \big\|_\infty I\{\|\widetilde\Lambda_{jk}\|_\cS \geq  \lambda, \|\widecheck \Lambda_{jk}\|_\cS \geq  \lambda\}
        \\
        & +\sum_{k = 1}^p  \Big\|\big[s_{\lambda}( \widetilde\Lambda_{jk}) - \widecheck \Lambda_{jk} \big]\widetilde \Psi_{jk}^{1/2} \Big\|_\cS I\{\|\widetilde\Lambda_{jk}\|_\cS \geq  \lambda, \|\widecheck \Lambda_{jk}\|_\cS <  \lambda\} + \sum_{k = 1}^p\| \Sigma_{jk}\|_\cS I\{\|\widecheck \Lambda_{jk}\|_\cS < 2\lambda\}
        \\
        \leq &\sum_{k=1}^p 2\lambda \big\|\widetilde \Psi_{jk}^{1/2} \big\|_\infty I\{  \| \widecheck \Lambda_{jk}\|_\cS \geq  \lambda\} + \sum_{k=1}^p (1+c) {\|\widecheck \Lambda_{jk}\|_\cS} \|\widetilde \Psi_{jk}^{1/2}\|_{\infty}I\{  \|\widecheck \Lambda_{jk}\|_\cS <  \lambda\}
        \\
        &+ \sum_{k=1}^p \|\widecheck \Lambda_{jk}\|_\cS \big\|\widetilde \Psi_{jk}^{1/2} \big\|_\infty I\{  \|\widecheck \Lambda_{jk}\|_\cS <2  \lambda\}
        \\
        \lesssim &~ \lambda^{1-q} \sum_{k=1}^p \big\|\widetilde \Psi_{jk}\big\|_\infty^{1/2} \|\widecheck \Lambda_{jk}\|_\cS^q
        \lesssim \lambda^{1-q}  \sum_{k=1}^p \big\|\Psi_{jk} \big\|_\infty^{(1-q)/2} \left\|\Sigma_{jk}\right\|_\cS^q
      \lesssim  \tilde s_0(p)\left(\frac{\log p}{n^{2\gamma_1}}\right)^{\frac{1-q}{2}}.
    \end{split}
 \end{equation*}
Theorem~\ref{thm_sparse_2} follows.
$\square$

\textbf{Proof of Theorem~\ref{thm_supp_2}}.
Consider two sets: $\widetilde S_{n1} = \{(j,k): \|\widetilde \Sigma_{jk}^\A\|_\cS \neq 0 ~\text{and}~\| \Sigma_{jk}\|_\cS= 0 \}$ and $\widetilde S_{n2} = \{(j,k): \|\widetilde \Sigma_{jk}^\A\|_\cS = 0 ~\text{and}~\| \Sigma_{jk}\|_\cS \neq 0 \}.$ It suffices to prove that
$$
P\big(|\widetilde S_{n1}| >0 \big) + P\big(|\widetilde S_{n2}| >0 \big) \to 0,
$$
as $n,p \to \infty.$
By Conditions (i)-(iii) on $S_{\lambda}(Z),$
$$
\widetilde S_{n1} =  \left\{(j,k): \left\|\frac{\widetilde \Sigma_{jk}}{\widetilde \Psi_{jk}^{1/2}}\right\|_\cS > \lambda ~\text{and}~\| \Sigma_{jk}\|_\cS= 0 \right\} \subset \left\{(j,k): \left\|\frac{\widetilde \Sigma_{jk} - \Sigma_{jk}}{\widetilde \Psi_{jk}^{1/2}}\right\|_\cS > \lambda \right\}
$$
Therefore, with the choice $\lambda = \tilde \delta (\log p/n^{2\gamma_1})^{1/2},$ we obtain
\begin{equation} \label{eq_S1_sparse}
    \begin{split}
        P(|\widetilde S_{n1}|>0) \leq P\left\{\max_{j,k} \left\|\frac{\widetilde \Sigma_{jk} - \Sigma_{jk}}{\widetilde \Psi_{jk}^{1/2}}\right\|_\cS > \lambda\right\}  = o(1),
    \end{split}
\end{equation}
as stated in the proof of Theorem~\ref{thm_sparse_2}.
Similarly, we have
\begin{equation} \nonumber
    \begin{split}
        \widetilde S_{n2} &= \left\{(j,k): \left\|\frac{\widetilde \Sigma_{jk}}{\widetilde \Psi_{jk}^{1/2}}\right\|_\cS \leq \lambda ~\text{and}~\| \Sigma_{jk}\|_\cS \neq 0 \right\}.
    \end{split}
\end{equation}
Note that $\|\Sigma_{jk}\|_\cS \neq 0$ implies that
\begin{equation} \label{eq_S2_con_2}
    \begin{split}
        (2\tilde\delta +\tilde\gamma)\left(\frac{\log p}{n^{2\gamma_1}}\right)^{1/2}< \left\|\frac{  \Sigma_{jk}}{\Psi_{jk}^{1/2}}\right\|_\cS \leq \left[ \left\|\frac{\Sigma_{jk} - \widetilde \Sigma_{jk}}{\widetilde \Psi_{jk}^{1/2}}\right\|_\cS + \left\|\frac{\widetilde \Sigma_{jk}}{\widetilde \Psi_{jk}^{1/2}}\right\|_\cS\right]\left\|\frac{  \widetilde \Psi_{jk}^{1/2}}{\Psi_{jk}^{1/2}}\right\|_\infty.
    \end{split}
\end{equation}
Let $\widetilde \Omega_{n3} = \Big\{\|(\widetilde \Psi_{jk}^{1/2} - \Psi_{jk}^{1/2})/\widetilde \Psi_{jk}^{1/2}\|_\infty\leq\tilde \epsilon \Big\}$ for  some small constant $0 < \tilde\epsilon < \tilde  \gamma/(4\tilde\delta + 2\tilde\gamma).$ By Condition~\ref{cond_sparse_rate2},  $P\{(\widetilde \Omega_{n3})^C\} = o(1).$ Conditioning on the event of $\widetilde \Omega_{n3}$, we can see that $\|\widetilde \Psi_{jk}^{1/2} /\Psi_{jk}^{1/2} \|_\infty \leq 1/(1-\tilde \epsilon).$ This together with (\ref{eq_S2_con_2}) shows that
\begin{equation} \nonumber
    \begin{split}
        \widetilde S_{n2} \cap \widetilde \Omega_{n3} 
      & \subset \left\{(j,k): \left\|\frac{\widetilde \Sigma_{jk} - \Sigma_{jk}}{\widetilde \Psi_{jk}^{1/2}}\right\|_\cS >  \tilde\delta \left(\frac{\log p}{n^{2\gamma_1}}\right)^{1/2} \right\}.
    \end{split}
\end{equation}
As a result,
\begin{equation} \label{eq_S2_sparse}
    \begin{split}
        P(|\widetilde S_{n2}|>0)  \leq P\{(\widetilde \Omega_{n3})^{C}\} + P\left\{\max_{j,k} \left\|\frac{\widetilde \Sigma_{jk} - \Sigma_{jk}}{\widetilde \Psi_{jk}^{1/2}}\right\|_\cS > \tilde\delta \left(\frac{\log p}{n^{2\gamma_1}}\right)^{1/2}\right\} = o(1).
    \end{split}
\end{equation}
Combining (\ref{eq_S1_sparse}) and (\ref{eq_S2_sparse}), we complete our proof.
$\square$

\subsection{Heuristic verification of $I_{jk}$ in (\ref{In}) and Condition~\ref{cond_sparse_rate2}}

\label{supp.verify.In}
In this section we provide the heuristic verification of $I_{jk}$ in (\ref{In}) and Condition~\ref{cond_sparse_rate2} as their detailed proofs are not only long and challenging but also largely deviate from the current focus of the paper.

Recall that
\begin{equation*}\nonumber
    \widetilde \Psi_{jk}   = I_{jk}\sum_{i=1}^n\big( W_{1,jk}V_{00,ijk} +  W_{2,jk}V_{10,ijk}+  W_{3,jk}V_{01,ijk}\big)^2,
\end{equation*}
where, for $a,b=0,1,2,$
\begin{eqnarray*}
&&V_{ab,ijk}(u,v)  = \sum_{i = 1}^{L_{ij}} \sum_{m = 1}^{L_{ik}}g_{ab}\big(h_C,(u,v),(U_{ijl},U_{ikm})\big)\big\{Z_{ijl}Z_{ikm} - \widetilde \Sigma_{jk}(u,v)\big\},\\
&& g_{ab}\big\{h,(u,v),(U_{ijl},U_{ikm})\big\} =   K_{h}({U_{ijl}-u}) K_{h}({U_{ikm}-v})(U_{ijl}-u)^a(U_{ikm}-v)^b.
\end{eqnarray*}
The expression of $\widetilde \Psi_{jk}$ in (\ref{eq.kernel.var}) can be decomposed as
\begin{eqnarray} 
\widetilde \Psi_{jk} &=& I_{jk}W_{1,jk}^2\sum_{i=1}^n V_{00,ijk}^2 + I_{jk}W_{2,jk}^2\sum_{i=1}^nV_{10,ijk}^2+ I_{jk}W_{3,jk}^2\sum_{i=1}^nV_{01,ijk}^2  \nonumber\\
&& + 2I_{jk}W_{1,jk}W_{2,jk}\sum_{i=1}^n V_{00,ijk}V_{10,ijk} +
2I_{jk}W_{1,jk}W_{3,jk}\sum_{i=1}^n V_{00,ijk}V_{01,ijk}  \nonumber\\
&& + 2I_{jk}W_{2,jk}W_{3,jk}\sum_{i=1}^n V_{10,ijk}V_{01,ijk}  \nonumber\\
& = &  \widetilde \Psi_{jk}^{(1)} + \widetilde \Psi_{jk}^{(2)}+ \dots + \widetilde \Psi_{jk}^{(5)} + \widetilde \Psi_{jk}^{(6)}.\label{Psi.sum}
\end{eqnarray}
We first focus on the term $\widetilde \Psi_{jk}^{(1)}.$ 
For $a,b=0,1,2,$ define
\begin{eqnarray*}
V_{ab,ijk}^{(1)}(u,v)  &=& \sum_{l = 1}^{L_{ij}} \sum_{m = 1}^{L_{ik}}g_{ab}\big(h_C,(u,v),(U_{ijl},U_{ikm})\big)\big\{Z_{ijl}Z_{ikm} - \Sigma_{jk}(u,v)\big\},\\
V_{ab,ijk}^{(2)}(u,v)  &=& \sum_{l = 1}^{L_{ij}} \sum_{m = 1}^{L_{ik}}g_{ab}\big(h_C,(u,v),(U_{ijl},U_{ikm})\big)\big\{\Sigma_{jk}(u,v) - \widetilde \Sigma_{jk}(u,v)\big\}.
\end{eqnarray*}
The term $\widetilde \Psi_{jk}^{(1)}$ can be re-expressed as
\begin{eqnarray}
\nonumber
\widetilde \Psi_{jk}^{(1)} &=& I_{jk}W_{1,jk}^2\sum_{i=1}^n  \big\{V_{00,ijk}^{(1)}\big\}^2+I_{jk}W_{1,jk}^2\sum_{i=1}^n  \big\{V_{00,ijk}^{(2)}\big\}^2+2I_{jk}W_{1,jk}^2\sum_{i=1}^n V_{00,ijk}^{(1)}V_{00,ijk}^{(2)}\\
&=& D_{jk,1} + D_{jk,2}+D_{jk,3}.
\label{D.terms}
\end{eqnarray}

(a) {\it Verification of $I_{jk}$}. To show the rationale of imposing the rate $I_{jk}$ in (\ref{In}), we need to verify that  
\begin{equation}
\label{In_rate}
    I_{jk}\sum_{i=1}^n\big( W_{1,jk}V_{00,ijk} +  W_{2,jk}V_{10,ijk}+  W_{3,jk}V_{01,ijk}\big)^2 \asymp  1+o_P(1).
\end{equation}
for each $u,v \in \cU.$ 
Denote by $\tilde{n}_{jk} = \sum_{i=1}^n L_{ij} L_{ik}.$ Recall that
\begin{equation*}
S_{ab,jk}(u,v) =  \sum_{i = 1}^n \sum_{l = 1}^{L_{ij}} \sum_{m = 1}^{L_{ik}} g_{ab}\big\{h_C,(u,v),(U_{ijl},U_{ikm})\big\} .
\end{equation*}
It can be shown that
$S_{ab,jk}(u,v) \asymp \tilde{n}_{jk}h_C^{a+b} \{1 + o_P(1)\}$ for $a,b=0,1,2,$ which together with (\ref{eq.kernel.w}) implies that
\begin{equation}
\label{W_rate}
W_{1,jk}(u,v) \asymp \tilde{n}_{jk}^{-1}  \big\{1 + o_P(1)\big\},~~W_{2,jk}(u,v) \asymp W_{3,jk}(u,v) \asymp \tilde{n}_{jk}^{-1} h_C^{-1} \big\{1 + o_P(1)\big\}.    
\end{equation}
Similarly, we can also show that
\begin{eqnarray}
\label{G_sqaure_rate}
\nonumber
&&\sum_{i=1}^n \Big\{\sum_{l = 1}^{L_{ij}} \sum_{m = 1}^{L_{ik}}g_{00}\big(h_C,(u,v),(U_{ijl},U_{ikm})\big)\Big\}^2  \\ \nonumber
& = & \sum_{i=1}^n \sum_{l = 1}^{L_{ij}} \sum_{m = 1}^{L_{ik}} K_{h_C}^2({U_{ijl}-u}) K_{h_C}^2({U_{ikm}-v})  \\ \nonumber
&   & + \sum_{i=1}^n \sum_{l = 1}^{L_{ij}} \sum_{m' \neq m}^{L_{ik}} K_{h_C}^2({U_{ijl}-u}) K_{h_C}({U_{ikm}-v})K_{h_C}({U_{ikm'}-v}) \\ \nonumber
&   & + \sum_{i=1}^n \sum_{l \neq l'}\sum_{m = 1}^{L_{ik}} K_{h_C}({U_{ijl}-u})K_{h_C}({U_{ijl'}-u}) K_{h_C}^2({U_{ikm}-v}) \\ \nonumber
&   & + \sum_{i=1}^n \sum_{l \neq l'}^{L_{ij}} \sum_{m \neq m'}^{L_{ik}}  K_{h_C}({U_{ijl}-u}) K_{h_C}({U_{ikm}-v})K_{h_C}({U_{ijl'}-u}) K_{h_C}({U_{ikm'}-v})\\ 
&\asymp& \Big\{\sum_{i=1}^n \big(L_{ij} L_{ik} h_C^{-2} +  L_{ij}^2 L_{ik} h_C^{-1} + L_{ij} L_{ik}^2 h_C^{-1} +  L_{ij}^2 L_{ik}^2\big)\Big\}\big\{1+o_P(1)\big\}. 
\end{eqnarray}

By (\ref{W_rate}) and (\ref{G_sqaure_rate}), we obtain that
\begin{equation}
\label{W.rate2}
W_{1,jk}^2 \sum_{i=1}^n \Big\{\sum_{l = 1}^{L_{ij}} \sum_{m = 1}^{L_{ik}}g_{00}\big(h_C,(u,v),(U_{ijl},U_{ikm})\big)\Big\}^2 \asymp I_{jk}^{-1} \big\{ 1 + o_P(1)\big\},
\end{equation}
which together with  $\Sigma_{jk}(u,v) - \widetilde \Sigma_{jk}(u,v) = o_P(1)$ implies that
$D_{jk,2} = o_P(1)$ and $D_{jk,3} = o_P(1).$
Note that $\mathbb{E}\big\{Z_{ijl}Z_{ikm} - \Sigma_{jk}(u,v)\big\}^2$ is bounded. Together with (\ref{W_rate}) and (\ref{G_sqaure_rate}), we can also show that 
$
D_{jk,1}(u,v) \asymp 1+o_P(1).
$
Combining the above results yields that
$\widetilde \Psi_{jk}^{(1)} \asymp 1+o_P(1).$
In a similar fashion, we can also show that $\widetilde \Psi_{jk}^{(i)} \asymp 1+o_P(1)$ for $i = 2,\ldots,6$ in (\ref{Psi.sum}) and hence (\ref{In_rate}) follows. 

(b) {\it Verification of Condition~\ref{cond_sparse_rate2}}. To verify the uniform convergence rate in Condition~\ref{cond_sparse_rate2}, we need to refine our analysis above to construct the exponential type of tail bounds on $\widetilde{\Psi}_{jk}(u,v)-{\Psi}_{jk}(u,v)$ at each $(u,v) \in \cU^2$ rather than the consistency results in (a). 

Consider the first term $D_{jk,1}(u,v)=(\tilde n_{jk}W_{1,jk})^2 \times I_{jk} \tilde n_{jk}^{-2}\sum_{i=1}^n \{V_{00,ijk}^{(1)}\}^2$ in (\ref{D.terms}). Note that by (\ref{W_rate}) $\tilde{n}_{jk}|W_{1,jk}|$ is bounded with an overwhelming probability. Suppose that $X_{ij}(\cdot)$'s are sub-Gaussian processes and $\varepsilon_{ijl}$'s are independent sub-Gaussian errors. Since $\{V_{00,ijk}^{(1)}(u,v), i =1,\ldots,n\}$ forms an independent sequence, we can obtain the tail bound on $D_{jk,1}(u,v) - \mathbb{E} \{D_{jk,1}(u,v)\}$ by calculating all $q$-th moments of $\zeta_{ijk} = \big\{V_{00,ijk}^{(1)}(u,v)\big\}^2 - \mathbb{E}\big[\{V_{00,ijk}^{(1)}(u,v)\}^2\big]$ for $q = 2,3,4,\dots$ under regularity conditions. Since $\zeta_{ijk}$'s are either sub-Gaussian or sub-exponential, we can follow the similar techniques to prove Lemma~5 of \cite{qiao2020} by adopting a truncation technique and then applying Bernstein inequality (\textcolor{blue}{Boucheron et al., 2014}) to establish a rough exponential type of concentration inequality (i.e., the equipped tail bound is in the same form of the exponential tail bound in (\ref{coneq.L2})) for $D_{jk,1}(u,v)$ at each $(u,v) \in \cU^2.$ Similarly, we can also derive the exponential type of concentration inequality for the third term  $D_{jk,3}(u,v)$ in (\ref{D.terms}).

Consider the second term $D_{jk,2}(u,v)$ in (\ref{D.terms}), which can be re-expressed as
$$
D_{jk,2}(u,v)=2I_{jk}W_{1,jk}^2\sum_{i=1}^n \Big\{\sum_{l = 1}^{L_{ij}} \sum_{m = 1}^{L_{ik}}g_{ab}\big(h_C,(u,v),(U_{ijl},U_{ikm})\big)\Big\}^2\Big\{\Sigma_{jk}(u,v) - \widetilde \Sigma_{jk}(u,v)\Big\}^2. 
$$
Note that it follows from (\ref{W.rate2}) that $I_{jk}W_{1,jk}^2\sum_{i=1}^n \big\{\sum_{l = 1}^{L_{ij}} \sum_{m = 1}^{L_{ik}}g_{ab}(h_C,(u,v),(U_{ijl},U_{ikm}))\big\}^2$ is bounded with an overwhelming probability. Then the exponential type of concentration bound on $D_{jk,2}(u,v)$ at each $(u,v) \in \cU^2$ can be obtained through the exponential type tail bound on $\Sigma_{jk}(u,v) - \widetilde \Sigma_{jk}(u,v)$, which has been established in \cite{qiao2020}, see details in proofs of its Lemmas~4 and 5 under the sparse and dense designs, respectively.

To derive the uniform (i.e., over $\cU^2$) concentration inequality for $\widetilde \Psi_{jk}^{(1)}(u,v)$ in (\ref{D.terms}), we can apply the partition technique that reduces the problem from supremum over $\cU^2$ to the maximum over a grid of pairs and then follow the similar developments to prove the uniform concentration inequalities in Lemmas~4 and 5 of \cite{qiao2020}. In a similar fashion to the above procedure, we can develop the corresponding exponential type of uniform concentration inequality for $\widetilde \Psi_{jk}^{(i)}(u,v)$ for $i=2, \dots, 6.$ As a result, the exponential type of uniform concentration inequality for $\widetilde \Psi_{jk}(u,v)$ can be obtained. 

The uniform convergence rate in Condition~\ref{cond_sparse_rate2} is implied by the exponential type of uniform concentration inequalities for $\widetilde \Psi_{jk}(u,v)$ for each $j,k,$ which partially depend on the uniform concentration bounds on $\widetilde\Sigma_{jk}(u,v)$'s. 
In a similar spirit to the $L_2$ concentration bounds on $\widetilde \Sigma_{jk}(u,v)$'s implied by Condition~\ref{cond_sparse_rate1}, we consider the uniform convergence rate of $\widetilde \Sigma_{jk}(u,v),$
\begin{equation}
\label{coneq_uniform1}
\max_{1 \le j,k \le p } \underset{u,v \in \cU}{\sup}\Big|\widetilde{\Sigma}_{jk}(u,v) - \Sigma_{jk}(u,v)\Big| =O_P\left(\sqrt{\frac{\log p}{n^{2\gamma_1}}} + h^2 \right),
\end{equation}
which is satisfied if there exists some positive constants $c_i$ for $i=6, \dots, 9$ and $\gamma_1 \in (0,1/2]$ such that for each $j,k = 1, \dots, p$ and $t \in (0,1],$ 
\begin{equation}
\label{coneq_uniform2}
   P\Big\{{\sup}_{u,v \in \cU}|\widetilde{\Sigma}_{jk}(u,v) - \Sigma_{jk}(u,v)| \ge t  + c_8 h^2\Big\} \leq c_7 n^{c_9} \exp( - {c_6n^{2\gamma_1} t^2}). 
\end{equation}
Larger values of $\gamma_1$ correspond to a more frequent measurement schedule and hence faster rate in (\ref{coneq_uniform1}).
For sparsely sampled functional data, it follows from Lemma~4 of \cite{qiao2020} and the same proof technique for $j \neq k$ that (\ref{coneq_uniform2}) holds by choosing $\gamma_1=1/2-a$ and $c_9=1+2a$ with $h \asymp n^{-a}$ for some positive constant $a<1/2.$
For densely sampled functional data, it follows from Lemma~5 of \cite{qiao2020} and more efforts for $j \neq k$ that (\ref{coneq_uniform2}) holds with the choice of $\gamma_1=\min(1/2,1/3+b/6-\epsilon'/2-2a/3)$ and $c_9=\max(1,2/3-\epsilon'-b/3+4a/3)$ for some small constant $\epsilon'>0$ when $h \asymp n^{-a}$ and $L \asymp n^b$ for some positive constants $a,b.$

Following the proof procedure described above, we can establish exponential type of uniform concentration inequality for $\widetilde \Psi_{jk}(u,v)$ for each $j,k$ in the same form as (\ref{coneq_uniform2}) but with different positive constants and in particular $\gamma_2 \in (0,1/2],$ which will result in the uniform convergence rate in Condition~\ref{cond_sparse_rate2}. It is worth mentioning that such heuristic analysis can only help us establish uniform concentration inequalities for $\widetilde \Psi_{jk}(u,v)$'s leading to the sub-optimal rate. Investigating the corresponding optimal rate through the precise specification of the largest values of $\gamma_2$ under different measurement schedules or more generally through $n, h$ and possibly $L$ for the dense case is quite challenging and remains an open topic to be pursued in the future.

\section{Additional empirical results}
\label{supp.sec_emp}

\subsection{Simulation studies}
\subsubsection{Fully observed functional data}

\begin{figure}[!htbp]
\centering
\begin{subfigure}{.3\linewidth}
    \centering
    \centering
    \includegraphics[width=5cm]{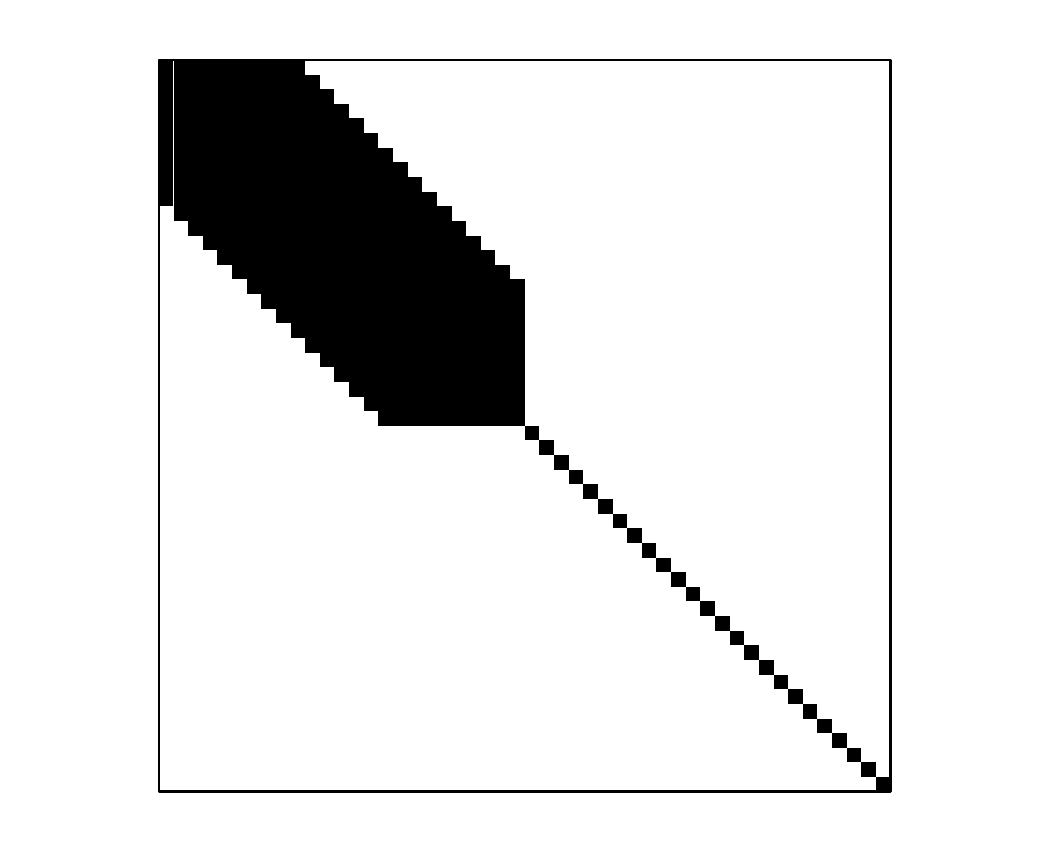}
    \caption{True}
\end{subfigure}
\begin{subfigure}{.3\linewidth}
    \centering
    \includegraphics[width=5cm]{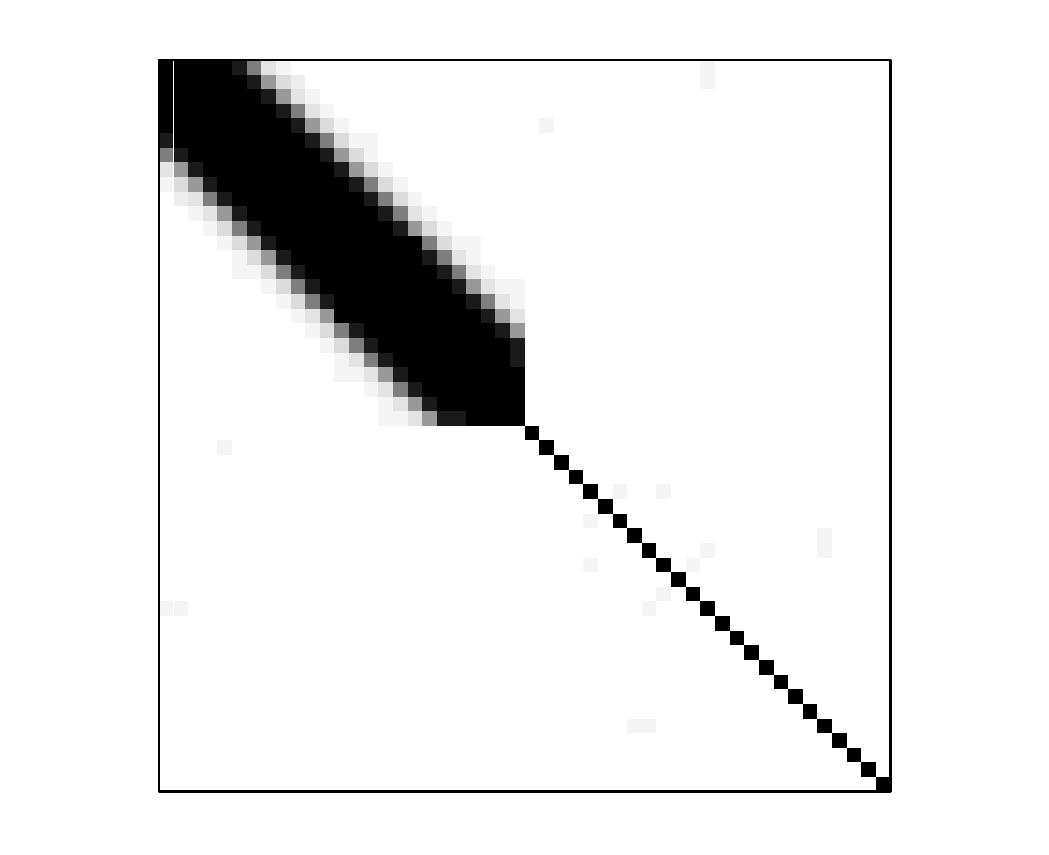}
    \caption{Hard $\wbSigma_\A$}
\end{subfigure}
\begin{subfigure}{.3\linewidth}
    \centering
    \includegraphics[width=5cm]{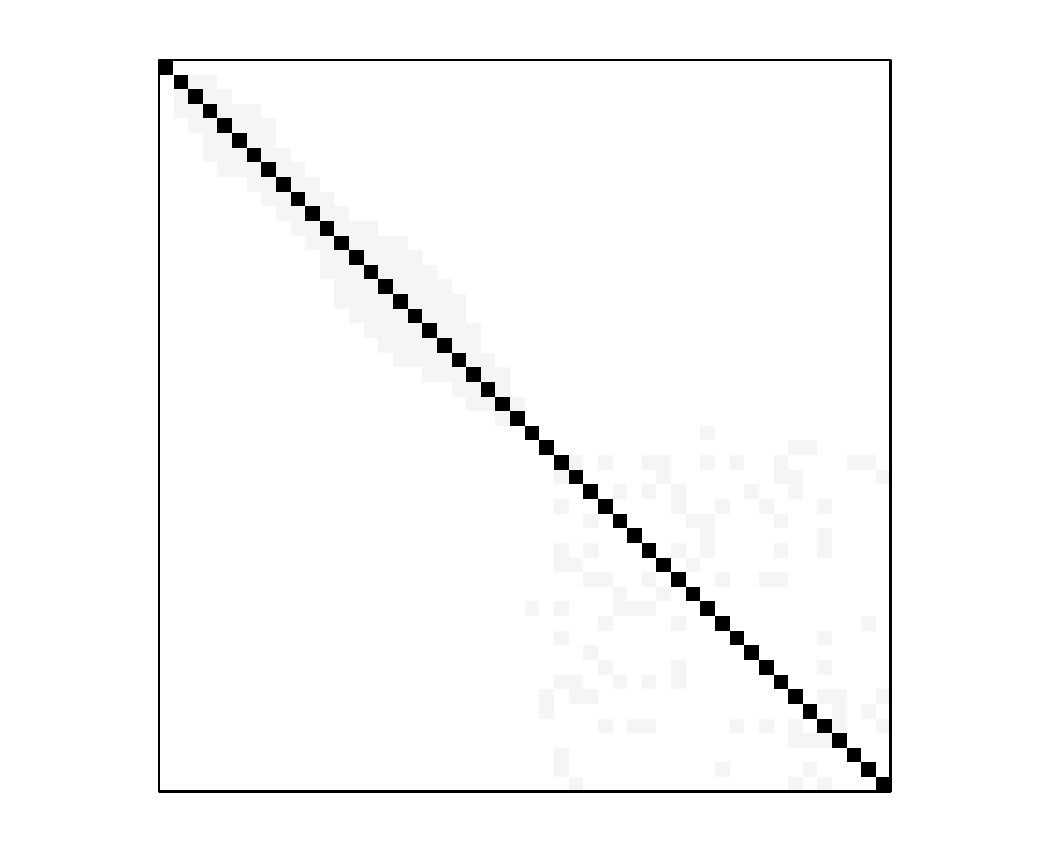}
    \caption{Hard $\wbSigma_\U$}
\end{subfigure}

\begin{subfigure}{.3\linewidth}
    \centering
    \centering
    \includegraphics[width=5cm]{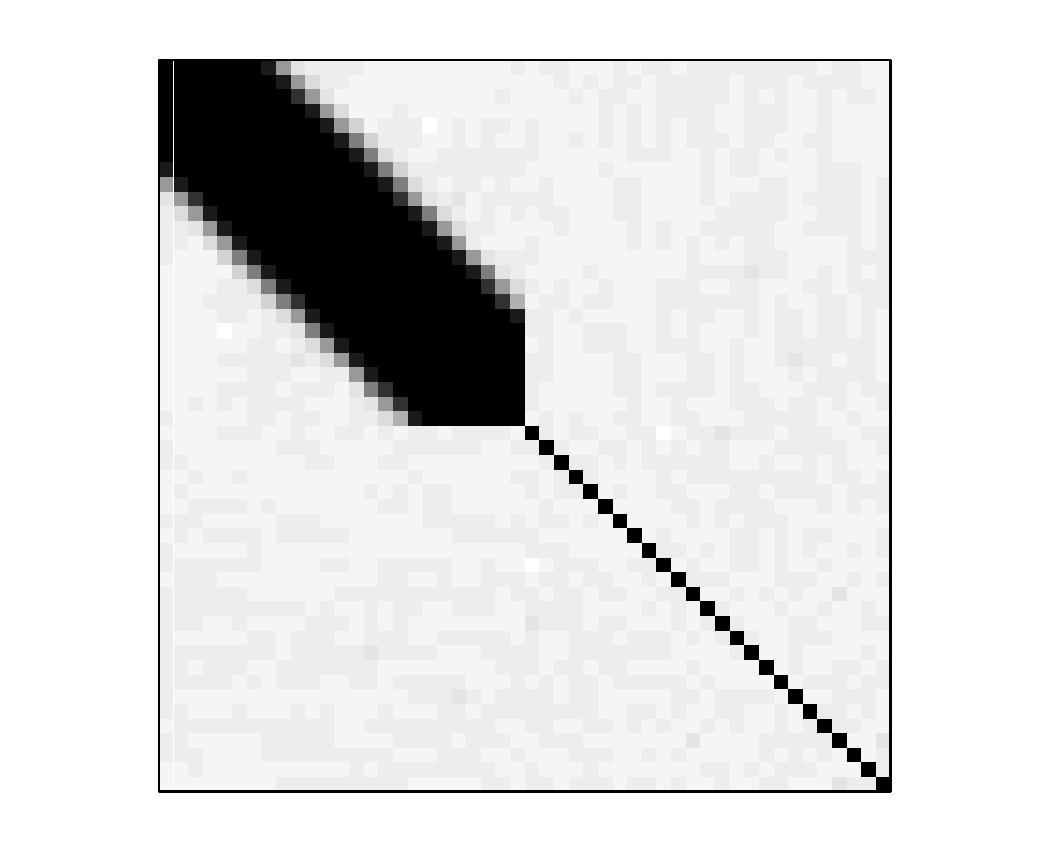}
    \caption{Soft $\wbSigma_\A$}
\end{subfigure}
\begin{subfigure}{.3\linewidth}
    \centering
    \includegraphics[width=5cm]{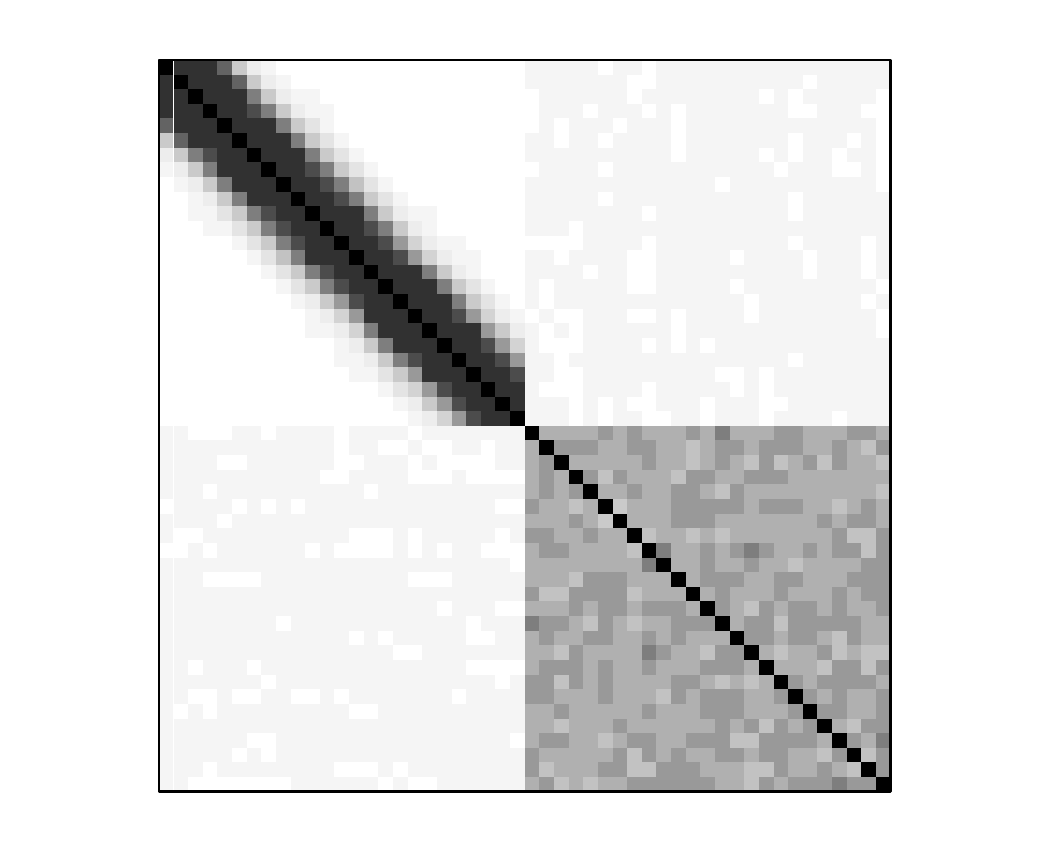}
    \caption{Soft $\wbSigma_\U$}
\end{subfigure}
\begin{subfigure}{.3\linewidth}
    \centering
    \includegraphics[width=5cm]{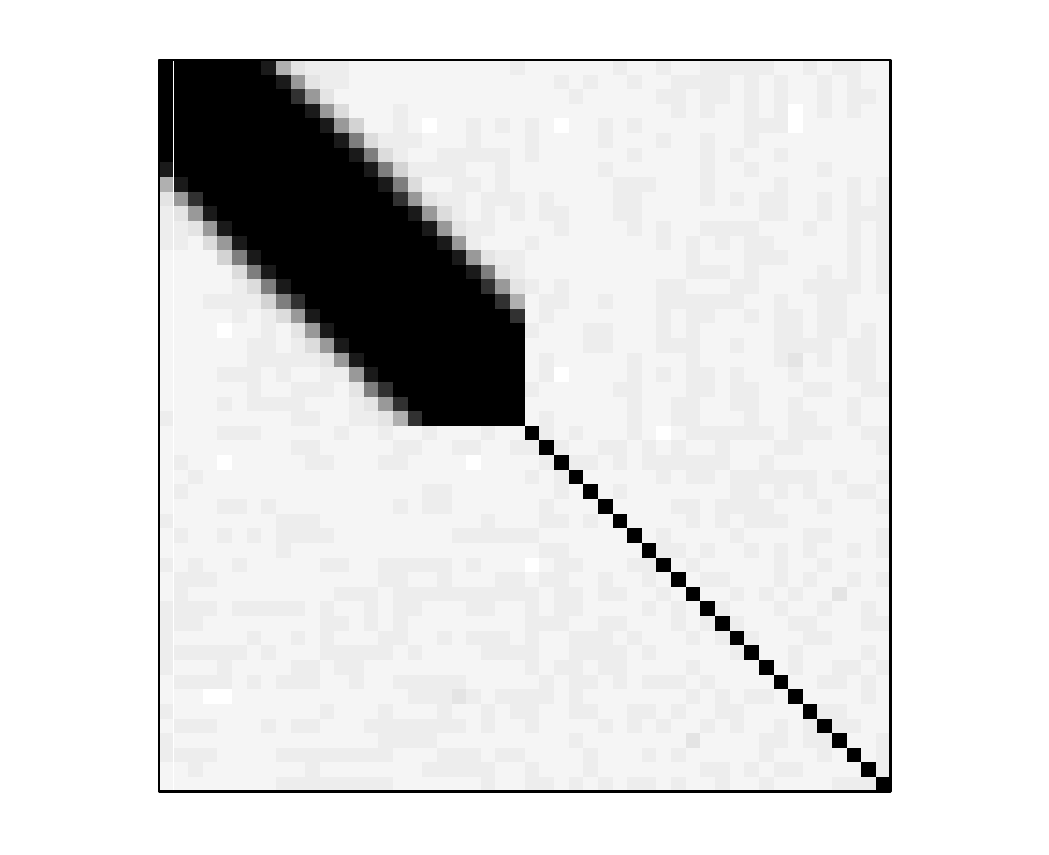}
    \caption{SCAD $\wbSigma_\A$}
\end{subfigure}

\begin{subfigure}{.3\linewidth}
    \centering
    \centering
    \includegraphics[width=5cm]{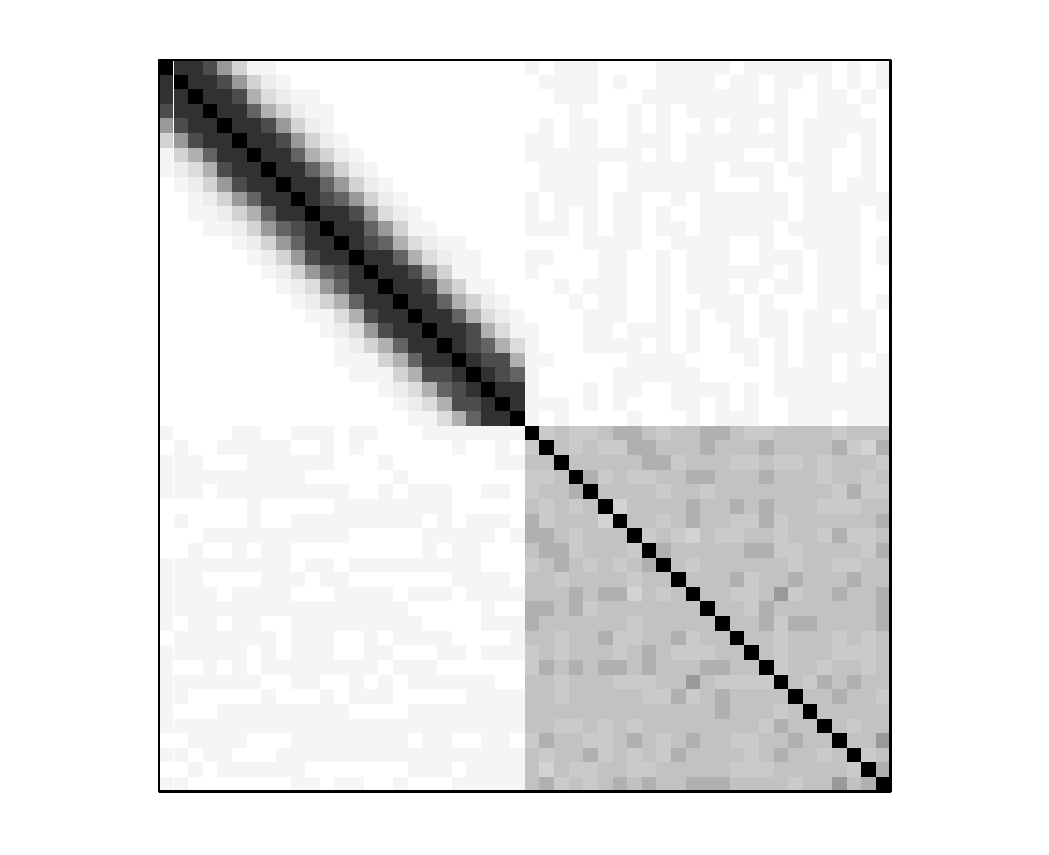}
    \caption{SCAD $\wbSigma_\U$}
\end{subfigure}
\begin{subfigure}{.3\linewidth}
    \centering
    \includegraphics[width=5cm]{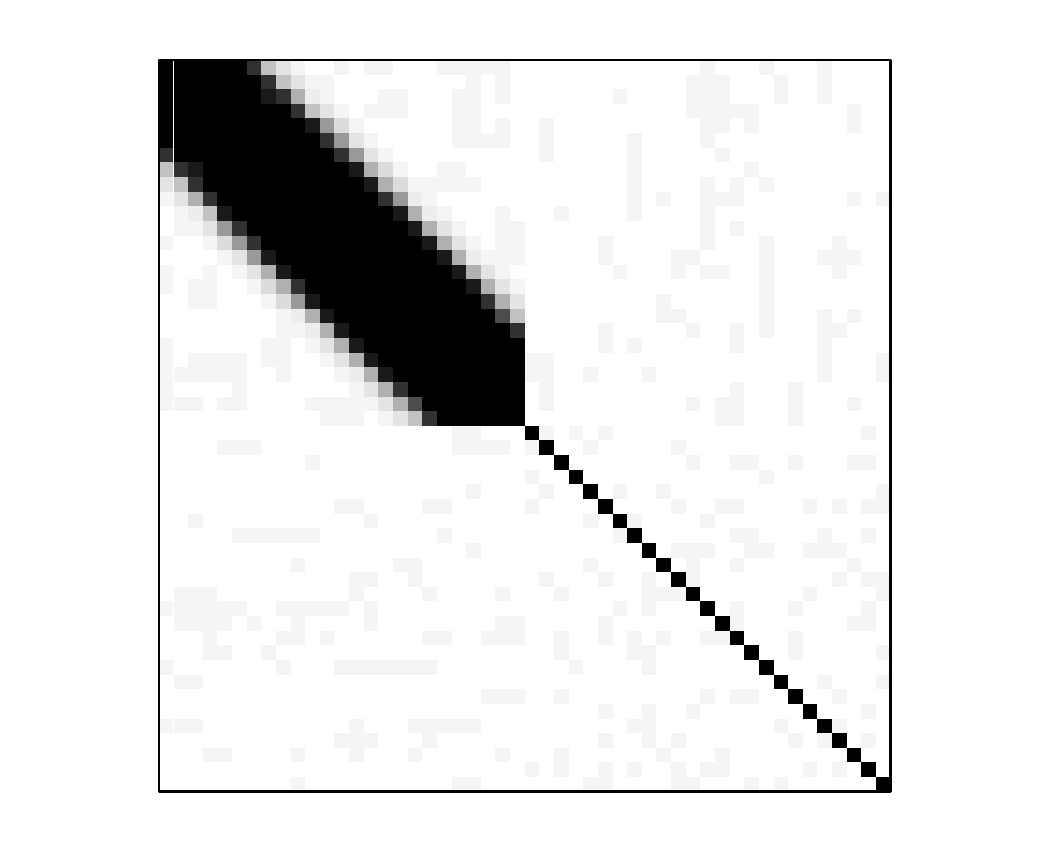}
    \caption{Adap. lasso $\wbSigma_\A$}
\end{subfigure}
\begin{subfigure}{.3\linewidth}
    \centering
    \includegraphics[width=5cm]{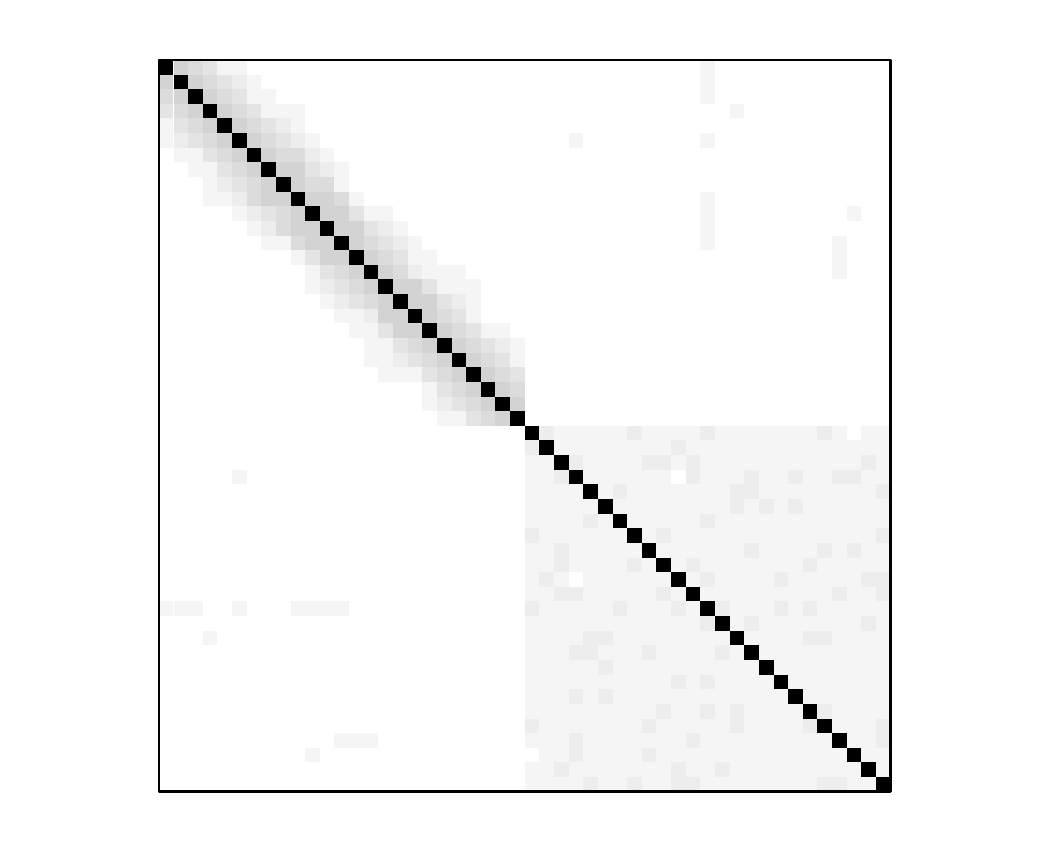}
    \caption{Adap. lasso $\wbSigma_\U$}
\end{subfigure}
\centering
\caption{\label{hm_band}{ Heat maps of the frequency of the zeros identified for the Hilbert--Schimidt norm of each entry of the estimated covariance function (when $p$ = 50) for Model 1 out of 100 simulation runs. White and black correspond to 100/100 and 0/100 zeros identified, respectively.  }}
\end{figure}

\begin{figure}[!htbp]
\centering
\begin{subfigure}{.3\linewidth}
    \centering
    \centering
    \includegraphics[width=5cm]{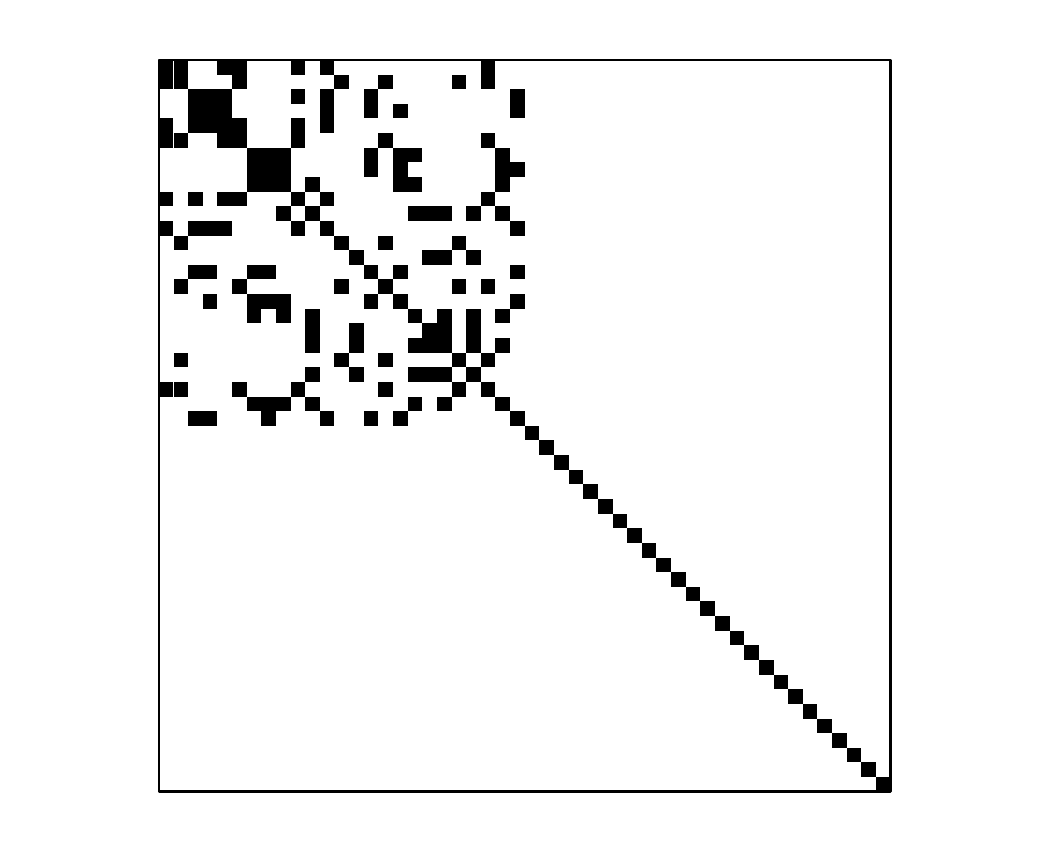}
    \caption{True}
\end{subfigure}
\begin{subfigure}{.3\linewidth}
    \centering
    \includegraphics[width=5cm]{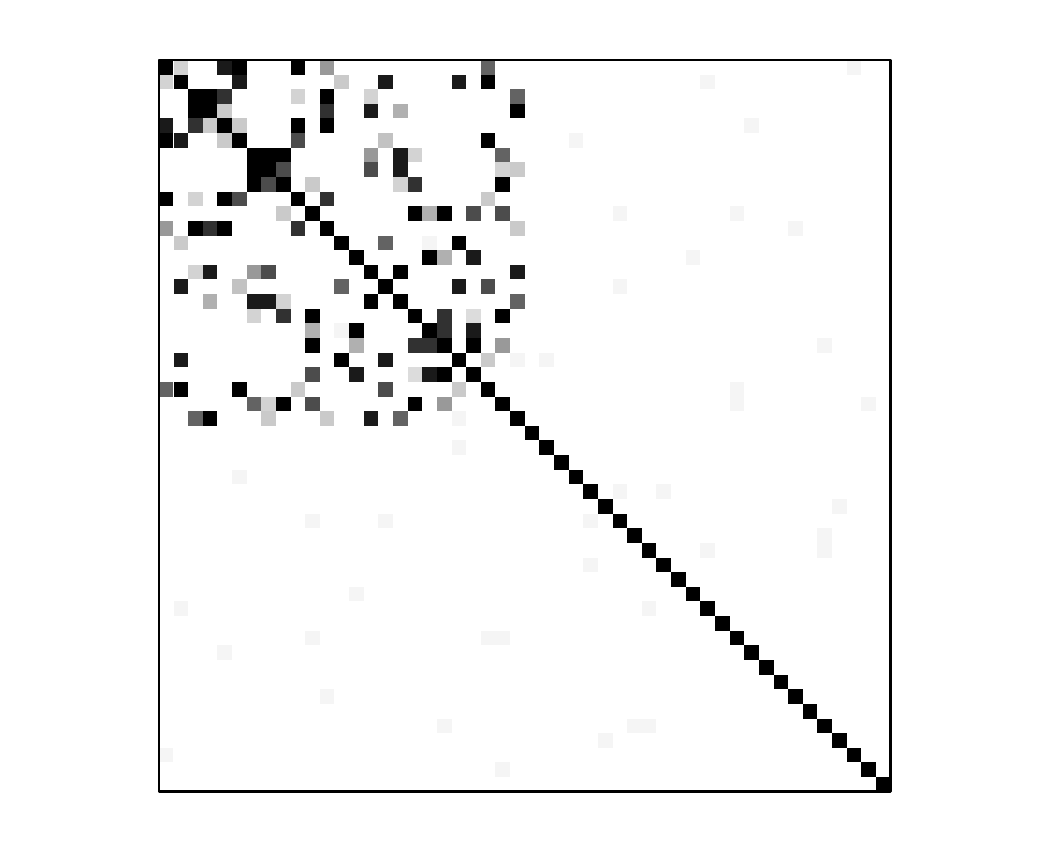}
    \caption{Hard $\wbSigma_\A$}
\end{subfigure}
\begin{subfigure}{.3\linewidth}
    \centering
    \includegraphics[width=5cm]{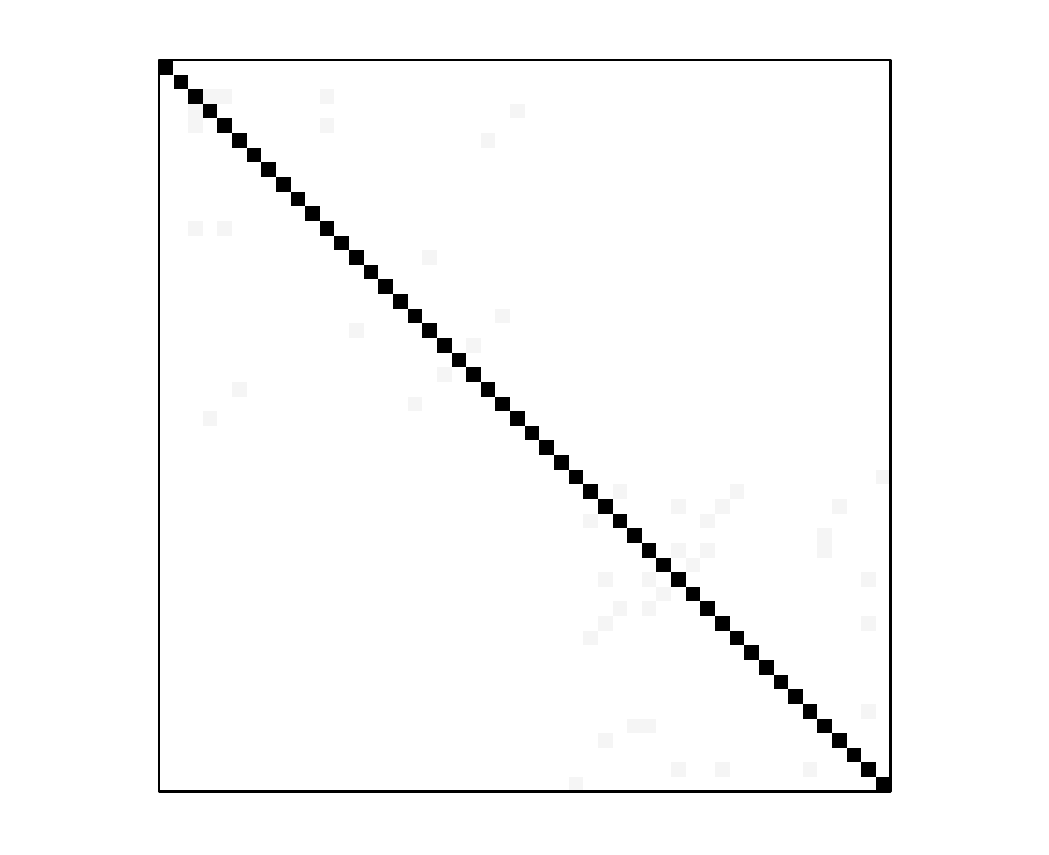}
    \caption{Hard $\wbSigma_\U$}
\end{subfigure}

\begin{subfigure}{.3\linewidth}
    \centering
    \centering
    \includegraphics[width=5cm]{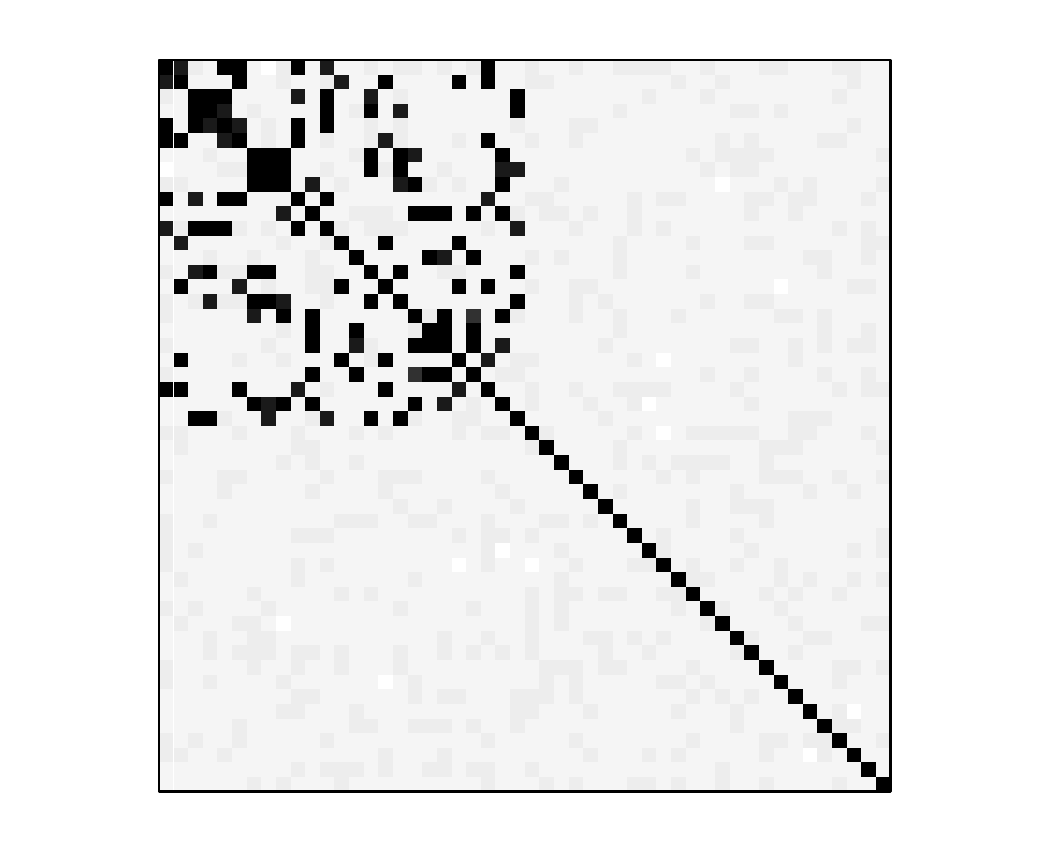}
    \caption{Soft $\wbSigma_\A$}
\end{subfigure}
\begin{subfigure}{.3\linewidth}
    \centering
    \includegraphics[width=5cm]{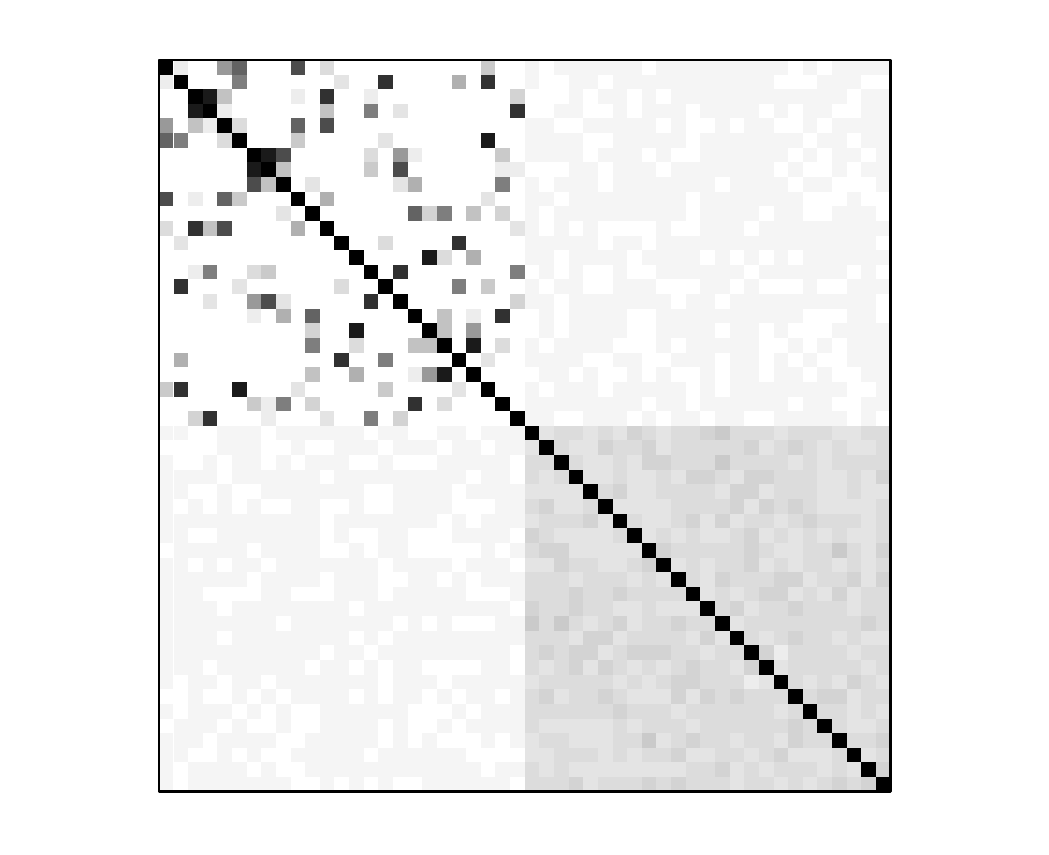}
    \caption{Soft $\wbSigma_\U$}
\end{subfigure}
\begin{subfigure}{.3\linewidth}
    \centering
    \includegraphics[width=5cm]{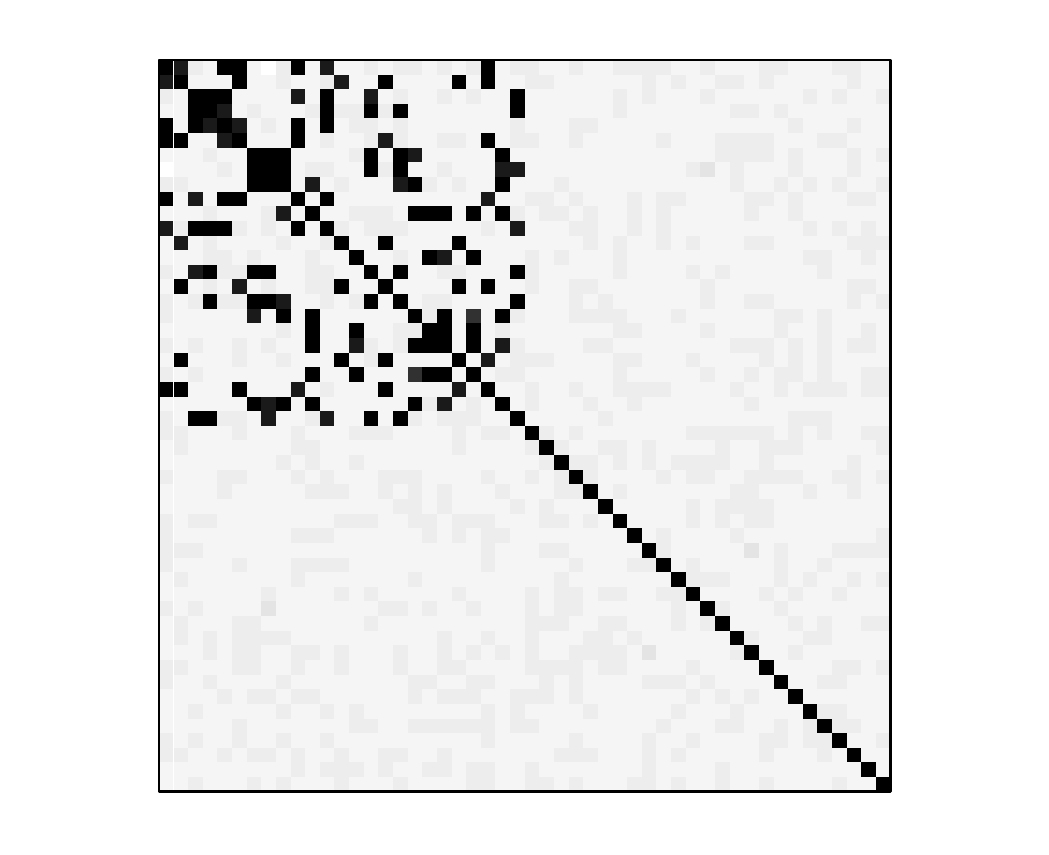}
    \caption{SCAD $\wbSigma_\A$}
\end{subfigure}

\begin{subfigure}{.3\linewidth}
    \centering
    \centering
    \includegraphics[width=5cm]{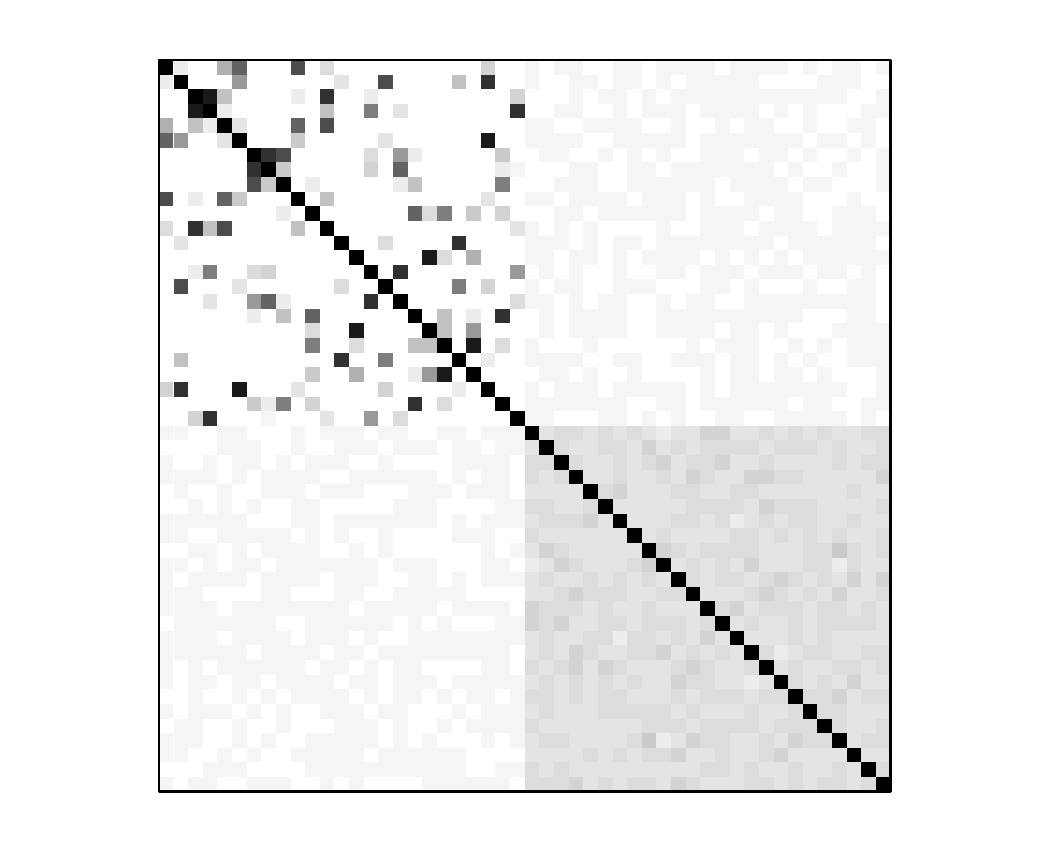}
    \caption{SCAD $\wbSigma_\U$}
\end{subfigure}
\begin{subfigure}{.3\linewidth}
    \centering
    \includegraphics[width=5cm]{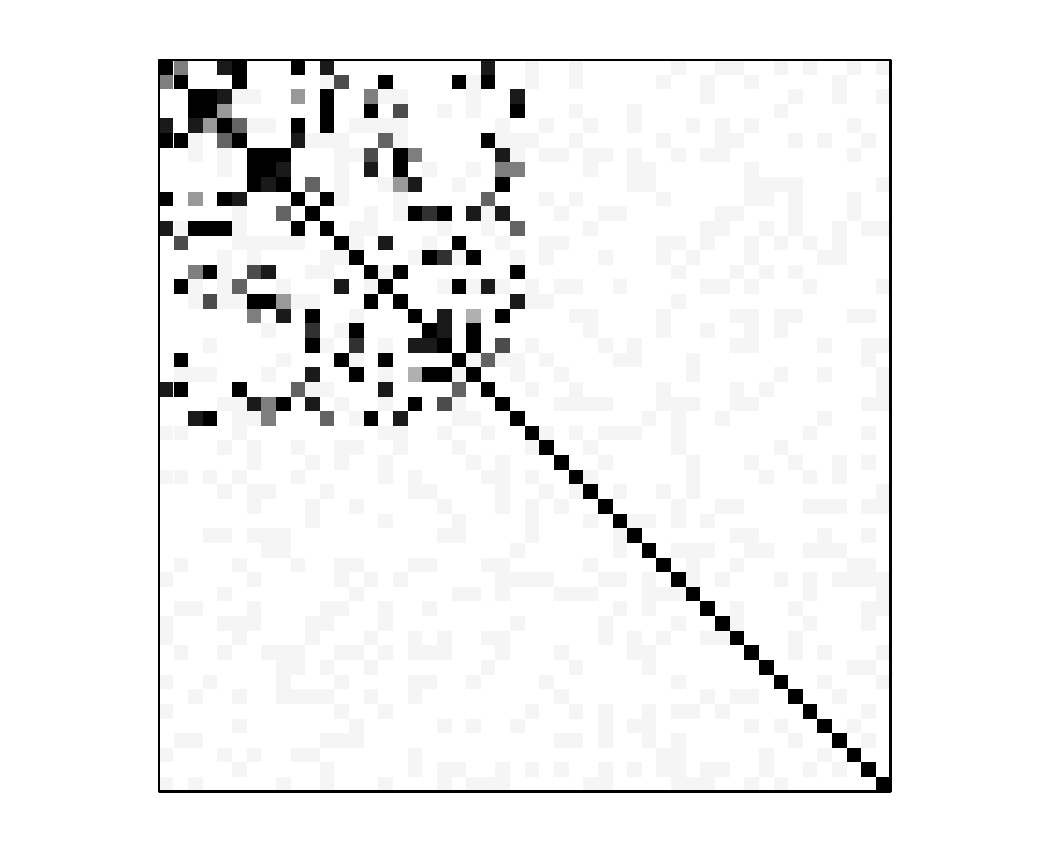}
    \caption{Adap. lasso $\wbSigma_\A$}
\end{subfigure}
\begin{subfigure}{.3\linewidth}
    \centering
    \includegraphics[width=5cm]{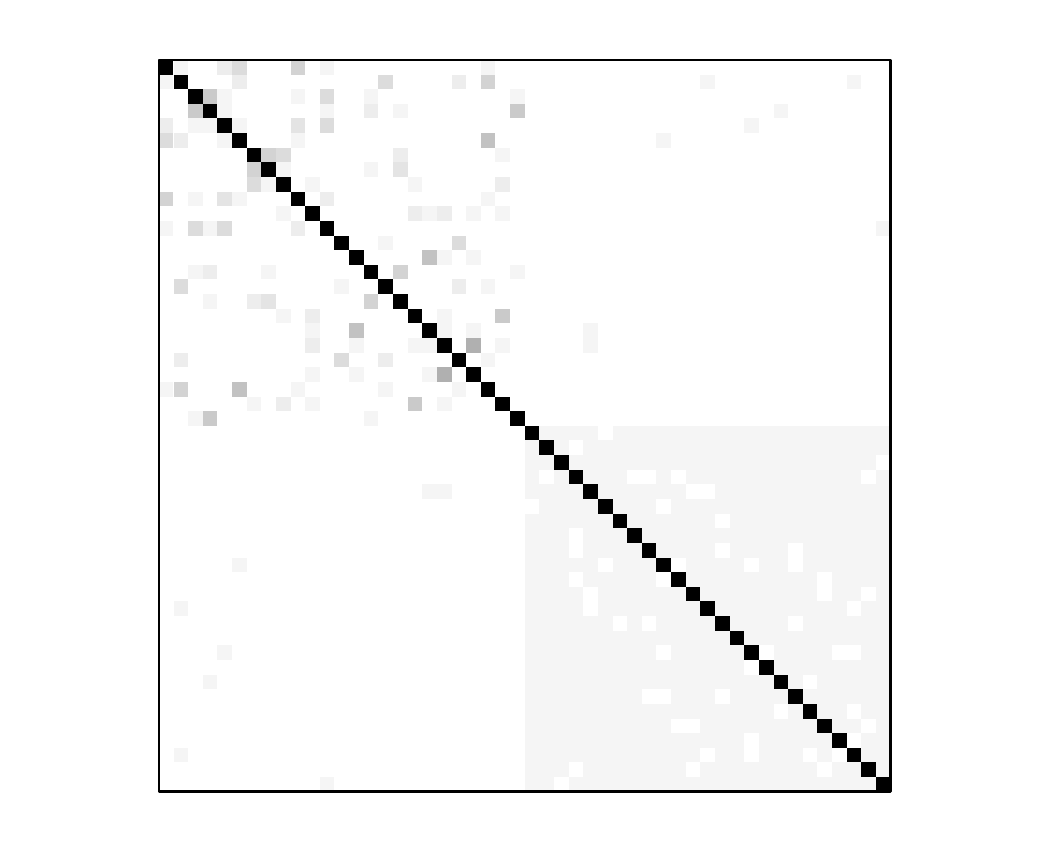}
    \caption{Adap. lasso $\wbSigma_\U$}
\end{subfigure}
\centering
\caption{\label{hm_ran}{ Heat maps of the frequency of the zeros identified for the Hilbert--Schimidt norm of each entry of the estimated covariance function (when $p$ = 50) for Model 2 out of 100 simulation runs. White and black correspond to 100/100 and 0/100 zeros identified, respectively. }}
\end{figure}

\begin{figure}[!htbp]
	\centering
	\includegraphics[width=1\textwidth]{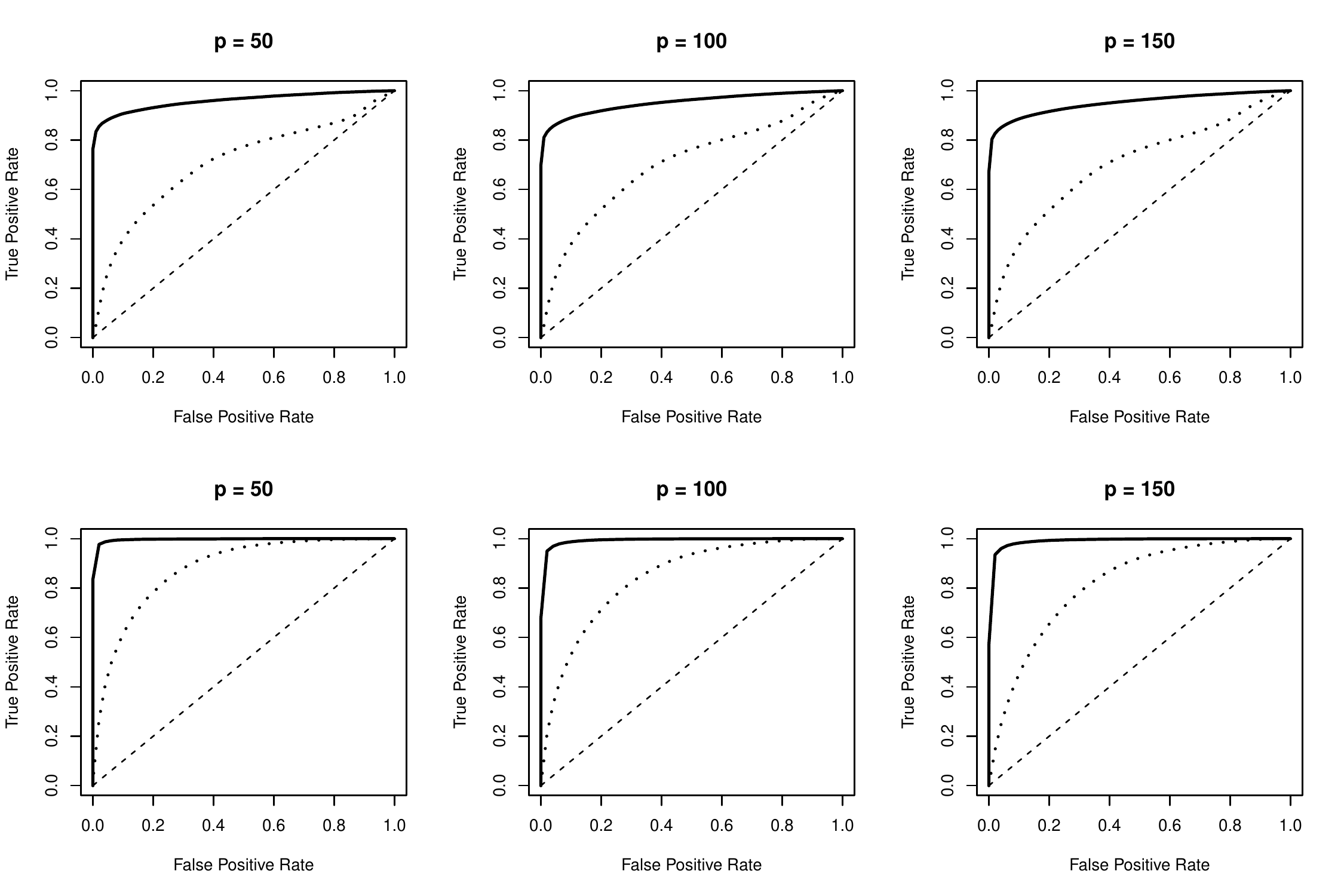}
	\caption{\label{ROC_n_100}{Model 1 (top row) and Model 2 (bottom row) for $p=50,100,150$: Comparison of the average ROC curves for adaptive functional thresholding (solid line) and universal functional thresholding (dotted line) over 100 simulation runs. }}
\end{figure}

Figures~\ref{hm_band} and \ref{hm_ran} plot the heat maps
 of the frequency of the zeros identified for the Hilbert--Schimidt norm of each entry of the estimated covariance function, when $p = 50,$ out of 100 simulation runs.  The true nonzero patterns of Model 1 and 2 are presented in Figures \ref{hm_band}(a) and \ref{hm_ran}(a), respectively. Figure \ref{ROC_n_100} displays the average receiver operating characteristic (ROC) curves (plots of true positive rates versus false positive rates over a sequence of $\lambda$ values) for both the adaptive functional thresholding and universal functional thresholding methods. These results again demonstrate the uniform superiority of the adaptive functional thresholding method in terms of graph selection consistency.

\subsubsection{Partially observed functional data}

Tables~\ref{err.table.kernel.100} and \ref{TPR/FPR.table.kernel.100} summarize the estimation and support recovery performance of BinLLS-based adaptive and universal functional thresholding estimators for the setting of $p = 100$ satisfying Models~1 and 2 under different measurement schedules. The same patterns as those from Tables~\ref{err.table.kernel.50} and \ref{TPR/FPR.table.kernel.50} can be observed. 

\begin{table}[!htbp]
	\caption{\label{err.table.kernel.100}
	The average (standard error) functional matrix losses for partially observed functional scenarios and $p=100$ over 100 simulation runs.}
	\begin{center}
		\resizebox{6.2in}{!}{
			\begin{tabular}{ccrrrrrrrr}
	\hline																			
	&		&	\multicolumn{2}{c}{$L_i = 11$}			&	\multicolumn{2}{c}{$L_i = 21$}			&	\multicolumn{2}{c}{$L_i = 51$}			&	\multicolumn{2}{c}{$L_i = 101$}			\\	
{Model} 	&	{Method}	&	\multicolumn{1}{c}{$\widecheck \bSigma_\A$}	&	\multicolumn{1}{c}{$\widecheck \bSigma_\U$}	&	\multicolumn{1}{c}{$\widecheck \bSigma_\A$}	&	\multicolumn{1}{c}{$\widecheck \bSigma_\U$}	&	\multicolumn{1}{c}{$\widecheck \bSigma_\A$}	&	\multicolumn{1}{c}{$\widecheck \bSigma_\U$}	&	\multicolumn{1}{c}{$\widecheck \bSigma_\A$}	&	\multicolumn{1}{c}{$\widecheck \bSigma_\U$}	\\	\hline
\multirow{12}{*}{1}	&	 \multicolumn{9}{c}{Functional Frobenius norm}																\\		
	&	Hard	&	11.40(0.03)	&	18.34(0.01)	&	9.63(0.03)	&	17.80(0.01)	&	8.55(0.04)	&	17.51(0.01)	&	8.17(0.04)	&	17.42(0.01)	\\	
	&	Soft	&	12.79(0.05)	&	18.33(0.01)	&	11.28(0.05)	&	17.71(0.02)	&	10.33(0.05)	&	16.68(0.07)	&	10.01(0.05)	&	16.06(0.07)	\\	
	&	SCAD	&	12.41(0.05)	&	18.33(0.01)	&	10.58(0.05)	&	17.72(0.02)	&	9.42(0.05)	&	16.77(0.06)	&	9.01(0.05)	&	16.23(0.07)	\\	
	&	Adap. lasso	&	11.22(0.04)	&	18.33(0.01)	&	9.59(0.04)	&	17.79(0.01)	&	8.54(0.04)	&	17.49(0.01)	&	8.19(0.04)	&	17.34(0.03)	\\	
\cline{2-10}																		
	&	 \multicolumn{9}{c}{Functional matrix $\ell_1$ norm}																	\\	
	&	Hard	&	5.97(0.05)	&	9.41(0.01)	&	5.15(0.05)	&	9.35(0.01)	&	4.70(0.05)	&	9.33(0.01)	&	4.53(0.05)	&	9.32(0.01)	\\	
	&	Soft	&	7.06(0.04)	&	9.41(0.01)	&	6.55(0.05)	&	9.34(0.01)	&	6.23(0.05)	&	9.19(0.02)	&	6.12(0.05)	&	9.02(0.03)	\\	
	&	SCAD	&	6.93(0.05)	&	9.41(0.01)	&	6.20(0.05)	&	9.34(0.01)	&	5.74(0.05)	&	9.23(0.02)	&	5.56(0.05)	&	9.11(0.03)	\\	
	&	Adap.lasso	&	6.00(0.05)	&	9.41(0.01)	&	5.32(0.06)	&	9.35(0.01)	&	4.89(0.06)	&	9.32(0.01)	&	4.74(0.06)	&	9.32(0.01)	\\	
\hline																		
\multirow{12}{*}{2}	&	 \multicolumn{9}{c}{Functional Frobenius norm}																\\		
	&	Hard	&	13.21(0.04)	&	17.03(0.01)	&	11.33(0.04)	&	16.40(0.01)	&	10.06(0.04)	&	16.06(0.01)	&	9.60(0.04)	&	15.96(0.01)	\\	
	&	Soft	&	13.54(0.04)	&	17.01(0.01)	&	12.06(0.04)	&	16.26(0.02)	&	11.10(0.04)	&	15.32(0.05)	&	10.75(0.04)	&	14.86(0.05)	\\	
	&	SCAD	&	13.50(0.04)	&	17.01(0.01)	&	11.90(0.04)	&	16.26(0.02)	&	10.78(0.04)	&	15.35(0.05)	&	10.36(0.04)	&	14.93(0.05)	\\	
	&	Adap. lasso	&	12.61(0.04)	&	17.01(0.01)	&	10.94(0.04)	&	16.39(0.01)	&	9.80(0.04)	&	15.99(0.02)	&	9.37(0.04)	&	15.81(0.03)	\\	
\cline{2-10}																			
	&	 \multicolumn{9}{c}{Functional matrix $\ell_1$ norm}																	\\	
	&	Hard	&	6.14(0.04)	&	7.27(0.01)	&	5.49(0.04)	&	7.19(0.01)	&	5.01(0.05)	&	7.16(0.01)	&	4.83(0.05)	&	7.15(0.01)	\\	
	&	Soft	&	6.22(0.02)	&	7.26(0.01)	&	5.90(0.03)	&	7.16(0.01)	&	5.65(0.03)	&	7.03(0.02)	&	5.55(0.03)	&	6.97(0.02)	\\	
	&	SCAD	&	6.21(0.02)	&	7.26(0.01)	&	5.87(0.03)	&	7.16(0.01)	&	5.58(0.03)	&	7.04(0.02)	&	5.45(0.03)	&	6.99(0.02)	\\	
	&	Adap. lasso	&	5.88(0.04)	&	7.26(0.01)	&	5.42(0.04)	&	7.19(0.01)	&	5.04(0.04)	&	7.15(0.01)	&	4.87(0.04)	&	7.14(0.01)	\\	
\hline

			\end{tabular}
		}	
	\end{center}
\end{table}

\begin{table}[!htbp]
	\caption{\label{TPR/FPR.table.kernel.100} 
	The average TPRs/ FPRs for partially observed functional scenarios and $p=100$  over 100 simulation runs.}
	\begin{center}
		\resizebox{6.2in}{!}{
			\begin{tabular}{cccccccccc}
\hline																					
	&			&	\multicolumn{2}{c}{$L_i = 11$}			&	\multicolumn{2}{c}{$L_i = 21$}			&	\multicolumn{2}{c}{$L_i=51$}			&	\multicolumn{2}{c}{$L_i=101$}			\\	
Model	&	Method	&	$\widecheck \bSigma_\A$	&	$\widecheck \bSigma_\U$	&	$\widecheck \bSigma_\A$	&	$\widecheck \bSigma_\U$	&	$\widecheck \bSigma_\A$	&	$\widecheck \bSigma_\U$	&	$\widecheck \bSigma_\A$	&	$\widecheck \bSigma_\U$	\\	\hline	
																					
\multirow{4}{*}{1}	&	Hard	&	0.57/0.00	&	0.00/0.00	&	0.62/0.00	&	0.00/0.00	&	0.65/0.00	&	0.00/0.00	&	0.66/0.00	&	0.00/0.00	\\		
	&	Soft	&	0.80/0.03	&	0.00/0.00	&	0.83/0.04	&	0.03/0.01	&	0.85/0.04	&	0.24/0.05	&	0.85/0.04	&	0.36/0.07	\\		
	&	SCAD	&	0.81/0.03	&	0.00/0.00	&	0.84/0.04	&	0.03/0.01	&	0.85/0.04	&	0.22/0.04	&	0.85/0.04	&	0.32/0.06	\\		
	&	Adap. lasso	&	0.67/0.00	&	0.00/0.00	&	0.71/0.00	&	0.00/0.00	&	0.73/0.00	&	0.00/0.00	&	0.74/0.00	&	0.01/0.00	\\	\hline	
\multirow{4}{*}{2}	&	Hard	&	0.48/0.00	&	0.00/0.00	&	0.57/0.00	&	0.00/0.00	&	0.65/0.00	&	0.00/0.00	&	0.68/0.00	&	0.00/0.00	\\		
	&	Soft	&	0.90/0.03	&	0.00/0.00	&	0.94/0.04	&	0.07/0.01	&	0.96/0.04	&	0.29/0.04	&	0.97/0.04	&	0.40/0.05	\\		
	&	SCAD	&	0.90/0.03	&	0.00/0.00	&	0.95/0.04	&	0.06/0.01	&	0.96/0.05	&	0.28/0.03	&	0.97/0.05	&	0.37/0.04	\\		
	&	Adap. lasso	&	0.70/0.00	&	0.00/0.00	&	0.78/0.00	&	0.00/0.00	&	0.83/0.00	&	0.02/0.00	&	0.85/0.00	&	0.03/0.00	\\	\hline	
											
			\end{tabular}
		}	
	\end{center}
\end{table}

\subsection{Real data}
\begin{figure}[!htbp]
	\centering
	\includegraphics[width=1\textwidth]{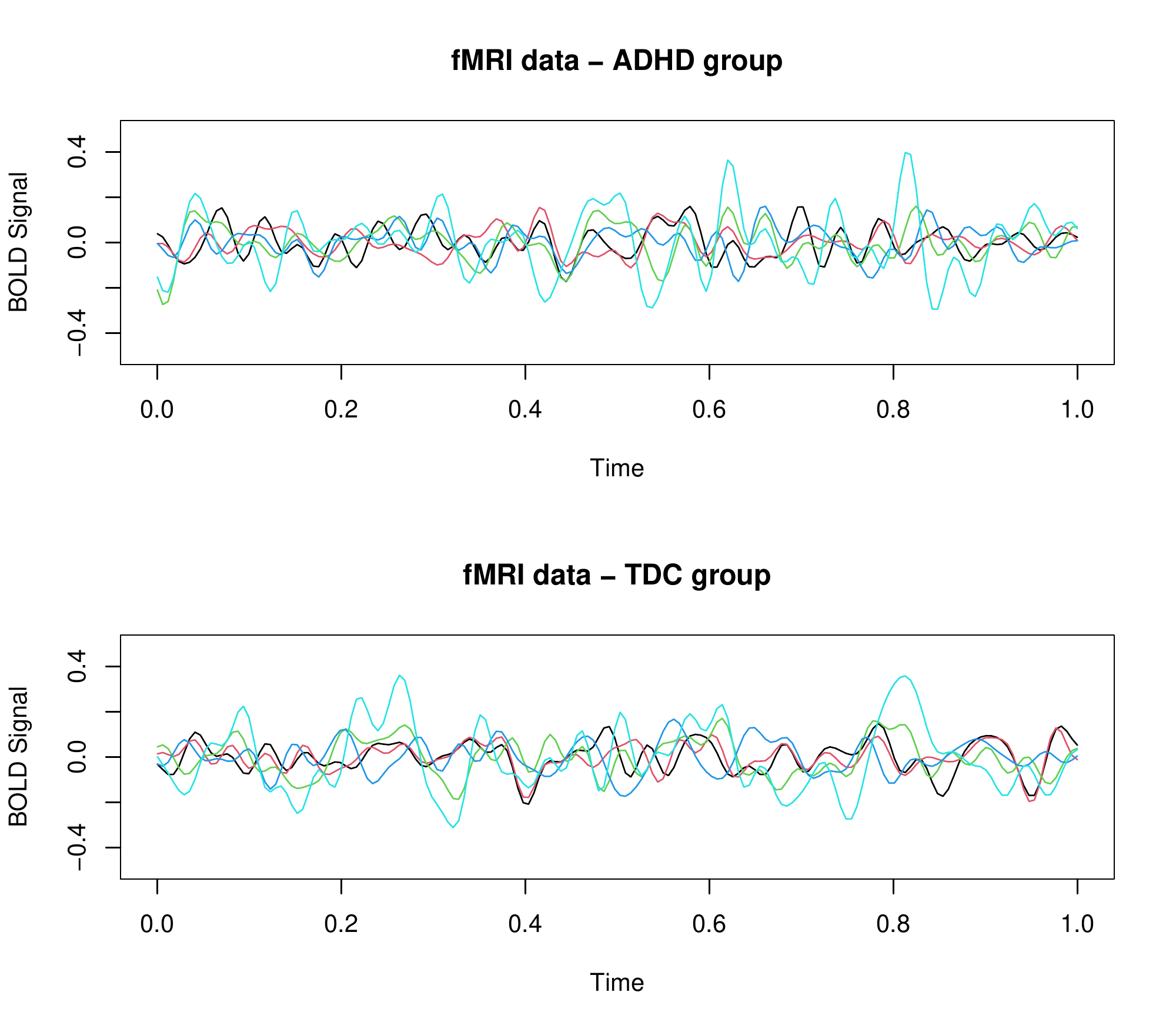}
	\caption{\label{data_ADHD}{
	ADHD dataset: the smoothed BOLD signals at the first $5$ ROIs of two subjects in ADHD and TDC groups respectively. The $5.73$-minute interval with $172$ scanning points is rescaled to $[0,1].$}}
\end{figure}

\begin{figure}[!htbp]
	\centering
	\includegraphics[width=1\textwidth]{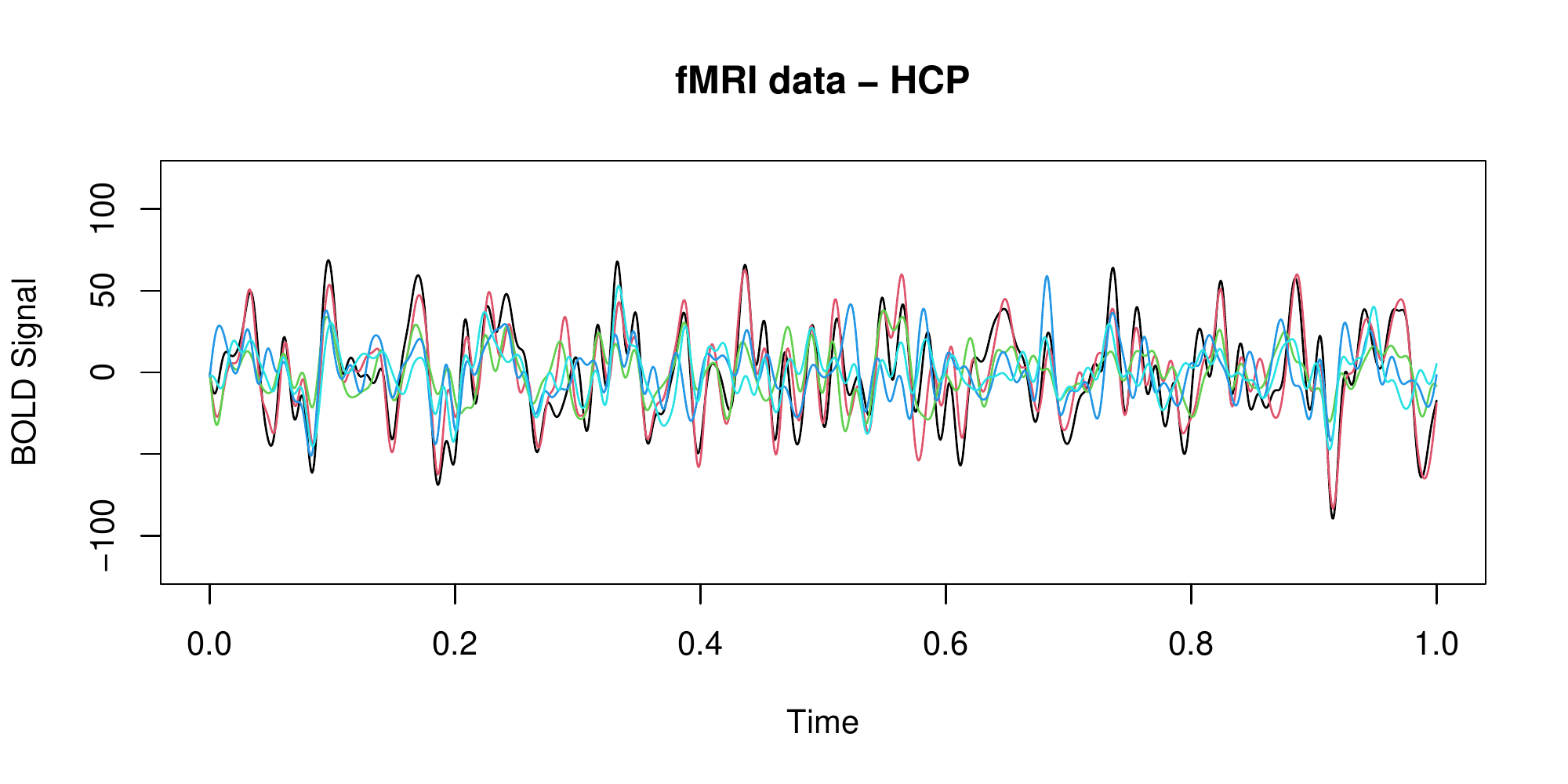}
	\caption{\label{data_HCP}{HCP dataset: the smoothed BOLD signals at the first $5$ ROIs of one subject. The $14.40$-minute interval with $1200$ scanning points ($14.40$ mins) is rescaled to $[0,1].$}}
\end{figure}

\begin{figure}[!htbp]
\centering
\begin{subfigure}{0.49\linewidth}
  \includegraphics[width=7cm]{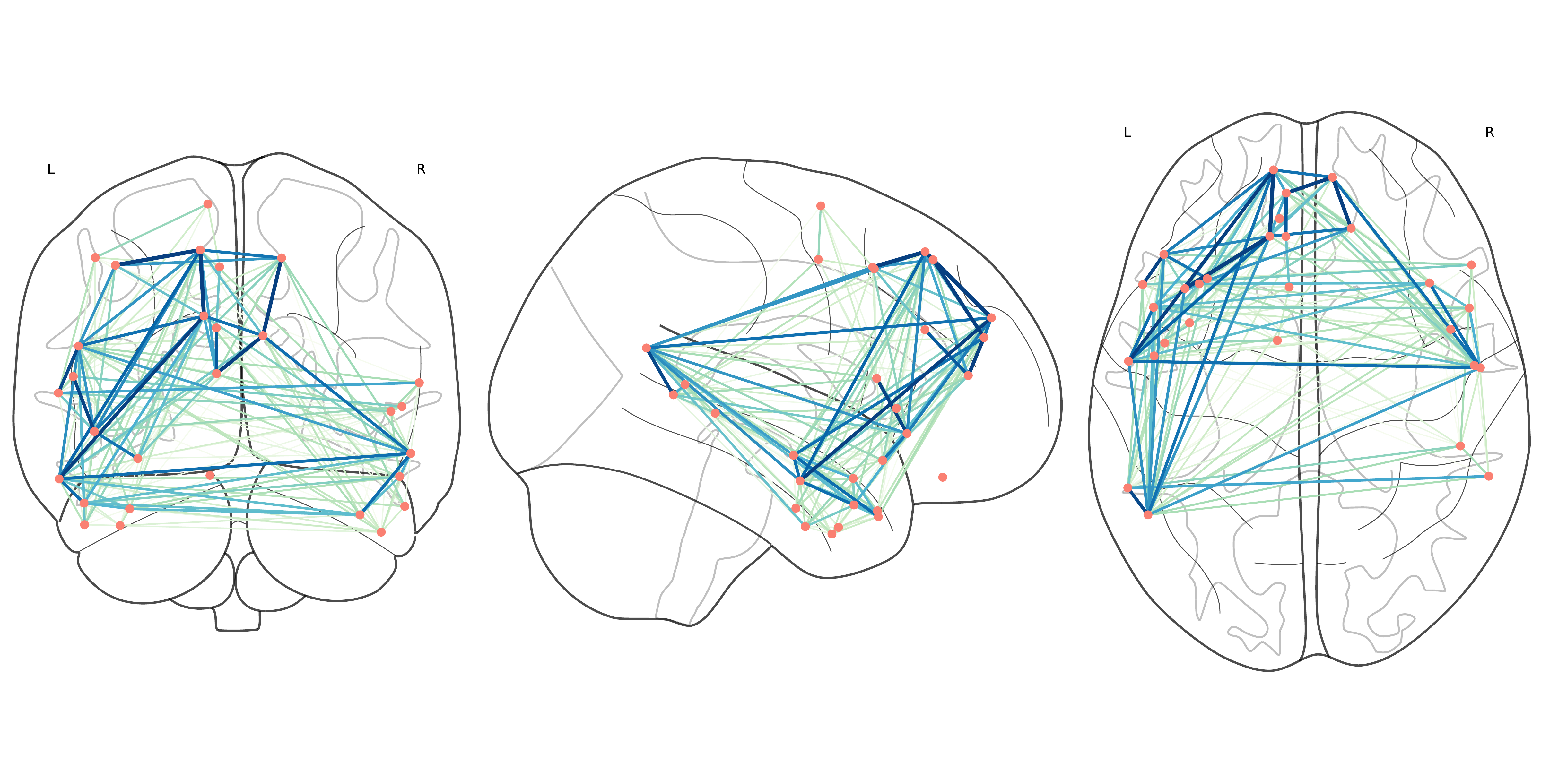}
  \caption{$\textit{gF}  \leq 8$: the medial frontal module in Fig.~\ref{hm_hcp}(a)} 
  \par\medskip 
  \includegraphics[width=7cm]{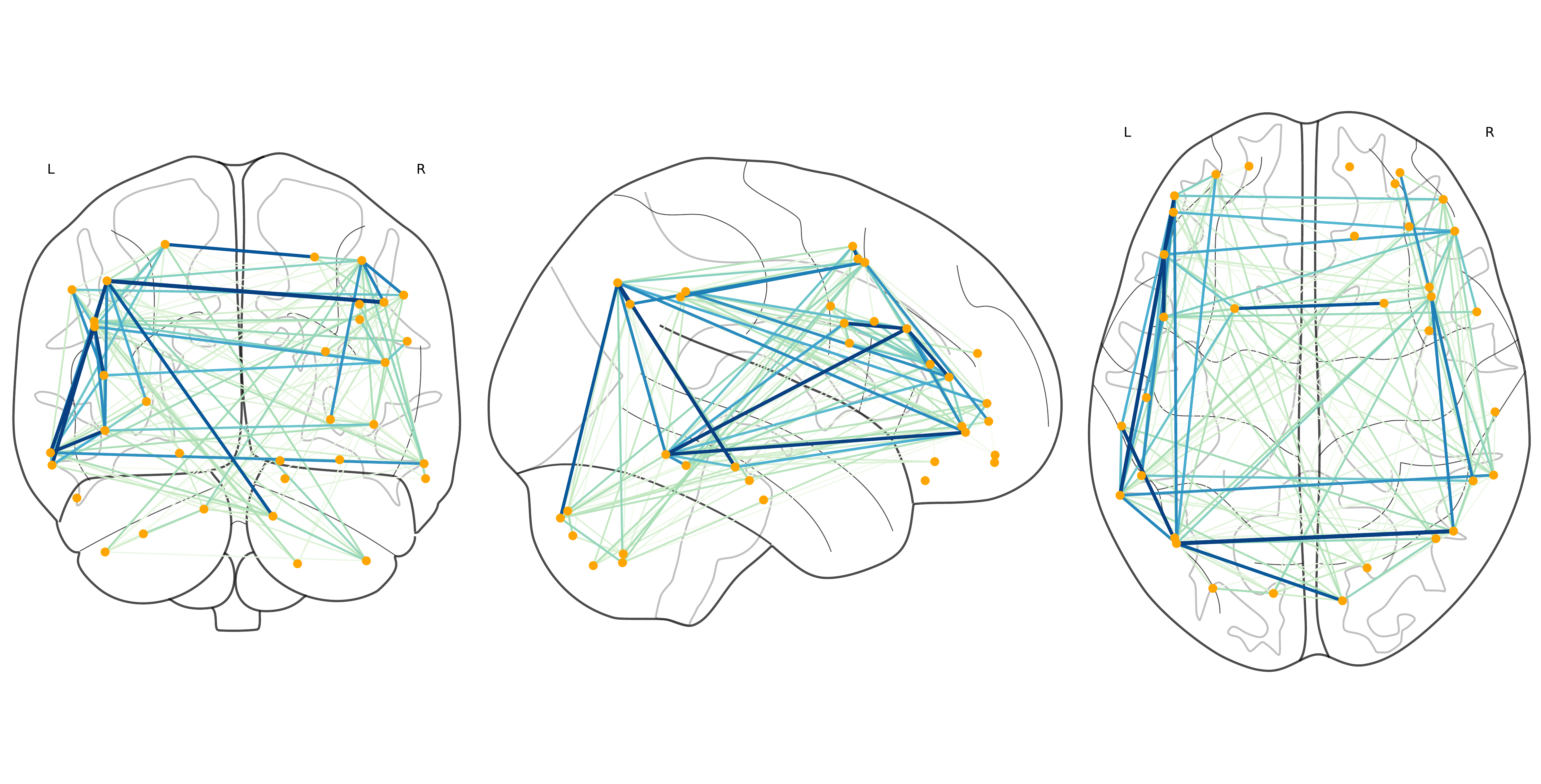}
  \caption{$\textit{gF}  \leq 8$: the frontoparietal module in Fig.~\ref{hm_hcp}(a)} 
  \par\medskip 
  \includegraphics[width=7cm]{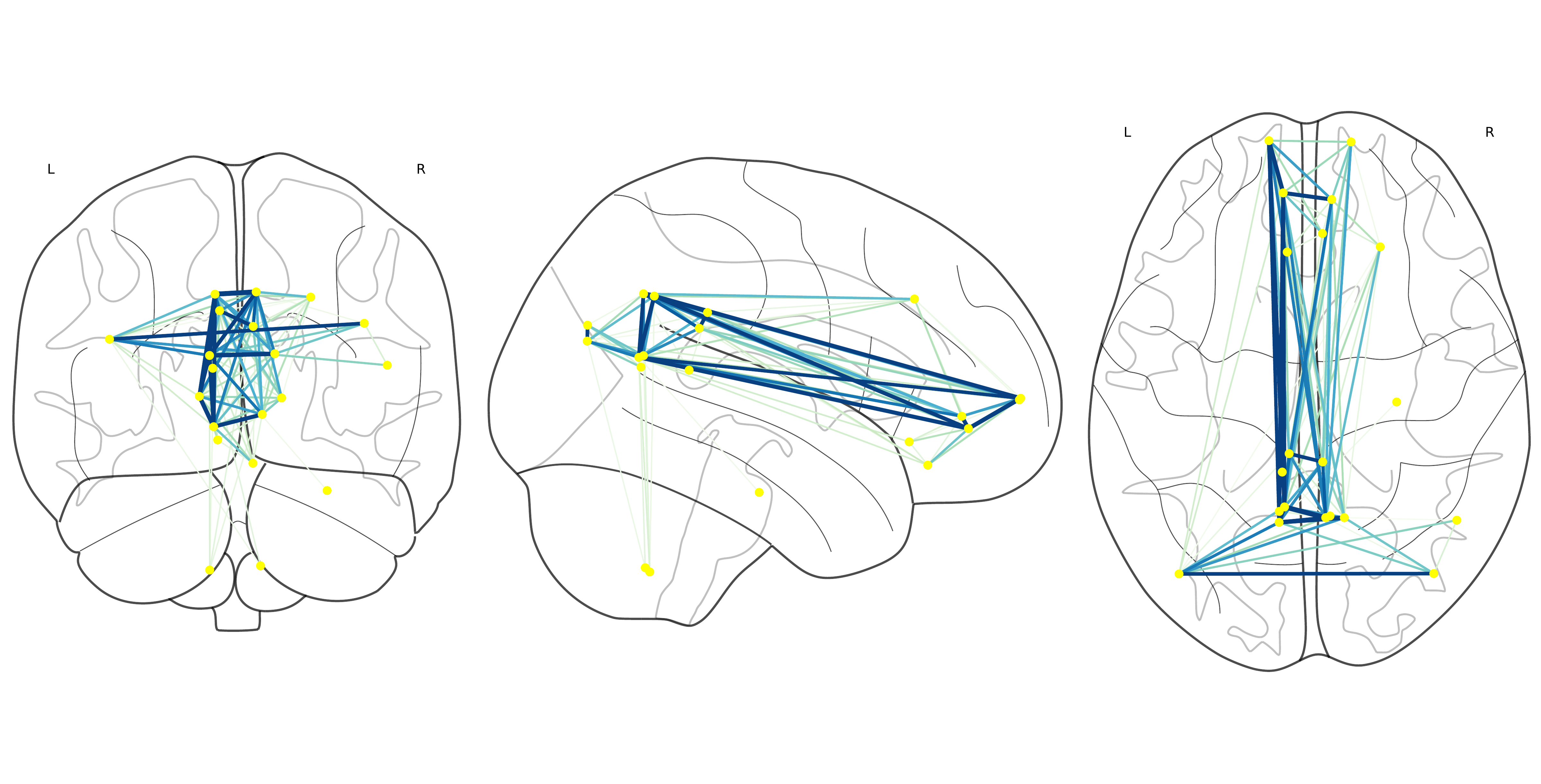}
  \caption{$\textit{gF}  \leq 8$: the default mode module in Fig.~\ref{hm_hcp}(a)} 
\end{subfigure}
\centering
\begin{subfigure}{0.49\linewidth}
  \includegraphics[width=7cm]{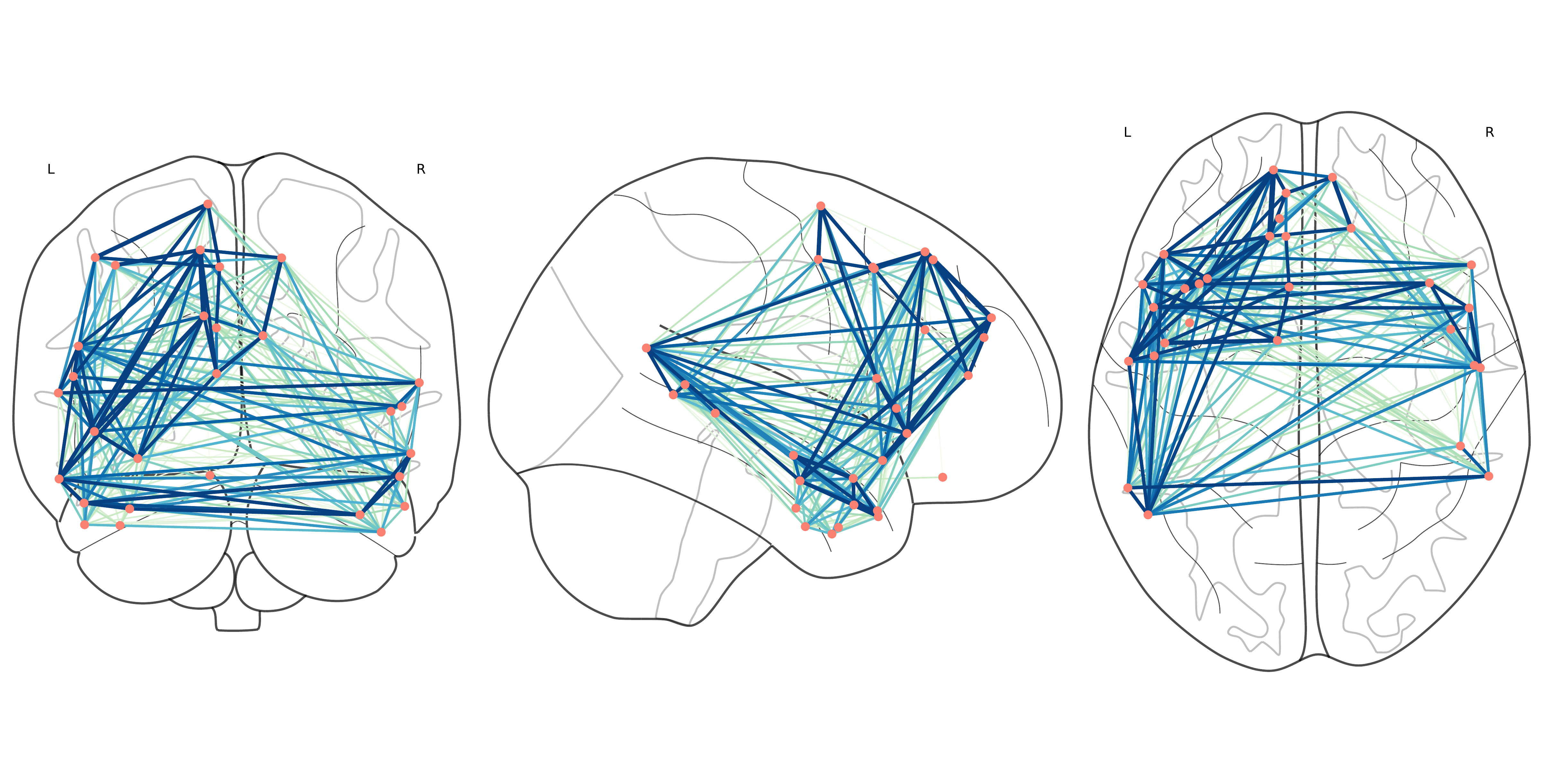}
  \caption{$\textit{gF}   \geq 23$: the medial frontal module in Fig.~\ref{hm_hcp}(b)} 
  \par\medskip 
  \includegraphics[width=7cm]{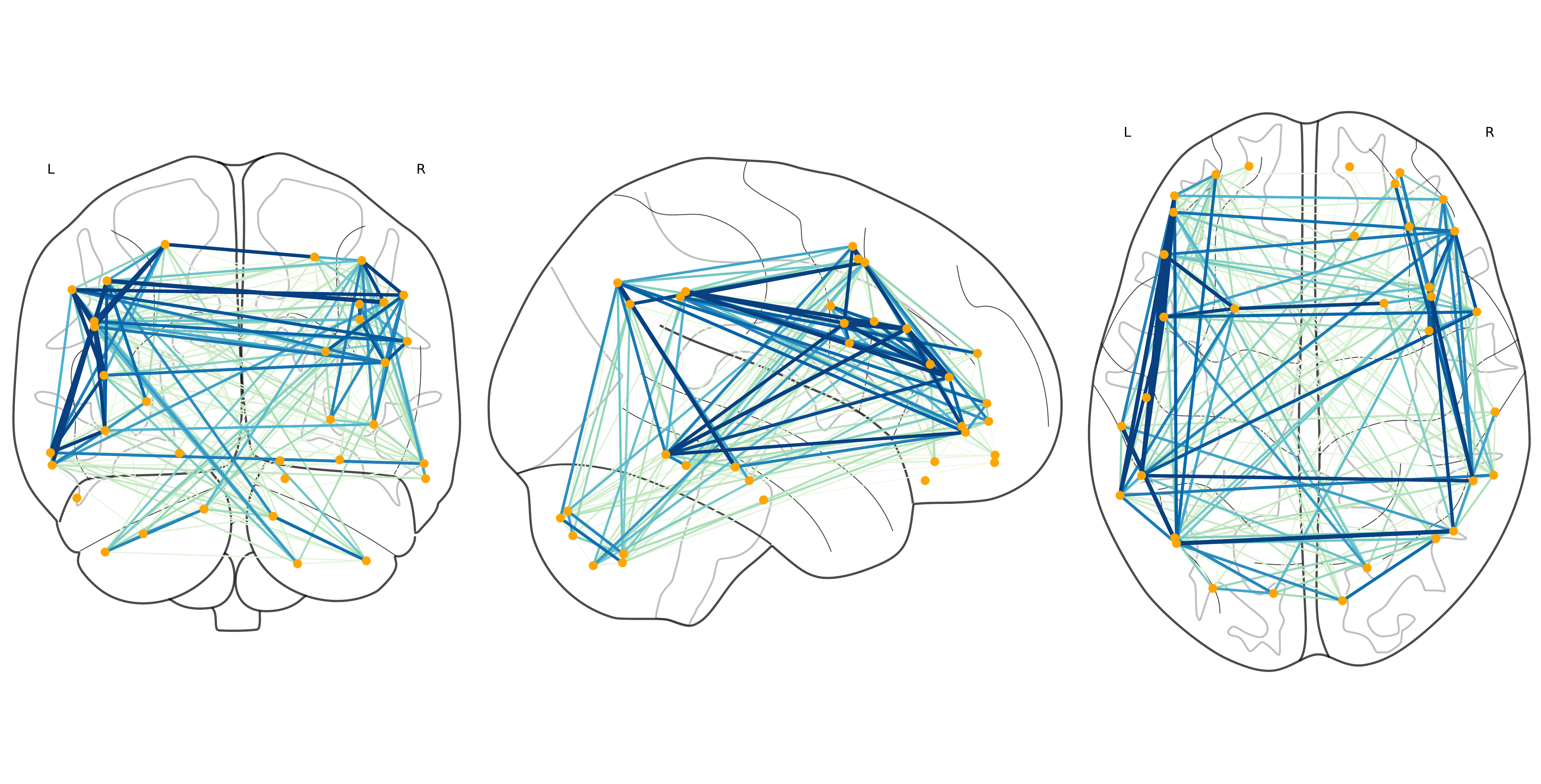}
  \caption{$\textit{gF}   \geq 23$: the frontoparietal module in Fig.~\ref{hm_hcp}(b)} 
  \par\medskip 
  \includegraphics[width=7cm]{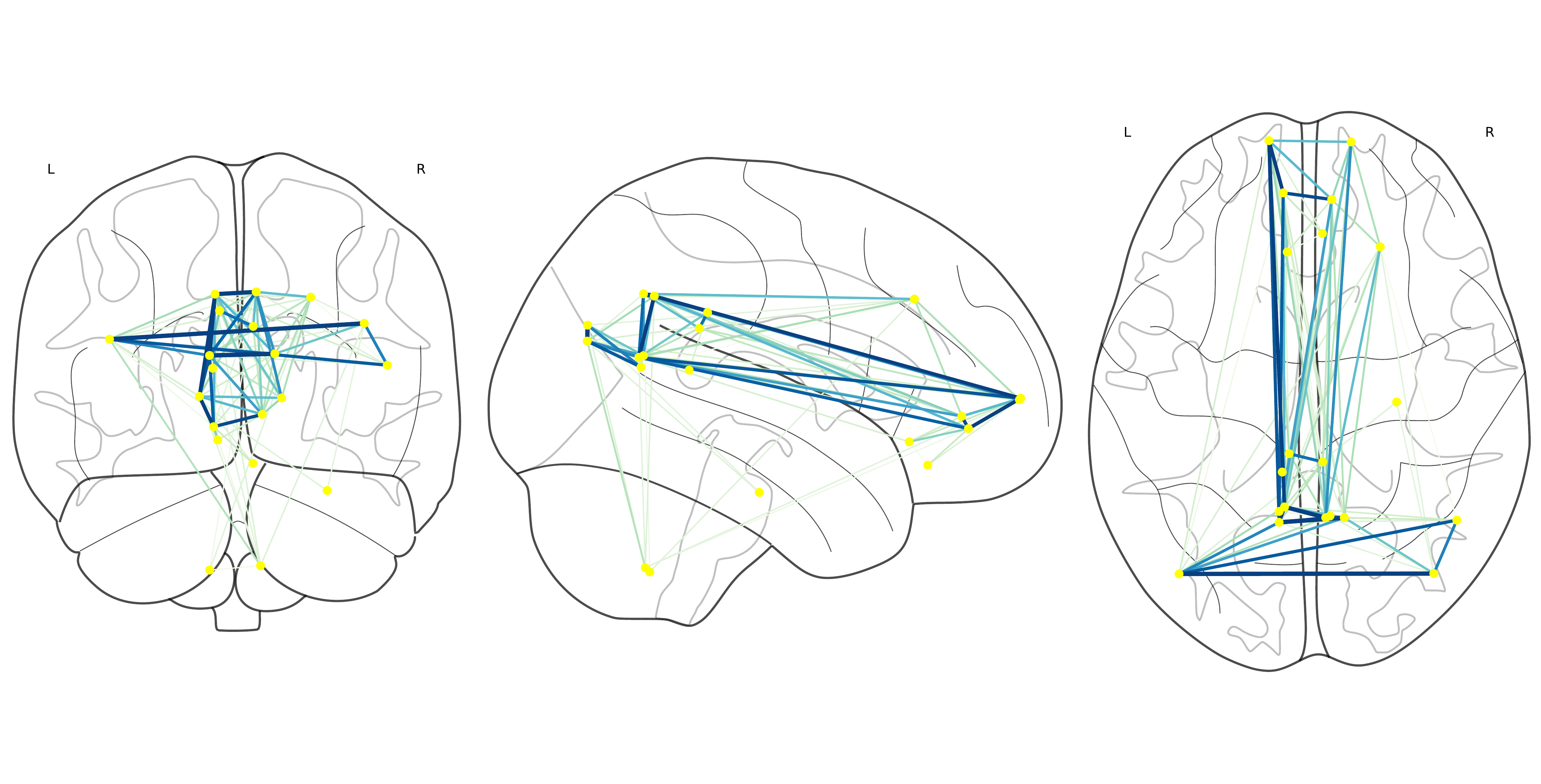}
  \caption{$\textit{gF}   \geq 23$: the default mode module in Fig.~\ref{hm_hcp}(b)} 
\end{subfigure}

\centering
\caption{\label{hcp_network_cv}{ The connectivity strengths in Fig.~\ref{hm_hcp}(a)--(b) at fluid intelligence $\textit{gF} \leq 8$ and $\textit{gF} \geq 23$.
Salmon, orange and yellow nodes represent the ROIs in the medial frontal, frontoparietal and  default mode modules, respectively. 
The edge color from cyan to blue corresponds to the value of $\|\widehat \Sigma_{jk}^\A\|_\cS/\{\|\widehat \Sigma_{jj}^\A\|_\cS\|\widehat \Sigma_{kk}^\A\|_\cS\}^{1/2}$ from small to large.}}
\end{figure}

\begin{figure}[!htbp]
\centering
\begin{subfigure}{0.49\linewidth}
  \includegraphics[width=7cm]{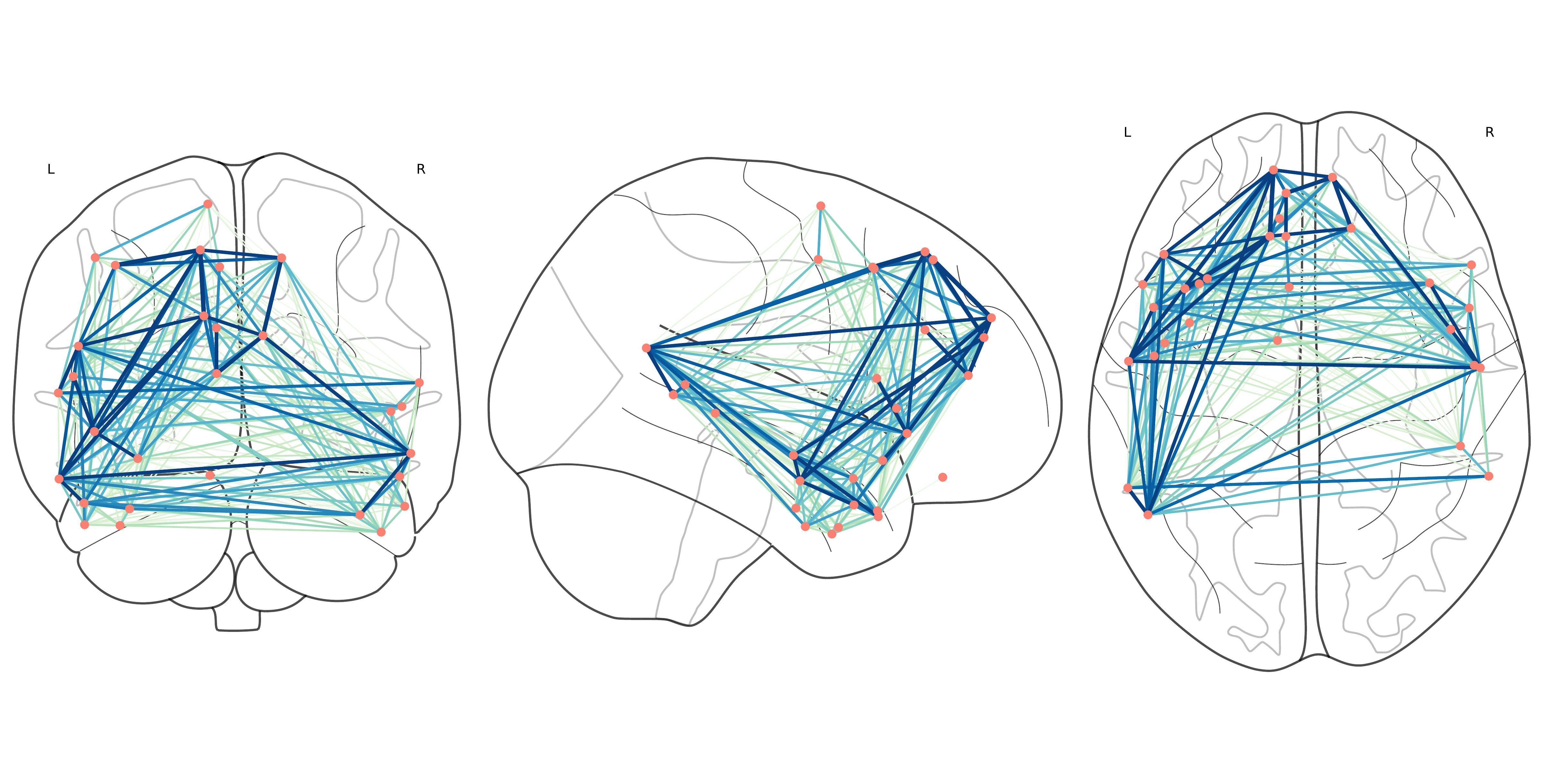}
  \caption{$\textit{gF}  \leq 8$: the medial frontal module in Fig.~\ref{hm_hcp}(c)} 
  \par\medskip 
  \includegraphics[width=7cm]{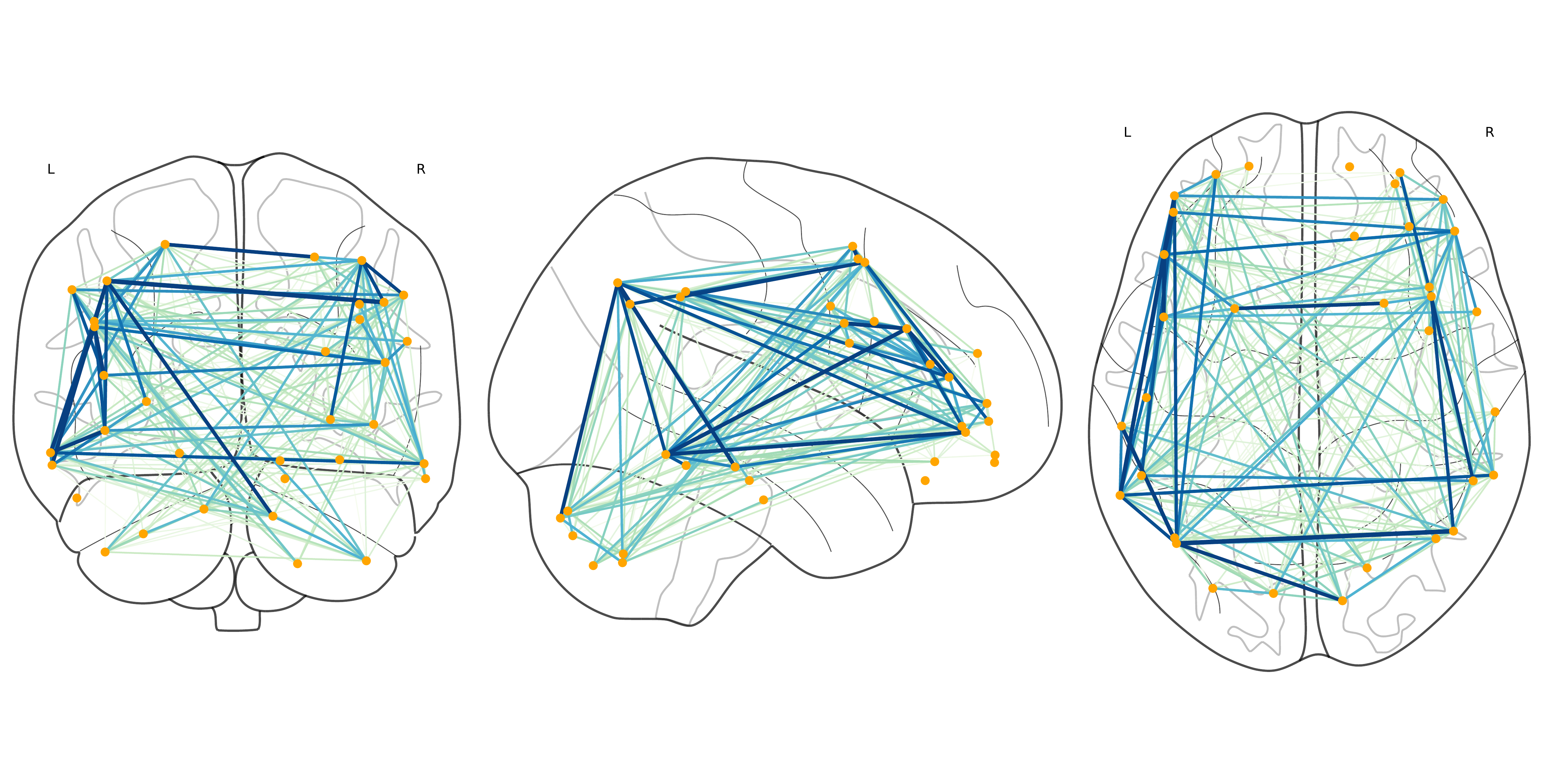}
  \caption{$\textit{gF}  \leq 8$: the frontoparietal module in Fig.~\ref{hm_hcp}(c)} 
  \par\medskip 
  \includegraphics[width=7cm]{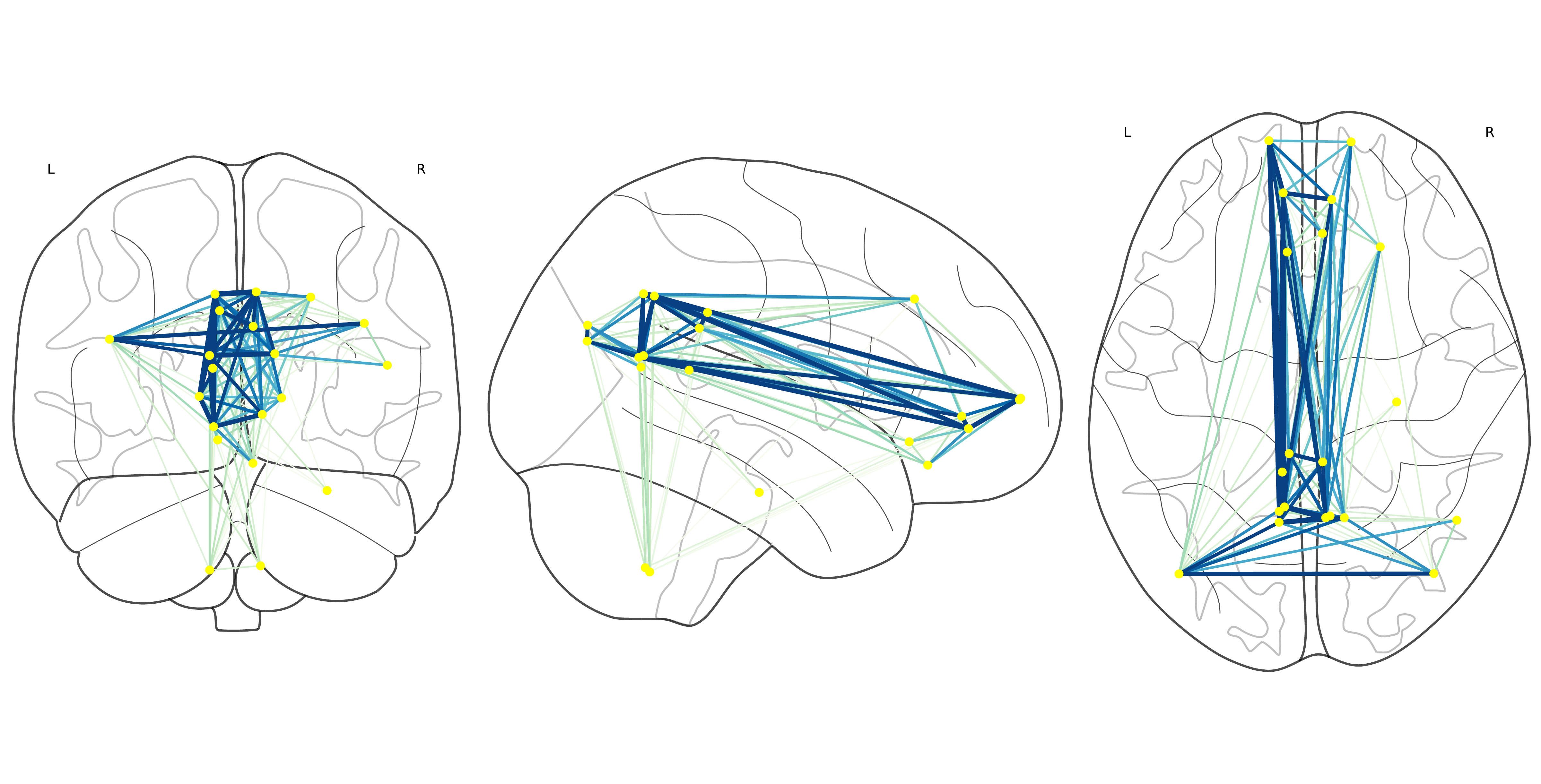}
  \caption{$\textit{gF}  \leq 8$: the default mode module in Fig.~\ref{hm_hcp}(c)} 
\end{subfigure}
\centering
\begin{subfigure}{0.49\linewidth}
  \includegraphics[width=7cm]{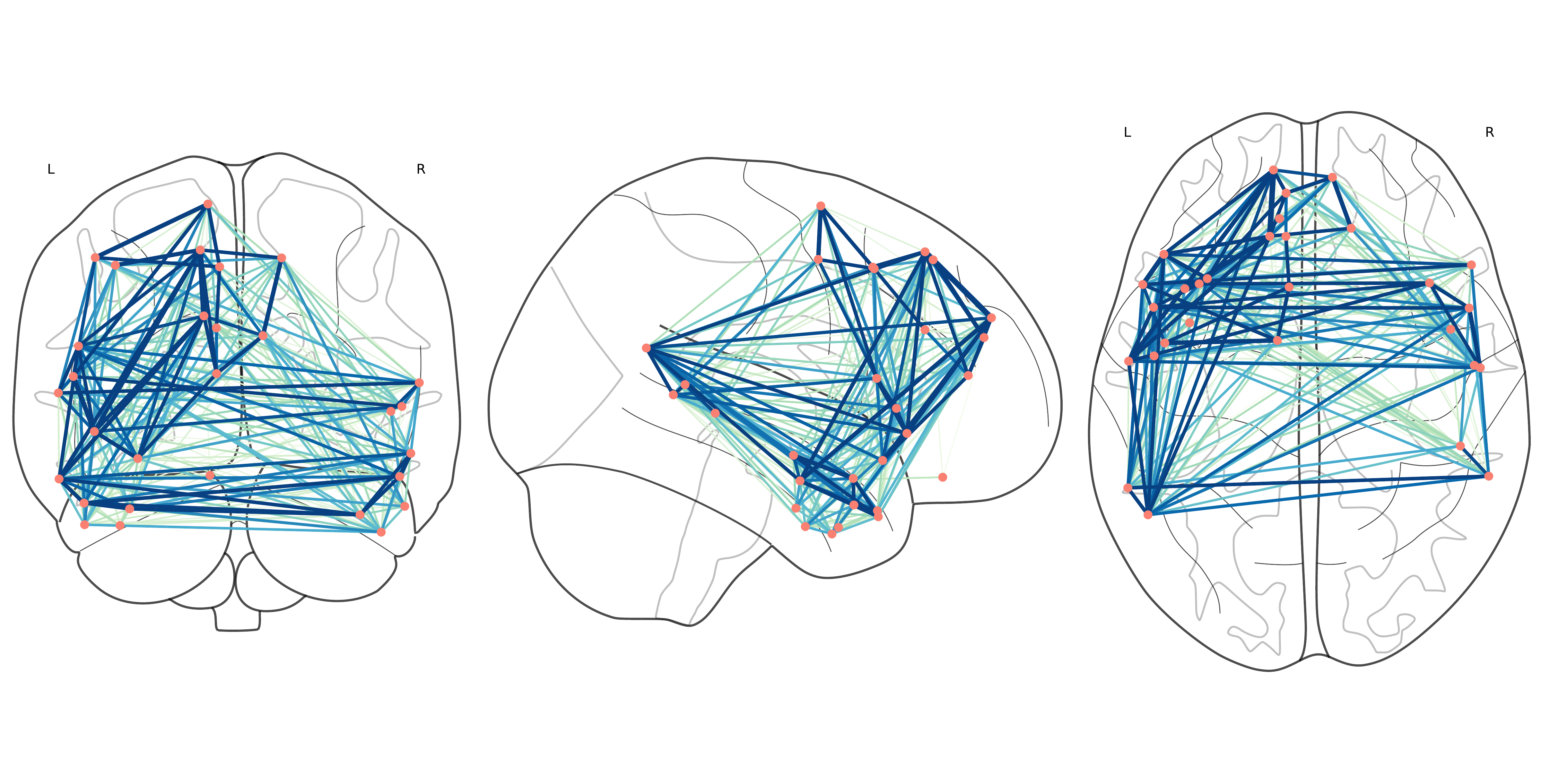}
  \caption{$\textit{gF}   \geq 23$: the medial frontal module in Fig.~\ref{hm_hcp}(d)} 
  \par\medskip 
  \includegraphics[width=7cm]{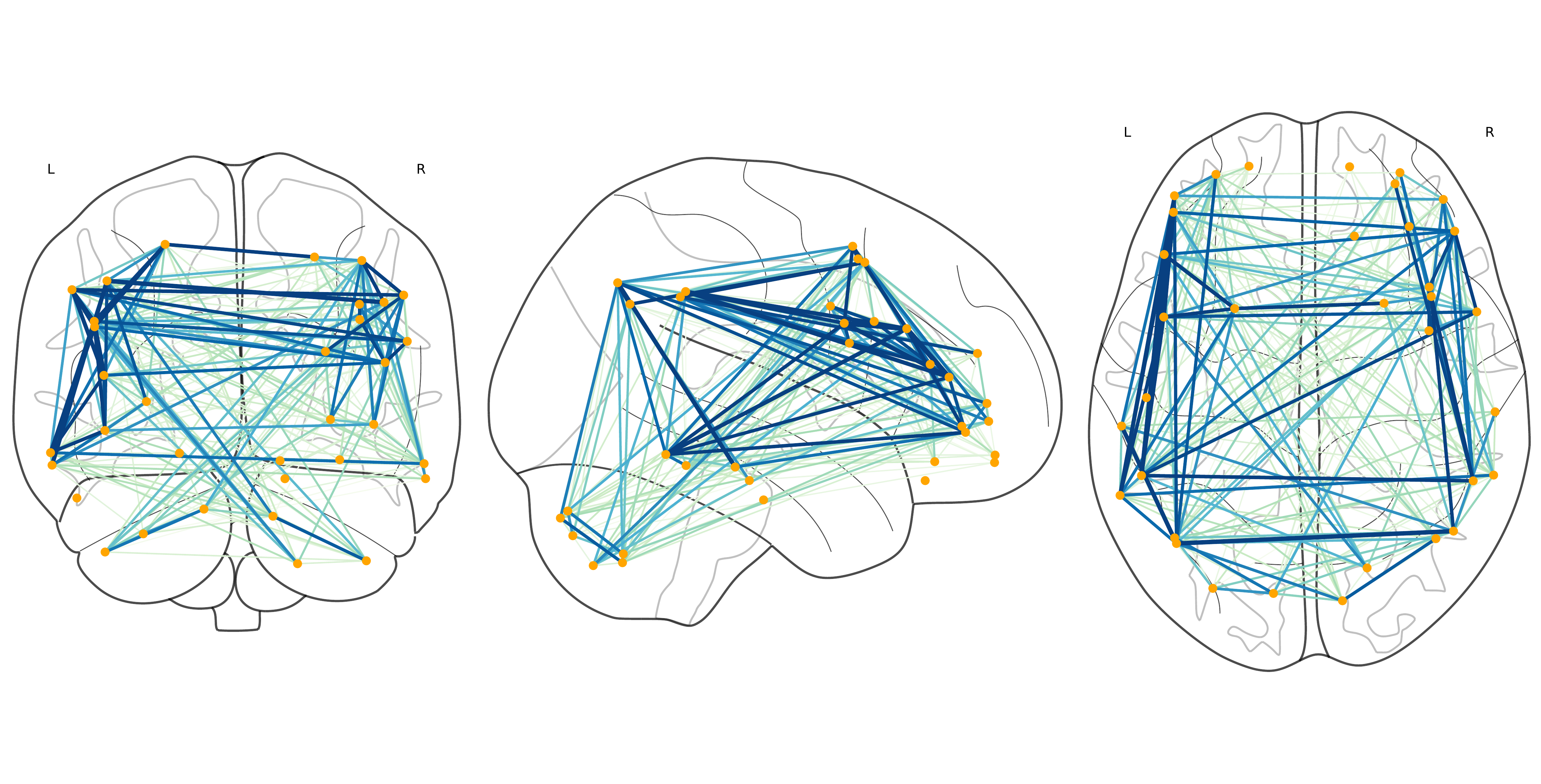}
  \caption{$\textit{gF}   \geq 23$: the frontoparietal module in Fig.~\ref{hm_hcp}(d)} 
  \par\medskip 
  \includegraphics[width=7cm]{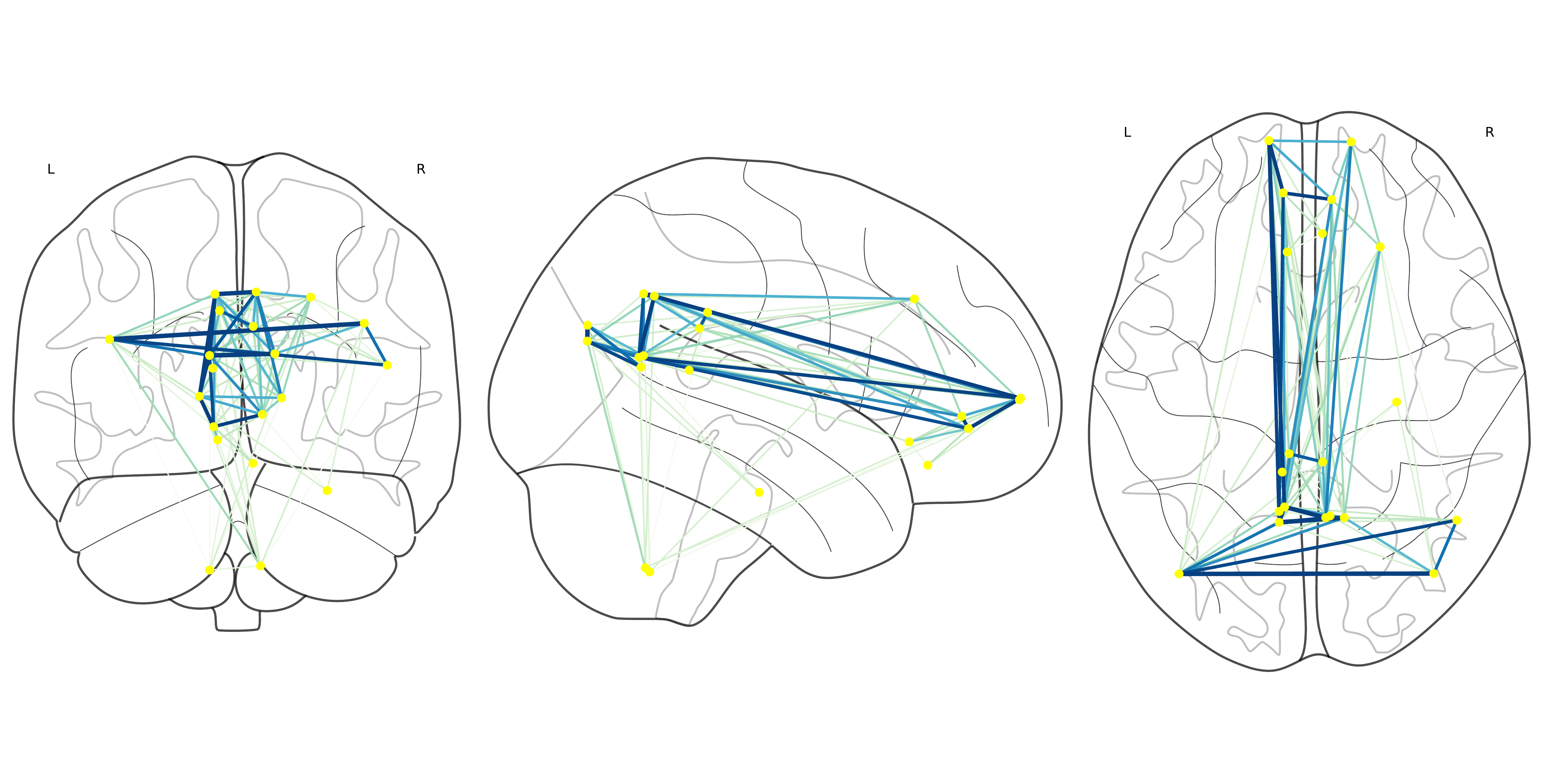}
  \caption{$\textit{gF}   \geq 23$: the default mode module in Fig.~\ref{hm_hcp}(d)} 
\end{subfigure}

\centering
\caption{\label{hcp_network30}{ The connectivity strengths in Fig.~\ref{hm_hcp}(c)--(d) at fluid intelligence $\textit{gF} \leq 8$ and $\textit{gF} \geq 23$.
Salmon, orange and yellow nodes represent the ROIs in the medial frontal, frontoparietal and  default mode modules, respectively. 
The edge color from cyan to blue corresponds to the value of $\|\widehat \Sigma_{jk}^\A\|_\cS/\{\|\widehat \Sigma_{jj}^\A\|_\cS\|\widehat \Sigma_{kk}^\A\|_\cS\}^{1/2}$ from small to large.}}
\end{figure}

Figures~\ref{data_ADHD} and \ref{data_HCP} display the pre-smoothed BOLD signal trajectories at a selection of ROIs of subjects from the ADHD and HCP datesets, respectively. 
Figures \ref{hcp_network_cv} and \ref{hcp_network30} plot the connectivity strengths at fluid intelligence $\textit{gF} \leq 8$ and $\textit{gF} \geq 23$ in Fig.~\ref{hm_hcp}(a)--(b) and Fig.~\ref{hm_hcp}(c)--(d), respectively. We observe that as \textit{gF} increases, the connectivity strengths in the  medial frontal and frontoparietal modules tend to increase while those in the default mode module decrease, which is consistent with our finding in Section~\ref{sec.hcp}.

\newpage

\section*{References}
\begin{description}
 	\item
 	Bosq, D. (2000). {\it Linear Process in Function Spaces}. New York: Springer.   
 	 \item
 	Boucheron, S., Lugosi, G. and Massart, P. (2014). {\it Concentration Inequalities: A Nonasymptotic Theory of Independence}. Oxford University Press.
 	\item 
     Kosorok,  M.  R. (2008). {\it Introduction to empirical processes and semiparametric inference.}   Springer Series in Statistics. Springer, New York.
    
    \item 
    Qiao, X., Qian, C., James, G. M. and Guo, S. (2020). Doubly functional graphical models in high dimensions. {\it Biometrika}, {\bf 107}, 415–431.
    
    \item 
    Rothman, A. J., Levina, E. and Zhu, J. (2009). Generalized thresholding of large covariance matrices. {\it Journal of the American Statistical Association}, {\bf 104}, 177–186.
\end{description}

\end{document}